\documentclass[reqno]{amsart}
\input{FancyMath.sty}
\usepackage{newunicodechar}
\usepackage[utf8]{inputenc}

\title{Twisting and foam inequalities for the upsilon invariant}

\makeatletter
\@namedef{subjclassname@1991}{2020 Mathematics Subject Classification}
\makeatother

\author{William Ballinger}
\address{Department of Mathematics, Harvard University}
\email{ballinger@harvard.edu}

\author{Gary D.~Dunkerley}
\address{Department of Mathematics, University of Georgia}
\email{garydunkerley@uga.edu}

\date{\today}


\begin{document}

\begin{abstract}
We describe two families of bounds on the Ozsv\'ath--Stipsicz--Szab\'o invariant $\Upsilon_K(t)$, one controlling its behavior under inserting full twists into a knot and one controlling a limited class of foams, or singular cobordisms between knots. 
By constructing certain explicit foams, we show that these two formal properties of $\Upsilon$ together exactly determine $\Upsilon_K(t)$ for all torus knots. This generalizes previously observed behavior at $t = 0$, Livingston's slice-torus invariants, and at $t=1$, Sato's unoriented analogue. 
Our foam bounds on $\Upsilon$ make significant use of Alishahi's generalization of $\Upsilon$ to knotted graphs. Along the way we compute Alishahi's invariant for all theta-graphs lying on the Heegaard torus in $S^3$. We also relate the evaluations of $\Upsilon$ for torus theta-graphs to a distinguished family of properly embedded 2-complexes in the 4-ball sharing many properties with Batson's surfaces.
\end{abstract}

\maketitle

\tableofcontents

\clearpage

\section{Introduction}

Knot Floer homology, defined by Ozsv\'{a}th--Szab\'{o} and independently by Rasmussen, is a powerful invariant of knots in closed, orientable 3-manifolds.
In particular, the knot Floer homology of a knot $K$ in $S^{3}$ can be used to obtain invariants which constrain the complexity of slice surfaces of ${ K }$--- that is, properly embedded surfaces ${ \Sigma }$ in the 4-ball $B^4$ such that $\Sigma \cap \partial B^{4} = K$. 
Where $\mathscr{C}$ denotes the smooth concordance group of knots in $S^3$, 
Ozsv\'{a}th--Stipsicz--Szab\'{o} \cite{OSS17a} use knot Floer homology to define a homomorphism
$$
\Upsilon : \mathscr{C} \rightarrow \mathrm{PL}[0,2], \qquad 
[K] \mapsto \Upsilon_{K} 
$$
where $[K]$ denotes the smooth concordance class of a knot $K\subset S^3$, and $\mathrm{PL}[0,2]$ the group of continuous piecewise linear functions over the interval $[0,2]$.
As a concordance homomorphism $\Upsilon$ appears to be quite powerful, particularly for studying linear independence among elements of $\mathscr{C}$
\cite{OSS17a, Wa18,KP18,FPR19}. 
However for a given knot $K$ it is somewhat unclear if evaluations of $\Upsilon_K$ communicate geometric information --- e.g. does the fact that $\Upsilon_{T(5,7)}(\tfrac{3}{4}) = -\tfrac{29}{4}$  communicate anything interesting or particular about $T(5,7)$?
It is known that $\Upsilon_K$ imposes constraints on smooth slice surfaces bounded by $K$; Ozsv\'{a}th--Stipsicz--Szab\'{o} showed \cite{OSS17a} that for any value of $t$ there is the genus bound  
\begin{equation}\label{eqn:genusBound}
t\cdot g_4(K) \ge \left|\Upsilon_K(t)\right|
\end{equation}
where $g_4$ denotes the orientable smooth slice genus.
Moreover the slope of $\Upsilon_K$ at $t=0$ is equal to the negative of the $\tau$ invariant, which is a \textit{slice-torus invariant} in the sense of \cite{Liv04}. In later work Ozsv\'{a}th--Stipsicz--Szab\'{o} proved that $\upsilon(K):=\Upsilon_K(1)$ constrains the complexity of nonorientable slice surfaces for $K$ as measured by the Euler characteristic and normal Euler number \cite{OSS17b}.

When $K$ is the $(p,q)$-torus knot $T(p,q)$, one can excavate further geometric meaning from outputs of $\Upsilon_K$.
For $t$ near $0$, $\Upsilon_{T(p,q)}(t)$ is as small as it could be without violating the genus bound (\ref{eqn:genusBound})
while $\Upsilon_{T(p,q)}(1)$ is as large as it could be without violating a bound associated to a distinguished, nonorientable surface described by Batson \cite{Bat14, JVC21}. In fact, $\Upsilon_{T(p,q)}(1)$ is completely determined by topological information: $\Upsilon_{T(p,q)}(t)$ is also controlled by surfaces in $\#^k \overline{\mathbb{C} P^2}$ (\cite{Zem19a}, Theorem 1.1), and for  $t = 1$ the upper bound from Batson's surface and the lower bound from a surface in $\#^k \overline{\mathbb{C} P^2}$ exactly match. 

In this note, we study how the behavior described in the previous paragraph generalizes to other values of $t$.
For $t \neq 1$, the r\^{o}le played by nonorientable surfaces in controlling $\Upsilon_K(t)$ will instead be filled by foams, cobordisms with certain allowed singularities. As we move Morse theoretically through a foam, a knot will be replaced by a knotted graph, and so we will rely on a generalization of $\Upsilon$ defined by Alishahi for such objects \cite{Ali22}. Given an embedding of a bipartite and balanced graph $G : \Gamma \hookrightarrow Y^3$, Alishahi defines a function $\Upsilon_G$ on a certain polytope of labelings of edges of $G$ by real numbers. Alternatively, Alishahi's $\Upsilon$ assigns a single real number\footnote{Things are this clean only when $G$ is a $\Theta_n$-graph; Alishahi's $\Upsilon$ is a multivalued invariant in general.} to a spatial graph whose edges are $\mathbb{R}_{\geq 0}$-colored by colors which sum to $2$ around each vertex.

\begin{figure}[htbp]
    \centering

    \begin{subfigure}[t]{0.25\textwidth}
        \centering
        \includegraphics[width=\textwidth]{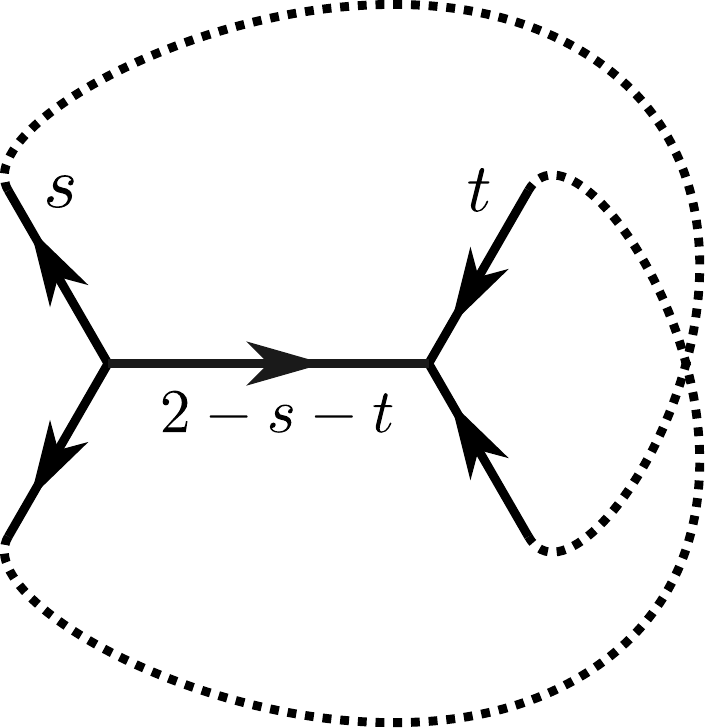}
        \caption{$G$}
        \label{fig:Gamma}
    \end{subfigure}
    \hfill
    \begin{subfigure}[t]{0.25\textwidth}
        \centering
        \includegraphics[width=\textwidth]{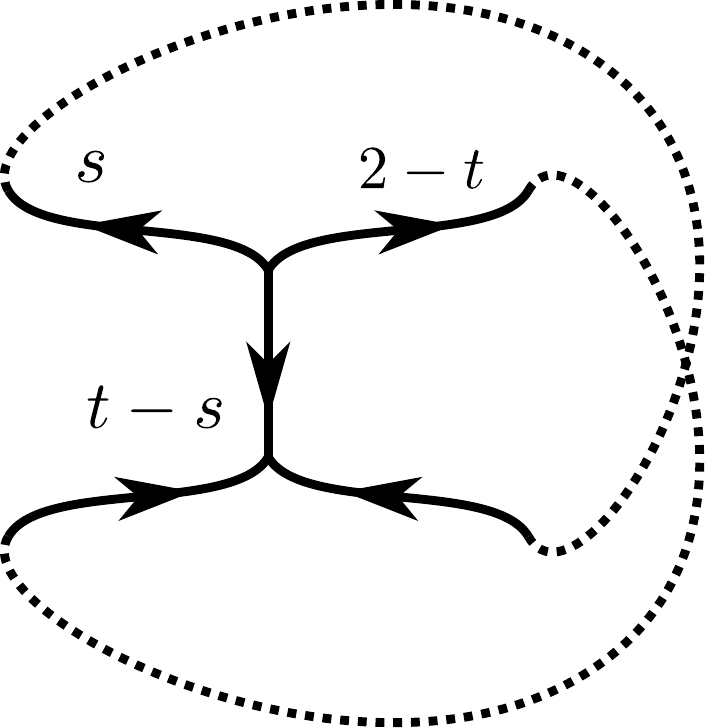}
        \caption{$G'$}
        \label{fig:Gamma_prime}
    \end{subfigure}
    \hfill
    \begin{subfigure}[t]{0.28\textwidth}
        \centering
        \includegraphics[width=\textwidth]{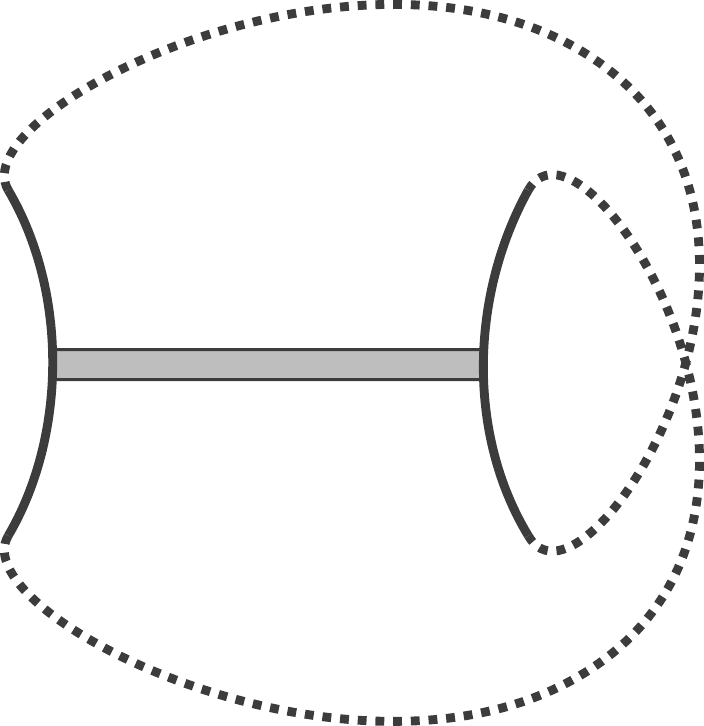}
        \caption{The associated unoriented band}
        \label{fig:linking_curves}
    \end{subfigure}

    \caption{}
    \label{fig:theorem_1.1_figures}
\end{figure}

\begin{theorem}\label{thm:FoamBound1}
    Letting $0\leq s \leq t \leq 2$ and $s+t \leq 2$, 
    suppose that $G$ and $G'$ are knotted, colored $\Theta_3$-graphs differing as in Figures \ref{fig:Gamma} and \ref{fig:Gamma_prime}. Consider the knots $K,K'$ formed by removing the central short edge in each of $G,G'$, and let $b$ be the single-band unoriented cobordism connecting $K$ and $K'$ (Figure~\ref{fig:linking_curves}). Then 
    \begin{equation*}
        \left|\Upsilon_{G'} - \Upsilon_{G} - \frac{s}{4} e(b)\right|\leq \frac{s}{2}.
    \end{equation*}
\end{theorem}

This is our first analogue of the Ozsv\'ath--Stipsicz--Szab\'o unoriented surface bound for general $t$; we will discuss generalizations in \S~\ref{sec:FoamPrelims}. 
For any of the invariants of colored $\Theta_3$-graphs we consider, removing a 0-colored edge does not change the invariant. When $s = t = 1$ the horizontal and vertical edges in Figures \ref{fig:Gamma} \& \ref{fig:Gamma_prime} respectively will be colored zero. Since removing a zero-colored edge doesn't change $\Upsilon$, the graphs reduce to knots, $\Upsilon_G = \upsilon(K)$, $\Upsilon_{G'} = \upsilon(K'),$ and this reduces to precisely the Ozsv\'ath--Stipsicz--Szab\'o bound.

The bound in Theorem \ref{thm:FoamBound1} can be interpreted as associating a change in the $\Upsilon$ invariant of $\Theta_3$-graphs which differ by a \textit{tetrahedron move}; see Definition \ref{defn:tetrahedronMove}.
In \S~\ref{sec:torusThetaWebs} we apply this to the study of torus knots. Given a torus knot $T(p,q)$, we construct a distinguished foam --- the \textit{Batson foam} --- embedded in $B^4$ with boundary $T(p,q)$. With respect to the radial Morse function for $B^4$, the regular level sets of this foam are $\Theta_3$-graphs embedded in a standard Heegaard torus of $S^3$, i.e. each level set is a \textit{torus $\Theta_3$-graph} and with respect to the radial Morse function on $B^4$, adjacent level sets differ by a tetrahedral move as in Theorem~\ref{thm:FoamBound1}.

\begin{theorem}\label{thm:foamBoundTight}
Any torus knot, or any torus $\Theta_n$-graph, saturates the bound on $\Upsilon$ coming from Theorem~\ref{thm:FoamBound1} and the Batson foam.
\end{theorem}

\begin{remark}
    In \S~\ref{sec:torusThetaGraphHomology}, we give a description of the tangle Floer homology of the torus $\Theta_n$-graphs in terms of subsemigroups of $\Z^2$. However, it is not straightforward to read off the formulae for $\Upsilon$ from Theorems~\ref{thm:foamBoundTight} or~\ref{thm:twistBoundTight} from this homology. Instead, we obtain an inductive relationship between the $\Upsilon$ invariants of torus $\Theta_n$-graphs which is mediated by the Farey tree.
\end{remark}

As a final ingredient, we generalize Zemke's proof that $\Upsilon$ controls cobordisms in $\#^k \overline{\mathbb{C} P^2}$ between embedded graphs in the 3-sphere, obtaining
\begin{theorem}\label{thm:twistBound1}
Let $\Theta_n$ denote the loopless, oriented graph with two vertices and $n$ edges oriented so that one vertex is a source and the other is a sink. 
Let $G, G' : \Theta_n \hookrightarrow S^3$ be such that $G'$ is obtained from $G$ by a positive full twist in a disk that meets the $i$th edge algebraically $a_i$ times. 
Equipping $\Theta_n$ with an $\mathbb{R}_{\geq 0}$-coloring of its edges satisfying $\sum_i t_i = 2$, we have
$$-\min\limits_{k \in \Z} \left( \binom{a_1+k}{2} t_1 + \cdots + \binom{a_n + k}{2}t_n \right) \leq \Upsilon_{G'} - \Upsilon_G. $$
\end{theorem}
Again, this bound is always tight:
\begin{theorem}\label{thm:twistBoundTight}
    Any positive torus knot or torus $\Theta_n$-graph may be untied by a series of full twists following the Euclidean algorithm; its $\Upsilon$ invariant saturates the bound provided by Theorem~\ref{thm:twistBound1}.
    \newline 
    In particular, this gives a recursive formula for the value of $\Upsilon$ on torus knots and torus $\Theta_3$-graphs, generalizing a recursive formula for the $\Upsilon$ invariant of torus knots described by Feller and Krcatovich \cite{FK17}.
\end{theorem}

Note that we now must have matching upper and lower bounds on any $\Upsilon_{T(p,q)}(t)$. This is the sense we mean when saying these invariants are determined by topological information. On its own, this has no computational value --- there are much easier ways to calculate $\Upsilon_{T(p,q)}(t)$.
On the other hand, these techniques likely extend somewhat beyond the case of the $\Upsilon$ invariant of torus knots. In future work, we plan to show that certain invariants derived from colored $\mathfrak{sl}_N$ homology satisfy analogues of Theorems~\ref{thm:FoamBound1} and~\ref{thm:twistBound1}, allowing these invariants to be calculated in some situations without computing the (much more intricate) underlying knot homologies.

\section{Acknowledgments}

GDD thanks his advisor, Akram Alishahi, for her support and guidance. 
We thank the organizers of the 2024 \textit{New Structures in Low-Dimensional Topology} conference where the collaboration which produced this work began. WB was partially supported by the Simons Collaboration grant on New Structures in Low-Dimensional Topology. GDD was partially supported by the National Science Foundation under grants DMS-2238103 and DMS-2342252. 
\vspace{5pt}

\noindent \textbf{AI Acknowledgment}  
The authors acknowledge use of large language models for assistance in formatting some figures and for assistance with line edits, spelling, and grammar.

\section{Heegaard Floer homology of spatial graphs}\label{sec:tangleFloerHomology}

Tangle Floer homology, defined by Alishahi and Eftekhary \cite{AE15, AE20}, is an invariant of systems of disjoint, properly embedded arcs $T$ in a compact, oriented 3-manifold with boundary $M$. It can be thought of as a ``minus'' version of the sutured Floer homology of \cite{Juh06}.
When the boundary is a disjoint union of spheres, one can interpret $(M,T)$ as describing a spatial graph and its tangle Floer homology is a straightforward generalization of knot Floer homology of \cite{Ras03} and \cite{OS04b}. By way of analogy to how the knot Floer homology of a knot $K \subset Y^3$ can be identified with the sutured Floer homology of its meridionally sutured complement, one can likewise encode sutured graph complements via multipointed Heegaard diagrams with a basepoint corresponding to the suture associated to a given edge in the spatially embedded graph.
When the multipointed Heegaard diagram happens to be a suitable generalization of a grid diagram, Alishahi--Eftekhary's theory recovers the grid homology for transverse spatial graphs described by Harvey--O'Donnol \cite{HO17}.  
\begin{remark}
In the literature, tangle Floer homology sometimes refers to a bordered invariant for tangles in the 3-ball described by Zibrowius \cite{Zib16}. 
We emphasize that the Alishahi--Eftekhary invariants are quite different from those of Zibrowius--- Alishahi--Eftekhary's theory
does not possess bordered pairing theorems since boundary parametrizations are absent from the input data. On the other hand, the Zibrowius tangle invariants can be used to compute particular specializations of the Alishahi--Eftekhary invariants, analogous to how bordered techniques from \cite{LOT18, Han23} can be used to compute $\widehat{\mathrm{HFK}}$ and $\mathrm{HFK}_{UV=0}$ respectively.
\end{remark}
In \cite{Ali22}, Alishahi used the tangle Floer machinery to define the graph upsilon $\Upsilon$, an invariant of spatial graphs which specializes to similar constructions of Ozsv\'{a}th--Stipsicz--Szab\'{o} \cite{OSS17a},  Cavallo 
\cite{Cav24}, and Kubota \cite{Kub23}.

We review the construction of Alishahi--Eftekhary's tangle Floer complex in \S~\ref{sec:prelimTangleFloer}. 
In \S~\ref{sec:absGradings}, we describe a general method for equipping $\mathrm{CFT}$ with relative and absolute gradings and in \S~\ref{sec:Upsilon} we give a definition of Alishahi's graph $\Upsilon$. In \S~\ref{sec:stabilizations}, we describe two stabilization operations on abstract graphs and their effect on tangle Floer homology. 
One of these is \textit{edge stabilization} and is identical to the quasistabilization operations for knots and links studied in \cite{MO25,Zem17,Zem19b}. 
For graphs containing a vertex of valence $\geq 2$, we have a notion of \textit{vertex stabilization} which will be important for our purposes. We conclude in \S~\ref{sec:squareGradings} with a geometrically motivated grading on the tangle Floer complex of certain graphs which we call the \textit{square grading}.

\subsection{The tangle Floer complex for a spatial graph}\label{sec:prelimTangleFloer}
For our setting, an oriented tangle is a pair ${ (M,T) }$ such that:
\begin{itemize}
\item ${ M }$ is a 3-manifold with no closed components and whose boundary splits as ${ \partial M = -\partial_{-}M \sqcup \partial_{+}M }$.
\item ${ T = \bigsqcup_{i=1}^{k} T_{i} }$ is a union of disjoint, properly embedded, and oriented arcs ${ T_{i} }$.
\item  With respect to the orientation on the components of $T$, we assume for all $i$ that the starting point of $T_i$ is contained in $\partial_- M$ and its terminal point is contained in $\partial_+ M$. 
\end{itemize}

\begin{remark}
Although the theory permits $M$ to be disconnected, we focus on the connected case for simplicity. We will also assume throughout that our tangles are \textit{balanced}, meaning that $\chi(\partial_{+}M) = \chi(\partial_{-}M)$ and every connected component of $\partial M$ meets $\partial T$. 
\end{remark}

Let $(M,T)$ be an oriented tangle in a 3-manifold whose boundary is a disjoint union of spheres. Enumerating components of $T$ as $\{T_j\}_{j=1}^n$, we may associate to each component ${ T_{j} }$ a formal variable ${ Z_{j} }$ in the polynomial ring ${R_n := (\mathbb{Z} / 2\mathbb{Z})[Z_{1}, \dots, Z_{n}] }$. 
A \textit{Heegaard diagram for an ${ A}$-tangle} ${(M, T) }$ is a quintuple
$\mathcal{H} = (\Sigma,\bm{\alpha},\bm{\beta}, \mathfrak{c}: \bm{z} \rightarrow A)$
where:
\begin{itemize}
\item ${ \Sigma }$ is an embedded surface in ${ M }$ which splits ${ M }$ into compression bodies ${ H_{\alpha} }$ and ${ H_{\beta} }$ and
$$ \partial H_{\alpha} = -(\partial_{-}M ) \cup \Sigma \qquad \text{ and } \qquad \partial H_{\beta} = -\Sigma \cup \partial_{+}M .   $$
\item ${ \bm{\alpha} = \bigsqcup_{i=1}^{\ell} \alpha_{i} }$ is a union of embedded simple closed curves in ${ \Sigma }$ such that each ${ \alpha_{i} }$ bounds a disk in ${ H_{\alpha} }$ and, compressing ${ \Sigma }$ along these disks, we obtain a surface isotopic to ${ \partial_{-} M }$ inside of ${ M }$. 
\item ${ \bm{\beta} }$ is defined identically to $\bm{\alpha}$ except that the result of compressing along disks bounded by the ${ \beta_{i} }$ gives a surface isotopic to ${ \partial_{+}M }$ inside of ${ M }$.
\item ${ \bm{z}= \bigsqcup_{i=1}^{n} z_{i} = \bigsqcup_{i=1}^{n} (T_{i}\cap \Sigma) }$ is a family of points in ${ \Sigma \setminus (\bm{\alpha}\cup \bm{\beta}) }$. These correspond to intersections of $T$ with the embedded Heegaard surface $\Sigma$; it will suffice for our purposes to insist that each component of ${ T }$ intersects ${ \Sigma }$ exactly once.
\item ${ \mathfrak{c}: \bm{z}\rightarrow A }$ is a coloring map to a ring $A$ which is a quotient of $R_n$. In particular, there is a universal coloring $\mathfrak{u}: \bm{z} \rightarrow R_n$ defined by ${ \mathfrak{u}(z_{i}) = Z_{i} }$.
\end{itemize}

There is a function $\mathfrak{s}: \mathbb{T}_{\alpha} \cap \mathbb{T}_{\beta} \to \mathrm{Spin}^c(M)$ and when ${ \mathcal{H} }$ is $\mathfrak{t}$-admissible for $\mathfrak{t} \in \mathrm{Spin}^c(M)$ (see \cite{AE15}, Definition 4.1) and $\mathfrak{c} : \bm{z} \rightarrow A$ is some coloring map, one may associate to $\mathcal{H}$ an $A$-module
${ \mathrm{CFT}(\mathcal{H},\mathfrak{t}) }$ freely generated by intersection points ${ \bm{x} \in \mathbb{T}_{\alpha} \cap \mathbb{T}_{\beta} \subset \mathrm{Sym}^{\ell}(\Sigma) }$ satisfying $\mathfrak{s}(\bm{x})=\mathfrak{t}$. This module possesses an endomorphism
$$ \partial_{\mathfrak{c}} (\bm{x}) = \sum\limits_{\bm{y} \in \mathbb{T}_{\alpha} \cap \mathbb{T}_{\beta} } \sum\limits_{\substack{\phi \in \pi_{2} (\bm{x},\bm{y})\\ \mu(\phi)=1}} \# \widehat{\mathcal{M}}(\phi)\mathfrak{c}(\phi)\cdot \bm{y}    $$ where ${ \mathfrak{c}(\phi) := \prod_{i=1}^{n} \mathfrak{c}(z_{i})^{n_{z_{i}}(\phi)} }$ and ${ \# \widehat{\mathcal{M}}(\phi) }$ is the $\mathbb{Z}/2\mathbb{Z}$-count of holomorphic representatives (modulo $\mathbb{R}$-translation) of the Whitney disk ${ \phi }$ connecting ${ \bm{x} }$ and ${ \bm{y} }$.

\begin{definition}[Curved chain complex]
A \textit{curved chain complex} over a ring $A$ is a pair $(C,\partial)$ where $C$ is an $A$-module and $\partial$ is an endomorphism of $C$ such that for a fixed $\omega \in A$ we have 
$\partial^2 = \omega \cdot \mathrm{Id}$. We may refer to $\omega$ as the \textit{curvature} of $C$.
\end{definition}

For a fixed coloring $\mathfrak{c}: \bm{z}\rightarrow A$, the $A$-module
$\mathrm{CFT}(\mathcal{H})$ with the endomorphism $\partial_{\mathfrak{c}}$ forms a curved complex whose curvature is obtained from incidence relations of the tangle components with components in the sutured boundary ${ \partial_{\pm}M }$:
\begin{equation}\label{eqn:curvature}
\omega_{\mathfrak{c}}
= \left(\sum\limits_{C_{+} \in \pi_0(\partial_{+}M)  } \prod\limits_{T_{j} \cap C_{+} \neq \emptyset  } \mathfrak{c}(z_j) \right) - \left(\sum\limits_{C_{-} \in \pi_0(\partial_{-}M)  } \prod\limits_{T_{k} \cap C_{-} \neq \emptyset  } \mathfrak{c}(z_k)\right) .
\end{equation}
If $\mathcal{H}$ encodes a spatial graph such that $\omega_{\mathfrak{u}} =0$, the uncolored complex $\mathrm{CFT}(\mathcal{H},\mathfrak{u})$ is an ordinary chain complex. More generally, one may obtain a chain complex by swapping  to a coloring $\mathfrak{c} : \bm{z} \rightarrow A$ which factors through a ring homomorphism $f: R_{n} \rightarrow A$ whose kernel contains the ideal $\langle \omega_{\mathfrak{u}} \rangle$.

Given a colored Heegaard diagram $\mathcal{H}$, it is shown in \cite{AE15} that the homotopy type of $\mathrm{CFT}(\mathcal{H})$ is an invariant of $(M_G,T_G)$. Moreover, if any pair of diagrams ${ \mathcal{H} }$ and ${ \mathcal{H}' }$ are $\mathfrak{c}$-colored diagrams for $(M_G,T_G)$, then they differ from one another by a sequence of Heegaard moves and, to this sequence, there is an induced homotopy equivalence between ${ \mathrm{CFT}(\mathcal{H}) }$ and ${ \mathrm{CFT}(\mathcal{H}') }$.

To an oriented, embedded graph $G: \Gamma \hookrightarrow Y^3$, one may associate an oriented tangle $T_{G}$ in a 3-manifold with boundary $M_G$ by removing neighborhoods of its vertices. When the choice of coloring is relevant, we write $\mathrm{CFT}(G,\mathfrak{c})$ for the homotopy type of a ${ \mathfrak{c} }$-colored tangle Floer complex associated to the triple $(M_G,T_{G},\mathfrak{c})$.

\subsection{Gradings and the $\Upsilon$ invariant}\label{sec:absGradings}

\subsubsection{Universal grading data}
Let $G : \Gamma \hookrightarrow S^3$ be an embedding of a balanced, bipartite graph with $n$ edges and $2k$ vertices. We suppress the distinguished $\mathrm{Spin}^c$-structure throughout; see \cite[\S~3.2]{Ali22} for the corresponding conditions in a general orientable, closed 3-manifold.
For each vertex $v$, consider the formal sum ${ s_{v} = \sum_{e \ni v } e \in \mathbb{Z}^{E(\Gamma)} }$.
We define the \textit{universal Heegaard Floer grading group} for the graph ${ \Gamma }$ to be 
$$ \mathcal{G}_{\Gamma} := \frac{\mathbb{Z} \oplus \mathbb{Z}^{E(\Gamma)} }{ \langle (0,s_{v}) - (2,0) : v \in V(\Gamma)    \rangle } .$$

We abuse notation and use $e \in E(\Gamma)$ to denote its image in $\mathcal{G}_{\Gamma}$. The tangle Floer coefficient ring for $\Gamma$, $R_{\Gamma} = \mathbb{F} [Z_e : e \in E(\Gamma)]$, is $\mathcal{G}_{\Gamma}$-graded by $\deg(Z_e) = -e$.

\begin{proposition}
The uncolored (and possibly curved) complex $\mathrm{CFT}(\mathcal{H})$ of a Heegaard diagram for $G: \Gamma \hookrightarrow S^3$ has a relative $\mathcal{G}_{\Gamma}$-grading such that, for any $\bm{x},\bm{y} \in \mathbb{T}_{\alpha} \cap \mathbb{T}_{\beta}$, it satisfies the conditions:
\begin{equation}\label{eqn:universalGraphGrading}
 \begin{aligned}
gr_{\Gamma}(Z_e\, \bm x)
 &=gr_{\Gamma}(\bm x)-e,\\                          
gr_{\Gamma}(\bm x)-gr_{\Gamma}(\bm y)      
 &=\mu(\phi)-\sum_{e\in E(\Gamma)}n_{z_e}(\phi)e,
 \qquad \phi\in\pi_2(\bm x,\bm y).
 \end{aligned}
 \end{equation}    
 Moreover, the differential has degree $(-1,0) \in \mathcal{G}_{\Gamma}$.
 \begin{proof}
The ambiguity in the second formula comes from periodic domains.
Since the ambient manifold is $S^3$, every periodic domain is an
integral combination of the components of $\Sigma\setminus\bm\alpha$
and $\Sigma\setminus\bm\beta$. The component corresponding to a
vertex $v$ has Maslov index $2$ and multiplicity vector $s_v$, so
its contribution is $(2,0)-(0,s_v)=0$ in $\mathcal{G}_\Gamma$.
Thus the formula is independent of $\phi$; compare
\cite[Lemma 3.7]{Ali22}. 
Additivity under concatenation of domains gives
the relative grading. The tangle Floer differential counts disks of Maslov index one, so has degree
$(-1,0)$.
\end{proof}
\end{proposition}

For an arbitrary abelian group ${ H }$, 
the \textit{$\Gamma$-compatible ${ H }$-grading data} are 
$$\mathrm{Grad}_{H}(\Gamma) := \mathrm{Hom}(\mathcal{G}_{\Gamma},H)$$
and, given some $\rho \in \mathrm{Grad}_H(\Gamma)$, one obtains a relative $H$-grading $gr_{\rho}$ on $\mathrm{CFT}(\mathcal{H})$ in the obvious manner. 
In particular, it satisfies
$$
\begin{aligned}
gr_{\rho}(Z_e\, \bm x)
 &=gr_{\rho}(\bm x)-\rho(e) \\ \forall \phi \in \pi_2(\bm{x},\bm{y}): gr_{\rho}(\bm x)-gr_{\rho}(\bm y)      
 &=\mu(\phi)\rho(1,0)-\sum_{e\in E(\Gamma)}n_{z_e}(\phi)\rho(e).
 \end{aligned}
 $$

Recall that a \textit{perfect matching} for an abstract graph ${ \Gamma }$ is a set of edges ${ \mathfrak{m} }$ meeting each vertex exactly once. 
Let ${ \mathfrak{M}_{\Gamma} }$ denote the set of perfect matchings for ${ \Gamma }$. 
Each matching defines a homomorphism
$$
\rho_{\mathfrak{m}}:\mathcal G_\Gamma\rightarrow \mathbb{Z},
\qquad
\rho_{\mathfrak{m}}(1,0)=1,
\qquad
\rho_{\mathfrak{m}}(0,e)=
\begin{cases}
2,&e \in \mathfrak{m}\\
0,&e \notin \mathfrak{m}
\end{cases}.
$$
The relative grading for $\rho_{\mathfrak{m}}$ has a preferred absolute lift
$\operatorname{gr}_{\mathfrak m}$ which is compatible with the absolute Maslov gradings for link Floer homology:
forget the basepoints outside
$\mathfrak m$, equivalently set $Z_e=1$ for $e\notin\mathfrak{m}$,
and use the absolute Maslov grading of the resulting multipointed
Heegaard diagram for $S^3$ \cite[\S~3.2.2]{Ali22}.
We use the integral normalization in which the top class of the hat
specialization has degree zero. 
For $h\in H$, the grading
$h\operatorname{gr}_{\mathfrak m}$ corresponds to the homomorphism
with $1 \mapsto h$ and $e\mapsto2h\mathbf1_{\mathfrak m}(e)$.

\subsubsection{$\mathbb{R}$-valued gradings and perfect matching polytopes}
Setting $H = \mathbb{R}$, we consider a distinguished subset of the set of $\mathbb{R}$-grading data:
$$C_{\Gamma} = \{ \rho \in \mathrm{Grad}_{\mathbb{R}}(\Gamma) : \rho(0,e) \geq 0\}.$$

\begin{lemma}\label{thm:perfectMatchingsGeneratePositiveCone} 
We have 
  ${ C_{\Gamma} = \mathrm{Cone}\left\{\rho_{\mathfrak{m}} \; | \; \mathfrak{m} \in \mathfrak{M}_{\Gamma}\right\}  }$.
 \begin{proof}
 This is essentially the same argument as that given in \cite[Lemma 3.3]{Ali22}; the only change is that one is varying the degree of the tangle Floer differential.
\end{proof}
\end{lemma}

Recall that the \textit{perfect matching polytope} for $\Gamma$, denoted $\mathcal{P}_{\Gamma}$, is the convex polytope obtained as the convex hull of the matching incidence vectors in $\mathbb{R}^n$.
The hyperplane  
$\{\rho(1,0) = 1 \in \mathbb{R}\} \subset \mathrm{Grad}_{\mathbb{R}}(\Gamma)$ intersects the cone $C_{\Gamma}$ in the rescaled perfect matching polytope $2\mathcal{P}_{\Gamma}$. This intersection exactly corresponds to the set of relative gradings for nonnegative $\mathbb{R}$-valued edge weights which admit absolute lifts and for which the Floer differential has degree $-1$.

It is illuminating to consider how the number of matchings and linear dependency across their associated grading homomorphisms impact the existence and uniqueness of distinguished absolute $\mathbb{R}$-gradings; to this end we consider several examples.

\begin{example}\label{exmp:noMatchings}
Not all balanced bipartite graphs admit perfect matchings; below is a connected example 
$$\Gamma \qquad = \qquad  
\begin{tikzpicture}[baseline=(current bounding box.center),
                    x=1cm,y=0.7cm,
                    every node/.style={inner sep=1pt}]
  \coordinate (a) at (0, 1);
  \coordinate (b) at (0,-1);
  \coordinate (c) at (1, 0);
  \coordinate (d) at (2, 0);
  \coordinate (e) at (3, 1);
  \coordinate (f) at (3,-1);

  \draw (a)--node[midway,above right] {$e_1$} (c);
  \draw (b)--node[midway,below right] {$e_2$} (c);
  \draw (c)--node[midway,above]       {$e_3$} (d);
  \draw (d)--node[midway,above left]  {$e_4$} (e);
  \draw (d)--node[midway,below left]  {$e_5$} (f);

  \foreach \v in {a,b,d}
    \fill (\v) circle[radius=2pt];
  \foreach \v in {c,e,f}
    \draw[fill=white] (\v) circle[radius=2pt];
\end{tikzpicture}
$$
It is not hard to see that $\mathcal{G}_{\Gamma} \cong \mathbb{Z}$ and corresponds to (the image of) the $\mathbb{Z}$-summand in $\mathbb{Z}\oplus \mathbb{Z}^{E(\Gamma)}$. A bit of work shows the valid edge weights associated to an $\mathbb{R}$-grading are precisely $\{(2t,2t,-2t,2t,2t) \; | \; t\in \mathbb{R}\}$.
From the above, we see that $C_{\Gamma} = \{0\}$ and hence $C_{\Gamma} \cap \{\rho \; | \; \rho(1,0) =1\} = \emptyset = 2\mathcal{P}_{\Gamma}$. 
Each of the non-zero gradings has associated edge colorings containing at least one negative color; hence none of these relative gradings admit a preferred absolute lift.
\end{example}

\begin{example}\label{exmp:K_{3,3}}
In contrast to Example \ref{exmp:noMatchings}, one can also have ``too many matchings'' in the sense that ${|\mathfrak{M}_{\Gamma}| > \mathrm{dim}( \mathrm{Grad}_{\mathbb{R}}(\Gamma))}$; this introduces linear dependencies among the matching homomorphisms of $\mathrm{Hom}(\mathcal{G}_{\Gamma},\mathbb{R})$ and hence among
those morphisms contained in 
$C_{\Gamma} \cap \{\rho \; | \; \rho(1,0) =1\}$.
In this situation, a preferred lift is not distinguished by the weight vectors alone.

Consider the complete bipartite graph ${ K_{3,3} }$. The rescaled Birkhoff polytope ${ 2\mathcal{B}_{3} \subset \mathbb{R}^{9} }$ has vertices corresponding to matchings of ${ K_{3,3} }$; these are in bijection with elements of the symmetric group ${ S_{3} }$ and represented in $\mathbb{R}^9$-coordinates as $3\times 3$ permutation matrices, rescaled by 2, which we enumerate as
$$ \mathbf{t}_{123}, \; \mathbf{t}_{231}, \; \mathbf{t}_{312}, \; \mathbf{t}_{132}, \; \mathbf{t}_{213},\; \mathbf{t}_{321} \in \mathbb{R}^{9}.      $$ 
With a bit of work, one finds $\mathrm{dim}(\mathrm{Grad}_{\mathbb{R}}(K_{3,3}))=5$.
Indeed, there is a linear dependence among the ${ \mathbf{t}_{i} }$:
$$\mathbf{t}_{123} + \mathbf{t}_{231} + \mathbf{t}_{312} = \mathbf{t}_{132} + \mathbf{t}_{213} + \mathbf{t}_{321} = 2\cdot \mathbf{1} := 2\cdot (1,\cdots, 1) .  $$
In particular, ${ \mathbf{t}' := \frac{2}{3}  \mathbf{1} \in \mathrm{relint}(2\mathcal{B}_{3}) }$ arises from two distinct convex combinations
$$ \mathbf{t}' = \frac{1}{3}(\mathbf{t}_{123} + \mathbf{t}_{231} + \mathbf{t}_{312}   ) = \frac{1}{3}(\mathbf{t}_{132} + \mathbf{t}_{213} + \mathbf{t}_{321}   ) $$
and we obtain two absolute gradings at ${ \mathbf{t}' }$:
\begin{equation}\label{eqn:even_and_odd_gradings} 
gr_{\mathbf{t}'}^{ev} := \frac{1}{3}(gr_{123} + gr_{231} + gr_{312}   ) \qquad \text{ and } \qquad gr_{\mathbf{t}'}^{odd} := \frac{1}{3}(gr_{132} + gr_{213} + gr_{321}  )   .
\end{equation}

These gradings are different.
Figure \ref{fig:G3_full} displays Heegaard states and the (curved) Floer complex differential for a diagram $\mathcal{H}= (\mathbb{T}^2, \bm{\alpha}, \bm{\beta}, \bm{z}) $ which describes a particular embedding of ${ K_{3,3} \hookrightarrow  S^{3} }$. 
${ \mathcal{H} }$ is such that $\bm{\alpha}$ is a union of three parallel longitudes of $\mathbb{T}^2$, $\bm{\beta}$ is three parallel meridians, and every connected component of $\mathbb{T}^2 \setminus (\bm{\alpha}\cup \bm{\beta})$ contains a basepoint.\footnote{Heegaard diagrams of this type are often called \textit{grid diagrams}; see \S~\ref{sec:gridDiagrams} for further discussion.} The set of matchings $\mathfrak{M}_{K_{3,3}}$ and the Heegaard states of $\mathcal{H}$ are in bijection with elements of $S_3$ and one has $gr_{\mathfrak{m}_{\sigma}}(\bm{x}_{\sigma})=0$ for all $\sigma \in S_3$.
Table \ref{table:G_3full_gradings} displays the $gr_{\mathfrak{m}_{\sigma}}$ gradings of all Heegaard states and shows that the absolute gradings $gr_{\mathbf{t}'}^{odd}$ and $gr_{\mathbf{t}'}^{ev}$ described in Equation \ref{eqn:even_and_odd_gradings} differ from one another by $\tfrac{1}{3}$.

\begin{figure}
\includegraphics[width=\linewidth]{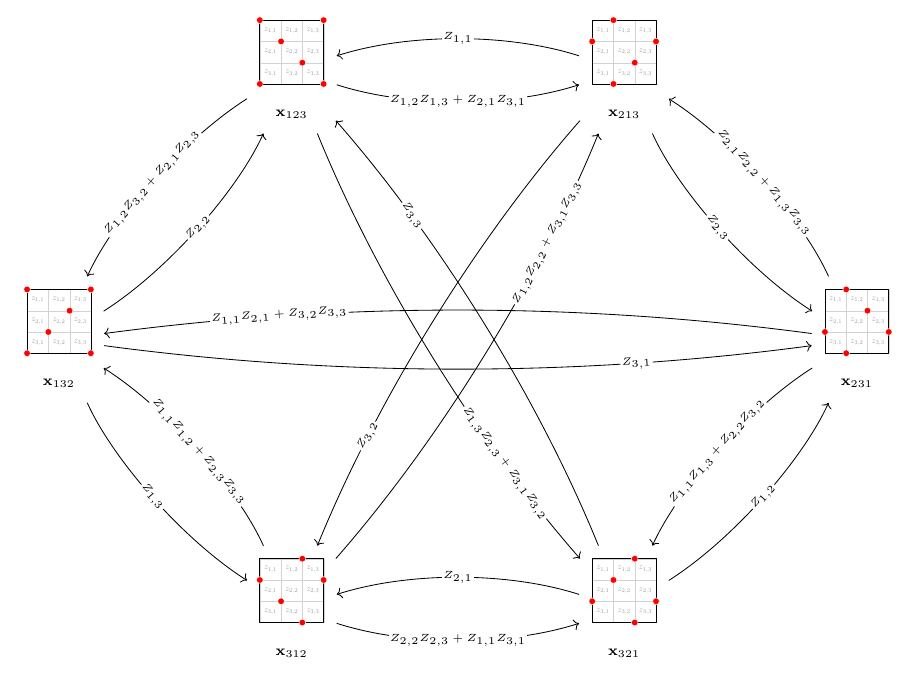}
\caption{
Example \ref{exmp:K_{3,3}} 
describes an embedding $K_{3,3} \hookrightarrow S^3$ encoded by a $3\times 3$ grid diagram $\mathcal{H}$.   
The above encodes the uncolored tangle Floer complex 
$\mathrm{CFT}(\mathcal{H})$ of the diagram; this is a curved complex 
over $\mathbb{F}[Z_{i,j}]_{i,j=1}^3$ 
with curvature 
$\omega_{K_{3,3}} = \sum_{i=1}^3 (Z_{i,1}Z_{i,2}Z_{i,3} + Z_{1,i}Z_{2,i}Z_{3,i}).$
}\label{fig:G3_full}
\end{figure}

\begin{table}[!h]
\centering
$$
\begin{array}{|c|cccccc|cc|} \hline
&
\operatorname{gr}_{123}
&
\operatorname{gr}_{213}
&
\operatorname{gr}_{132}
&
\operatorname{gr}_{231}
&
\operatorname{gr}_{312}
&
\operatorname{gr}_{321}
&
\phantom{\Big)} \operatorname{gr}_{\mathbf{t}'}^{\mathrm{ev}} \phantom{\Big)}
&
\operatorname{gr}_{\mathbf{t}'}^{\mathrm{odd}} 
\\ \hline
\mathbf{x}_{123}
& 0 & -1 & -1 & -2 & -2 & -1
& -\frac{4}{3} & -1 \phantom{\Big)}
\\ \hline
\mathbf{x}_{213}
& -1 & 0 & -2 & -1 & -1 & 0
& -1 & -\frac{2}{3} \phantom{\Big)}
\\ \hline
\mathbf{x}_{132}
& -1 & 0 & 0 & -1 & -1 & -2
& -1 & -\frac{2}{3} \phantom{\Big)}
\\ \hline
\mathbf{x}_{231}
& -2 & -1 & -1 & 0 & -2 & -1
& -\frac{4}{3} & -1 \phantom{\Big)}
\\ \hline
\mathbf{x}_{312}
& -2 & -1 & -1 & -2 & 0 & -1
& -\frac{4}{3} & -1 \phantom{\Big)}
\\ \hline
\mathbf{x}_{321}
& -1 & -2 & 0 & -1 & -1 & 0
& -1 & -\frac{2}{3} \phantom{\Big)} 
\\ \hline
\end{array}
$$
\caption{}\label{table:G_3full_gradings}
\end{table}

\end{example}

We conclude this section with an example where no complications arise whatsoever.
\begin{example}\label{exmp:thetaGraphs}
For $\Gamma = \Theta_n$, one has that the grading data $\mathrm{Grad}_{\mathbb{R}}(\Theta_n)$ is spanned by homomorphisms associated to perfect matchings. Moreover, there is no linear dependence among the homomorphisms for the matching gradings; hence each grading in $C_{\Gamma}$ admits a distinguished absolute lift.
The intersection of the hyperplane $\{\rho(1,0)=1\}$ with $C_{\Theta_n}$ is a rescaled $(n-1)$-simplex: $2 \Delta^{n-1}$. 
\end{example}

\subsubsection{The $\Upsilon$-invariant of knots and its generalization to spatial graphs}\label{sec:Upsilon}

Given a doubly pointed knot ${ (K,w, z) }$, Ozsv\'{a}th--Stipsicz--Szab\'{o} associate to each ${ t \in [0,2] }$ a ``${ t }$-modified''
version of the knot Floer complex ${ t\mathrm{CFK} }$ which possesses a grading ${ gr_{t} }$ \cite{OSS17a}.
Although they did not write about it in such terms, their ${ gr_{t} }$ grading can be defined via a formula
\begin{equation}\label{eqn:convex_maslov_grading}
gr_{t} = (1-\tfrac{t}{2})gr_{w} + \tfrac{t}{2} gr_{z}  
\end{equation}
where, as in \cite{Zem19a}, $gr_w$ and $gr_z$ denote the absolute Maslov gradings for $\mathcal{CFK}(K)$ respectively associated to the basepoints $w$ and $z$. 
A doubly pointed knot may be regarded as an embedding of the loopless graph $\Theta_2$ with two edges and two vertices whose matchings are ${ \mathfrak{m}_{w} }$ and ${ \mathfrak{m}_{z} }$. 
As indicated in (\ref{eqn:convex_maslov_grading}), one may regard ${ gr_{t} }$ as a convex combination of the two absolute Maslov gradings ${ gr_{w} = gr_{\mathfrak{m}_{w}} }$ and ${ gr_{z} = gr_{\mathfrak{m}_{z}} }$ and these convex combinations are indexed by positions in the rescaled perfect matching polytope ${ 2\mathcal{P}_{\Theta_2} =[0,2]}$. With this in hand, Ozsv\'{a}th--Stipsicz--Szab\'{o} define the $\Upsilon$ invariant to be 
$$\Upsilon_K(t) = \max \left\{ gr_t(x) \; \middle| \; \begin{matrix} x \text{ homogeneous} \\ x \in H_{\ast}\big(t\mathrm{CFK}(K)\big) \setminus \mathrm{Tors} \end{matrix}  \right\} .$$

Alishahi \cite{Ali22} observed that one can define an analogue of ${ t\mathrm{CFK} }$ for any bipartite, balanced spatial graph with perfect matchings. In particular, to any  $\bm{t} \in 2\mathcal{P}_{\Gamma}$, the previous subsection constructs a relative grading $gr_{\bm{t}}$ from a convex combination of the relative Maslov gradings ${ gr_{\mathfrak{m}} }$.
When ${ \mathcal{P}_{\Gamma} }$ is a simplex --- as is the case for ${ \Theta_{n} }$-graphs --- these $gr_{\bm{t}}$ gradings admit unique absolute lifts.

\begin{warning}
We caution the reader that when ${ \mathcal{P}_{\Gamma} }$ is not a simplex, Alishahi's ${ gr_{\bm{t}} }$ gradings are not determined intrinsically by the abstract graph alone. We have already seen this subtlety in Example \ref{exmp:K_{3,3}}. Alishahi works around this by fixing an ordering on the edges of $\Gamma$ and using an order-dependent triangulation of $2\mathcal{P}_{\Gamma}$ whose vertices are vectors corresponding to matchings; this suffices to specify a distinguished absolute lift of $gr_{\bm{t}}$ for any $\bm{t} \in 2\mathcal{P}_{\Gamma}$.
\end{warning}

In \cite[Lemmas 4.1 and 4.2]{Ali22} it is shown that if ${ G : \Gamma \hookrightarrow S^3 }$ is a spatial graph where $\Gamma$ is balanced, bipartite, and has perfect matchings, then one may obtain a colored tangle Floer complex $\bm{t}\mathrm{CFT}(G)$ for each $\bm{t}\in 2\mathcal{P}_{\Gamma}$.\footnote{A detail we have suppressed is that, in Alishahi's account, the complex $\bm{t}\mathrm{CFT}(G)$ is a module over the long power series ring $\mathcal{R}$.
} 
If $\Gamma$ has $2k$ vertices, then the homology
${ \bm{t}\mathrm{HFT}(G) }$ of this complex has the property that ${ \bm{t}\mathrm{HFT}(G)\, /\, \mathrm{Tors}  }$ is a free module of rank ${ 2^{k-1} }$ and the ${ gr_{\bm{t}}}$ gradings of this homogeneous basis are invariants of the spatial graph ${ G }$.

\begin{definition}[$\Upsilon$ for spatial graphs]\label{defn:GraphUpsilon} 
  Fix ${ \bm{t} \in 2\mathcal{P}_{\Gamma} }$ and let ${ \left\{e_{i}(\bm{t})\right\}_{i=1}^{2^{k-1}} }$ be a basis of homogeneous elements generating ${ \bm{t}\mathrm{HFT}(G)\, / \, \mathrm{Tors} }$ such that $gr_{\bm{t}}(e_{i}(\bm{t})) \leq gr_{\bm{t}}(e_{i+1}(\bm{t}))$. Then ${ \Upsilon_{G}(\bm{t}) }$ is a nondecreasing sequence of ${ 2^{k-1} }$ real numbers 
  $$ \Upsilon_{G}(\bm{t}) = (\Upsilon_{i} )_{i=1}^{2^{k-1} }  \qquad \Upsilon_{i} = gr_{\bm{t}} (e_{i} (\bm{t}))  . $$ 
 \end{definition}

In the case where $G$ is a $\Theta_n$-graph, the $\Upsilon$ invariant is a single, piecewise linear function over $2 \Delta^{n-1}$. Equivalently, it is an $\mathbb{R}$-valued invariant of $\mathbb{R}_{\geq 0}$-colored $\Theta_n$-graphs whose edge colorings sum to 2.
We conclude with a summary of properties which are relevant for this note.
\begin{enumerate}
\item The restriction of $\Upsilon_G$ to any 1-cell of  $2\Delta^{n-1}$ agrees with the knot $\Upsilon$ of an associated cycle of $G$. 
\item Given $G,G' : \Theta_n \hookrightarrow S^3$ together with an ordering on their edges, the \textit{vertex connect sum} $G \#_v G'$ is the $\Theta_n$-graph obtained by puncturing out a vertex from each and joining the resulting half-edges together according to the identification specified by the edge ordering. 
Alishahi shows that $\Upsilon_{G\#_v G'} = \Upsilon_G + \Upsilon_{G'}$. 
\item In fact, $\Upsilon$ is invariant under concordance of $\Theta_n$-graphs; hence the above tells us that it defines a concordance homomorphism from the smooth concordance group of $\Theta_n$-graphs to $\mathrm{PL}(2\Delta^{n-1})$.
\end{enumerate}

\begin{figure}[!h]
    \centering

    \begin{subfigure}[t]{0.48\linewidth}
        \centering
        \includegraphics[width=\linewidth]{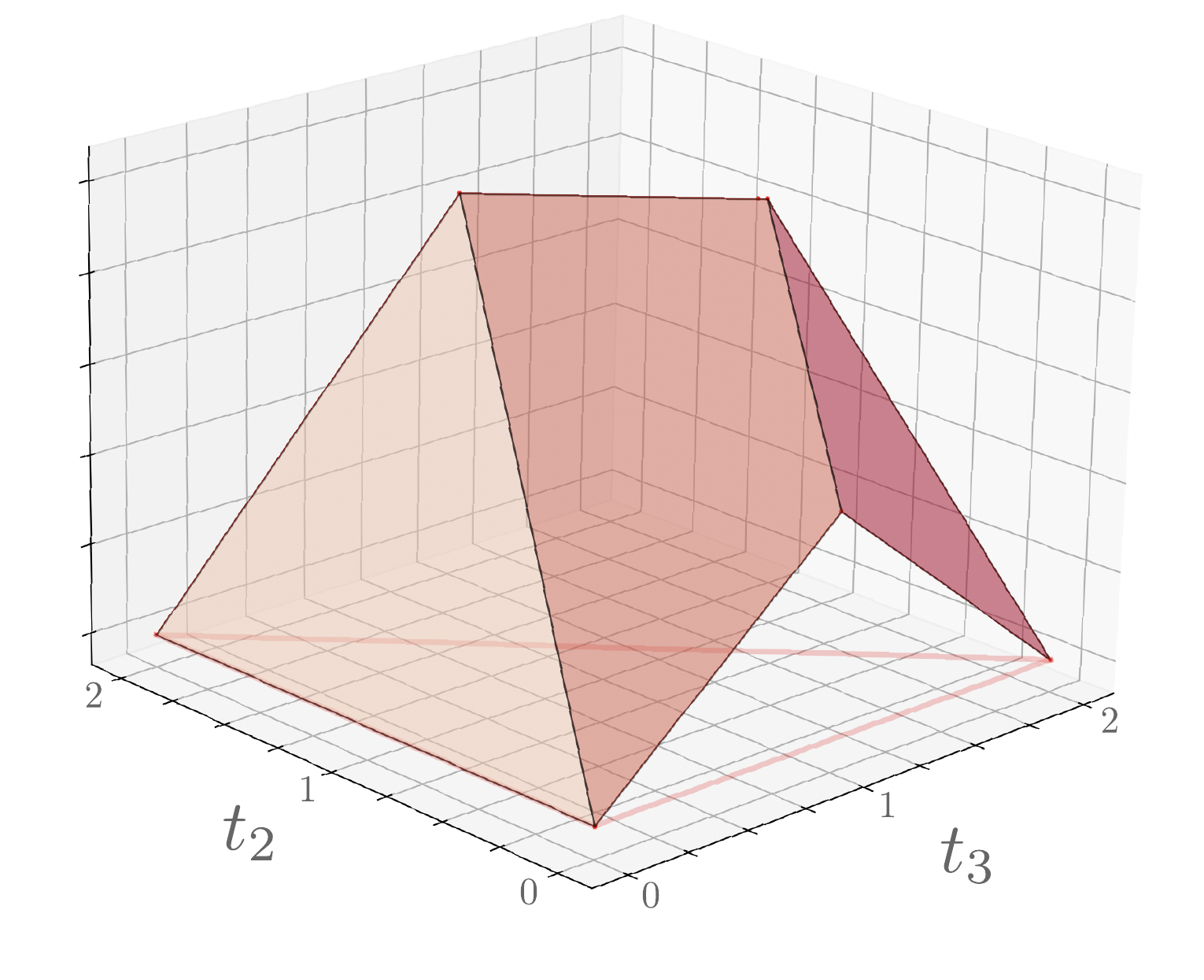}
        \caption{$\Theta_3\begin{psmallmatrix} 3 & 1 \\ 2 & 1 \end{psmallmatrix} \subset \Theta_4\begin{psmallmatrix} 3 & 4 \\ 2 & 3 \end{psmallmatrix}$ \\ ($t_1=0 \text{ facet}$)}
        \label{fig:3_4_2_3_smaller_upsilon}
    \end{subfigure}
    \hfill
    \begin{subfigure}[t]{0.48\linewidth}
        \centering
        \includegraphics[width=\linewidth]{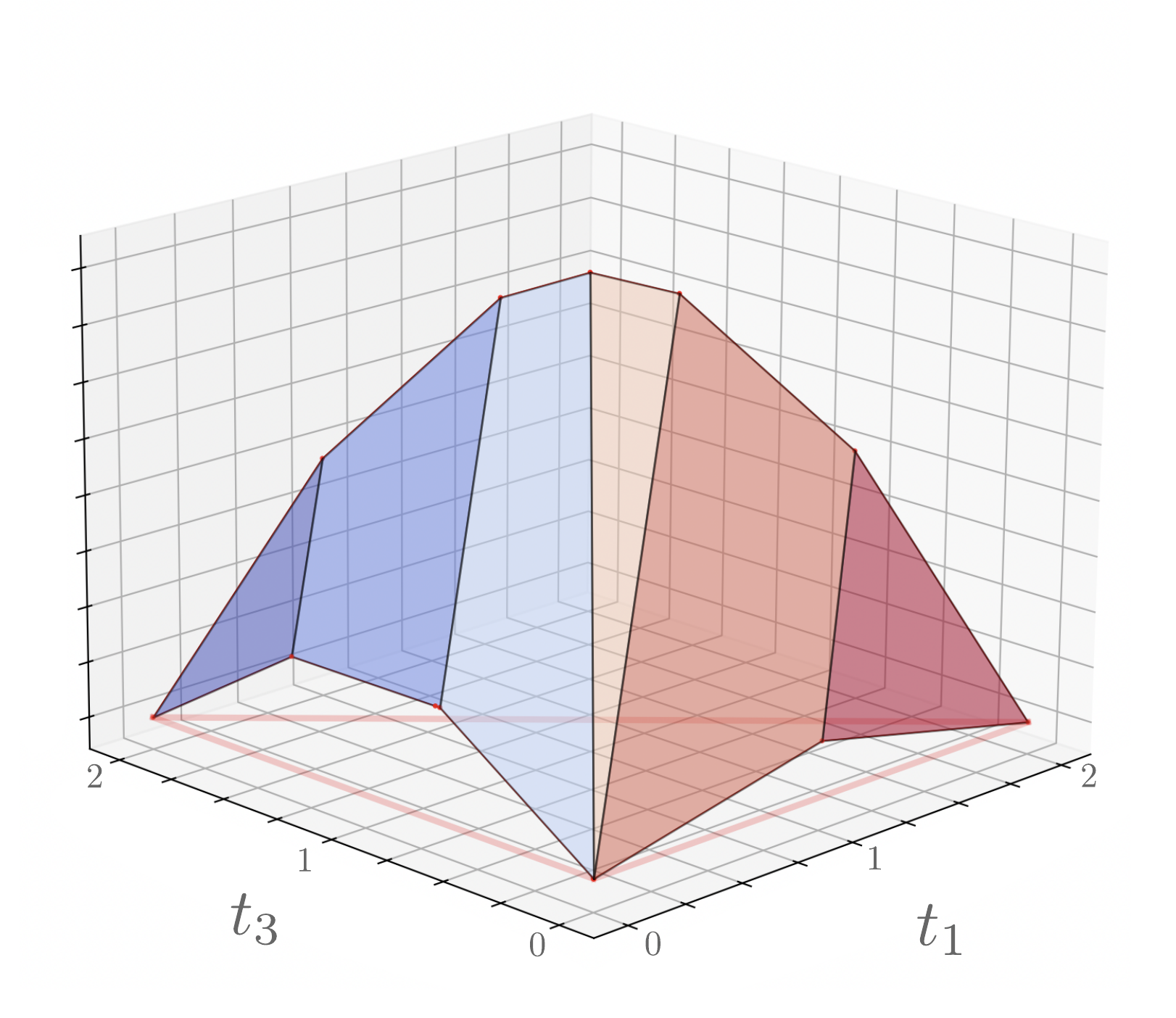}
        \caption{$\Theta_3\begin{psmallmatrix} 3 & 4 \\ 2 & 3 \end{psmallmatrix} \subset \Theta_4\begin{psmallmatrix} 3 & 4 \\ 2 & 3 \end{psmallmatrix}$ \\  ($t_4=0 \text{ facet}$)}
        \label{fig:3_4_2_3_bigger_upsilon}
    \end{subfigure}

    \caption{We computed $\mathrm{CFT}(\Theta_4
    \begin{psmallmatrix} 3 & 4 \\ 2 & 3 \end{psmallmatrix})$ in Example \ref{exmp:3_4_2_3_theta}. The above displays plots of $-\tfrac{1}{2}\Upsilon$ for the two distinct $\Theta_3$-subgraphs of $\Theta_4
    \begin{psmallmatrix} 3 & 4 \\ 2 & 3 \end{psmallmatrix}$; each may be obtained by restricting $-\frac{1}{2}\Upsilon_{\Theta_4
    \begin{psmallmatrix} 3 & 4 \\ 2 & 3 \end{psmallmatrix}}$ to two facets of $2\Delta^3$.}
    \label{fig:upsilon-comparison}
\end{figure}

In \S~\ref{sec:torusThetaGraphHomology} we give a complete description of the graded tangle Floer complex for $\Theta_n$-graphs embedded in the Heegaard torus of $S^3$. 
Figure \ref{fig:upsilon-comparison} displays plots of $-\tfrac{1}{2}\Upsilon$ for two $\Theta_3$-graphs embedded in the Heegaard torus of $S^3$. 
While this description is sufficient to obtain $\Upsilon$ of such objects, we obtain in \S~\ref{sec:twistingInequality} 
a recursive formula which generalizes the formula for $\Upsilon_{T(p,q)}$ described by Feller--Krcatovich \cite{FK17}.

\subsection{Graded maps for cobordisms of spatial graphs}

\begin{definition}
Suppose $(M, T )$ 
and $(M', T')$ are balanced tangles corresponding
to embedded graphs with the same number of vertices.
A stable cobordism from 
$(M, T )$ to $(M', T')$ is a pair $(W, F )$
satisfying the following:
\begin{enumerate}
\item $W$ is a smooth, oriented four-dimensional cobordism with boundary and corners from $M$ to $M'$. Its boundary decomposes as a union 
$$\partial W = \partial_v W \cup \partial_h W = (-M \sqcup M') \cup ( \partial M \times I)$$
with corner stratum $\partial_h W \cap \partial_v W = \partial M \sqcup \partial M'  $.
\item Similarly, $F$ is an oriented, properly embedded surface with boundary and
corners in $W$ which gives a cobordism from $T$ to $T'$ 
such that $\partial_v F = (-T) \sqcup T'$,
$\partial_h F \cap  \partial_h W = [0, 1] \times \partial T = [0,1] \times \partial T'$ under the above identifications, and
$\partial_h F \cap \partial_v F = \partial T \cup \partial T'$.
\item The orientation induced from $F$ on $\partial F$ is the orientation inherited from
$-T \sqcup T'$.
\item (Stability) Every connected component of $F$ which is not a disk intersects both
$T$ and $T'$ in more than one connected component.
\end{enumerate}
\end{definition}

As with the link Floer homology, functoriality of the graded invariant of a tangle requires decorations on the connected components of a stable cobordism. 

\begin{definition}
For a stable cobordism $(W,F) : (M,T) \rightarrow (M',T')$, an $A$-coloring of $(W,F)$ is a map $\mathfrak{c} : \pi_0(F) \rightarrow A$ such that $\mathfrak{c}$ induces $A$-colorings on the boundary tangles $T$ and $T'$.
\end{definition}

Suppose that $(M,T)$ and $(M',T')$ are adapted from  embeddings $G,G':\Theta_n \hookrightarrow S^3$ and hence both of their tangle Floer complexes possess absolute gradings indexed by the $n$ perfect matchings $\mathfrak{M}_{\Theta_n}$. 
A decorated, stable cobordism $(W,\mathcal{F}) : (M,T) \rightarrow (M',T')$ together with a choice $\mathfrak{s} \in \mathrm{Spin}^c(W)$ induces a chain map 
$$F_{W,\mathcal{F},\mathfrak{s}} : \mathrm{CFT}(M,T) \rightarrow \mathrm{CFT}(M',T')$$
which is absolutely $gr_{\mathfrak{m}_i}$-graded where $gr_{\mathfrak{m}_i}$ is the absolute grading associated to the matching $\mathfrak{m}_i \in \mathfrak{M}_{\Theta_n}$.
One should observe that any decorated, stable cobordism between $\Theta_n$-graphs $G$ and $G'$ induces decorated cobordisms between their cycles.
The conventions of \cite{AE15, AE20} are such that the grading shift for a given matching $\mathfrak{m} \in \mathfrak{M}_{\Gamma}$ will agree with the associated grading shift of the map coming from the knot Floer homology.
A straightforward adaptation of \cite[Theorem 1.4]{Zem19a} gives the following.

\begin{theorem}\label{thm:ZemkeGradingShifts}

Let $G_0,G_1 : \Theta_n \hookrightarrow  S^3$ and fix an ordering on the edges of $\Theta_n$, inducing an ordering on the matchings $\{\mathfrak{m}_i\}_{i=1}^n = \mathfrak{M}_{\Theta_n}$.
Let $(W,\mathcal{F}) : (S^3,\mathbb{G}_0) \rightarrow (S^3, \mathbb{G}_1)$ be a colored stable cobordism of graphs such that components of $\mathcal{F}$ are indexed by perfect matchings. Fix $\mathfrak{s} \in \mathrm{Spin}^c(W)$, and use the basepoints $\mathfrak{m}_1$ to index terms in the cobordism map by $\mathrm{Spin}^c$ structures.

\begin{enumerate}
\item The map $F_{W,\mathcal{F},\mathfrak{s}}$ is graded with respect to $gr_{\mathfrak{m}_1}$ and satisfies 
$$gr_{\mathfrak{m}_1}(F_{W,\mathcal{F},\mathfrak{s}}(\mathbf{x})) - gr_{\mathfrak{m}_{1}}(\mathbf{x}) = \frac{c_1(\mathfrak{s})^2-2\chi(W) - 3\sigma(W)}{4} + \widetilde{\chi}(\Sigma_{\mathfrak{m}_1}).$$
\item  For $i \neq 1$, the map $F_{W,\mathcal{F},\mathfrak{s}}$ is graded with respect to $gr_{\mathfrak{m}_i}$ and satisfies 
$$gr_{\mathfrak{m}_i}(F_{W,\mathcal{F},\mathfrak{s}}(\mathbf{x})) - gr_{\mathfrak{m}_i}(\mathbf{x}) = \frac{c_1(\mathfrak{s} - PD[\Sigma_{\mathfrak{m}_1} \cup -\Sigma_{\mathfrak{m}_i}])^2-2\chi(W) - 3\sigma(W)}{4} + \widetilde{\chi}(\Sigma_{\mathfrak{m}_i}).$$
where $PD[-]$ denotes the Poincar\'{e}--Lefschetz dual.
\end{enumerate}
\end{theorem}

\begin{remark}
Recall $S^3$ has a unique $\mathrm{Spin}^c$ structure which is torsion. For cobordisms $W$ with boundary components $Y_i$, the theorem above requires additional hypotheses that $\mathfrak{s}|_{Y_i}$ and $(\mathfrak{s}-PD[\Sigma_{\mathfrak{m}_1} \cup -\Sigma_{\mathfrak{m}_i}])|_{Y_i}$ are torsion for $Y_i$. 
\end{remark}

\subsection{Edge and vertex stabilizations of spatial graphs}\label{sec:stabilizations}

\begin{definition}[Edge stabilization]\label{def:edgeStabilization}
Let ${ \Gamma }$ be an embedded graph in a closed, orientable 3-manifold ${ Y }$. Selecting an edge ${ e \in \Gamma }$, we may perform an \textit{edge stabilization} of ${ \Gamma }$ along ${ e }$ to obtain a graph differing from ${ \Gamma }$ by replacing a single edge with three edges connected by two new vertices.

On the other hand, if a graph  contains two valence-two vertices connected by an edge, we can create a new graph by removing the vertices and merging edges. We call this latter operation an \textit{edge destabilization}.
\end{definition}

\begin{figure}[htbp]
\newcommand{\graphpic}[1]{%
    \vcenter{\hbox{#1}}%
}

\newcommand{\stabpair}[1]{%
    \vcenter{\hbox{$\displaystyle
        \begin{gathered}
            \xrightarrow{
                \substack{
                    \text{\Large #1}\\[-0.25ex]
                    \text{\Large stabilization}
                }
            }
            \\[3.4ex]
            \xleftarrow{
                \substack{
                    \text{\Large #1}\\[-0.25ex]
                    \text{\Large destabilization}
                }
            }
        \end{gathered}
    $}}%
}

\begin{equation*}
    \graphpic{\includegraphics[width=0.13\linewidth]{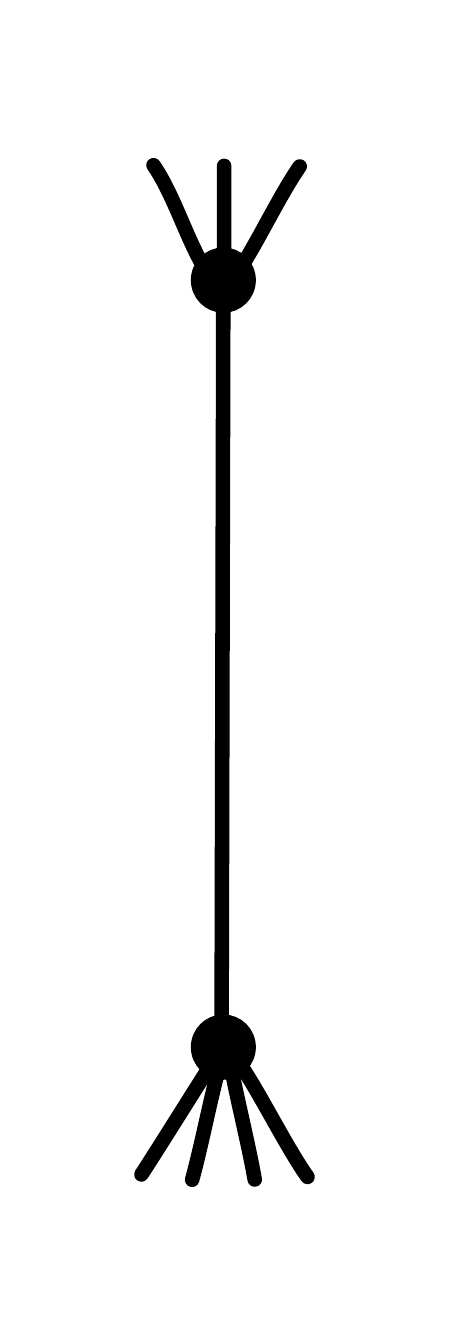}}
    \hspace{1.2em}
    \stabpair{edge}
    \hspace{1.2em}
    \graphpic{\includegraphics[width=0.13\linewidth]{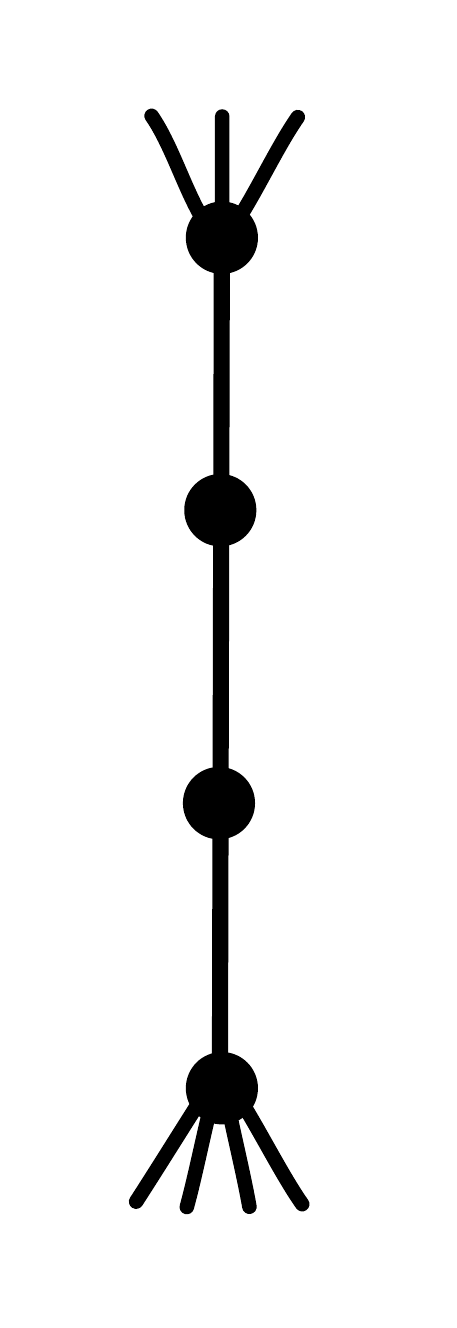}}
    \hspace{2em}
    \graphpic{\includegraphics[width=0.13\linewidth]{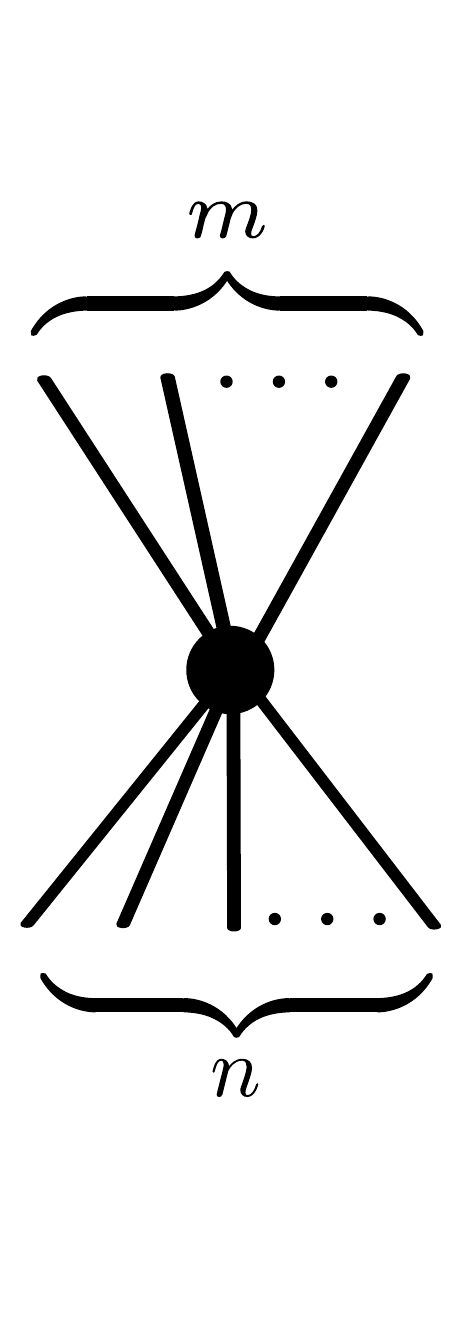}}
    \hspace{1.2em}
    \stabpair{vertex}
    \hspace{1.2em}
    \graphpic{\includegraphics[width=0.13\linewidth]{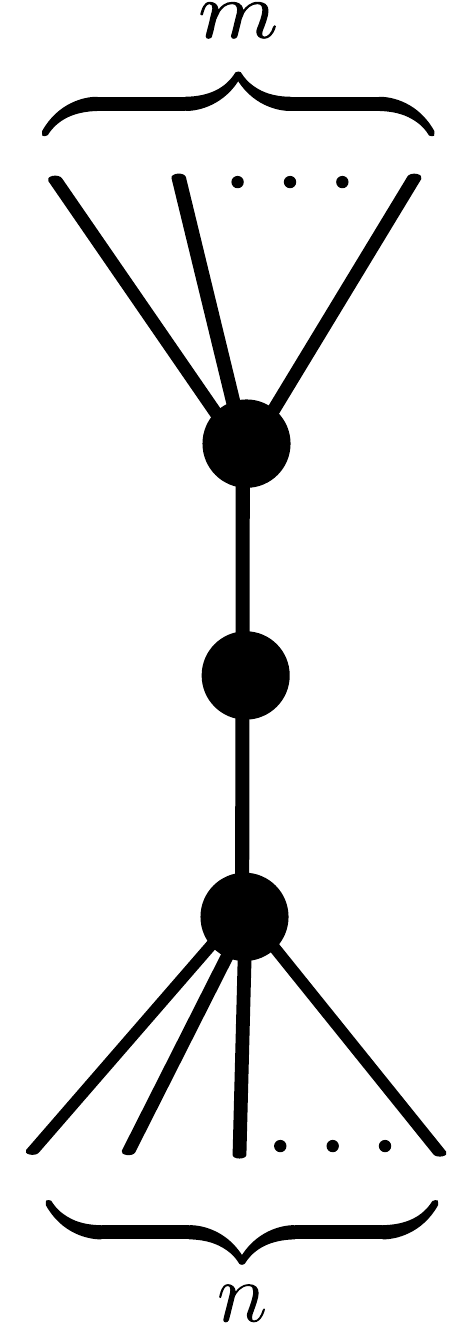}}
\end{equation*}
    \caption{
    Two (de)stabilization operations on abstract graphs which preserve the properties of being bipartite and balanced. There is also a bijective correspondence between perfect matchings in the initial graph and its (de)stabilization.}
    \label{fig:edgeAndVertexStabilizations}
\end{figure}

\begin{definition}[Vertex stabilization]\label{defn:vertexStab}
Let $\Gamma_{0}$ be an abstract graph containing an ${ (m+n) }$-valent vertex $v$. An ${ (m,n) }$-vertex stabilization of $\Gamma_{0}$ at ${ v }$ is a graph ${ \Gamma_{1} }$ obtained by replacing $v$ with a configuration consisting of one ${ (m+1) }$-valent vertex, one ${ (n+1) }$-valent vertex, and one bivalent vertex in the manner depicted in the right-hand side of Figure \ref{fig:edgeAndVertexStabilizations}.  We call the operation of going backwards a \textit{vertex destabilization}. 
\end{definition}

\begin{remark}\label{rem:edgeStabilizationFromVertexStabilization}
In the link Floer setting, the operation called ``quasistabilization'' in \cite{Zem17} and \cite{Zem19b} is equivalent to what we have chosen to call an edge stabilization.
The attentive reader might notice that an edge stabilization along an edge between two vertices $p$ and $q$ of valence $m$ and $n$ respectively may equivalently be described as either an $(m-1,1)$-vertex stabilization at $p$ or an $(n-1,1)$-vertex stabilization at $q$. We find it useful, however, to treat edge stabilization as a different kind of operation.
\end{remark}

\begin{lemma}\label{lem:matchingBijectionEdgeStab}
Let $\Gamma$ be a balanced, bipartite graph admitting perfect matchings. If ${ \widetilde{\Gamma}}$ differs from $\Gamma$  by an edge or vertex (de)stabilization, then there is a bijection between the perfect matchings for ${ \Gamma }$ and those of ${ \widetilde{\Gamma} }$.
\begin{proof}

Per Remark \ref{rem:edgeStabilizationFromVertexStabilization}, it suffices to prove the result for vertex stabilizations.
If ${ \widetilde{\Gamma} }$ is obtained by applying an ${ (m,n) }$-vertex stabilization to ${ \Gamma }$ at the vertex ${ v }$,
then as suggested by Figure \ref{fig:edgeAndVertexStabilizations}, we let ${ e_{1}, \dots, e_{m} }$ enumerate the first ${ m }$ edges of ${ \Gamma }$ about ${ v }$ and let ${ e_{m+1},\dots, e_{m+n} }$ enumerate the final ${ n }$ edges about ${ v }$. 
Each of these edges has an obvious counterpart in ${ \widetilde{\Gamma} }$. We let ${ e' }$ and ${ e'' }$ denote the two new edges of ${ \widetilde{\Gamma} }$ which have no counterpart in ${ \Gamma }$ so that ${ e' }$ shares a vertex with each element of ${ \left\{e_{i}\right\}_{i=1}^{m} }$ while ${ e'' }$ shares a vertex with each element of ${ \left\{e_{i}\right\}_{i=m+1}^{m+n} }$.
Note that every perfect matching ${\mathfrak{m} }$ of ${\Gamma }$ must contain exactly one of the edges ${ e_{1}, \cdots, e_{m+n} }$. We obtain a bijection between matchings for ${ \Gamma }$ and those for ${ \widetilde{\Gamma}}$ whose forward assignment is
$$ \mathfrak{m} \mapsto\begin{cases} \mathfrak{m} \cup \{e''\}  & \left\{e_{i} \right\}_{i=1}^{m}   \cap \mathfrak{m} \neq \emptyset \\ \mathfrak{m} \cup \{e'\} & \left\{e_{i} \right\}_{i=m+1}^{m+n} \cap \mathfrak{m} \neq \emptyset   \end{cases} .$$
\end{proof}
\end{lemma}

We conclude with a discussion of how these operations affect the homotopy type of the curved tangle Floer complexes. 
Recall that for a spatial graph $G : \Gamma \hookrightarrow S^3$, $\mathrm{CFT}(G)$ denotes the homotopy type of curved complexes over the polynomial ring $R_{E(\Gamma)} = \mathbb{F}[Z_e]_{e \in E(\Gamma)}$.

\begin{theorem}\label{thm:edgeStabilization}
Let ${ G : \Gamma \hookrightarrow Y^3 }$ be a spatial graph and suppose ${ \widetilde{G} }$ is obtained from ${ G }$ by an edge stabilization operation at an edge $e$ which connects to a positive vertex $v$. Let $\{e_i\}_{i=0}^{m}$ enumerate the edges of $\Gamma$ which are adjacent to $v$ with $e_0=e$. 
Let $w$ and $z$ denote two new edges created by the edge stabilization and use the label $e_0$ for the new edge which is adjacent to $v$.
Then we have a homotopy equivalence of curved complexes
$$
\begin{aligned}
\mathrm{CFT}(\widetilde{G}) \simeq \mathrm{CFT}(G)\otimes \left[
\begin{tikzcd}[ampersand replacement=\&, row sep=0pt, column sep=2.2em, baseline=(current bounding box.center)]
\Theta^+ \& \& \Theta^-
\arrow[bend left=15, from=1-1, to=1-3, "{ P_1 }"]
\arrow[bend left=15, from=1-3, to=1-1, "{ P_2 }"]
\end{tikzcd}
\right] & \qquad \qquad  \begin{matrix} P_1 = Z_{z} + \prod\limits_{i=1}^{m} Z_{e_i} \\ P_2 = Z_{e_0} + Z_{w} \phantom{\prod\limits_{i=1}^{m} }  \end{matrix}
\end{aligned} .
$$

\begin{proof}
Figure \ref{fig:edgeStabilizationHeegaardDiagram} displays two local Heegaard diagrams $\mathcal{H}$ and $\widetilde{\mathcal{H}}$ which respectively describe the original spatial graph $G$ and its edge stabilization $\widetilde{G}$. In both figures, the basepoint is given the name of the edge it represents.

\begin{figure}[!h]
\captionsetup{
        justification=raggedright,
        singlelinecheck=false,
        labelsep=newline
    }
\newcommand{\graphpic}[1]{%
    \vcenter{\hbox{#1}}%
}

\newcommand{\stabpair}[1]{%
    \vcenter{\hbox{$\displaystyle
        \begin{gathered}
            \xrightarrow{
                \substack{
                    \text{\Large #1}\\[-0.25ex]
                    \text{\Large stabilization}
                }
            }
            \\[3.4ex]
            \xleftarrow{
                \substack{
                    \text{\Large #1}\\[-0.25ex]
                    \text{\Large destabilization}
                }
            }
        \end{gathered}
    $}}%
}

\begin{equation*}
\graphpic{\includegraphics[width=0.25\linewidth]{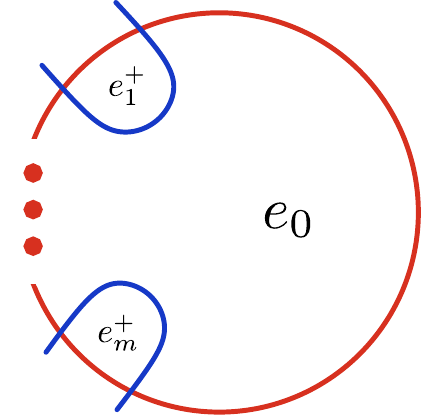}}
    \hspace{1.2em}
    \stabpair{edge}
    \hspace{1.2em}
    \graphpic{\includegraphics[width=0.25\linewidth]{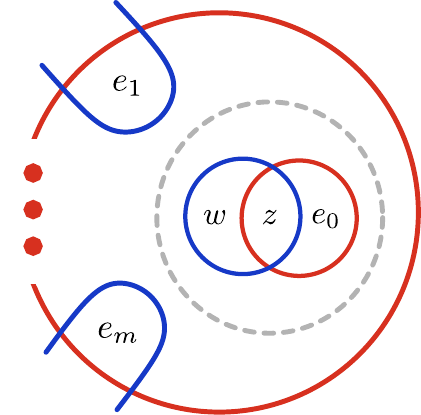}}
\end{equation*}
\caption{
\textbf{Left:} a local Heegaard diagram $\mathcal{H}$ describing all half-edges meeting a vertex $v$ of valence $m+1$. The edges $\{e_i\}_{i=1}^m$ project to arcs which traverse the $\beta$ curves. The edge $e_0$ we are stabilizing along projects to an arc that traverses the lone $\alpha$ curve. \newline
\textbf{Right:} a local Heegaard diagram $\widetilde{\mathcal{H}}$ describing the result of an edge stabilization along $e_0$. Two new vertices are introduced to the graph, corresponding to the new $\alpha$ and $\beta$ curves. We also introduce two basepoints, corresponding to new edges $w$ and $z$.}\label{fig:edgeStabilizationHeegaardDiagram}
\end{figure}

An argument identical to that presented in \S~4 of \cite{Zem19b} tells us that one may obtain the desired model for $\mathrm{CFT}(\widetilde{G})$ by choosing an almost complex structure which sufficiently ``stretches the neck'' at the dashed curve $\gamma$ seen in the depiction of $\widetilde{\mathcal{H}}$ displayed in Figure \ref{fig:edgeStabilizationHeegaardDiagram}. 
The only difference is that in Zemke's original argument the $\alpha$-boundary degeneration which traverses the annular neighborhood of $\gamma$ only meets a single basepoint whereas, in our setting, it meets the basepoints $e_1,\cdots , e_m$.
\end{proof}
\end{theorem}

\begin{theorem}\label{thm:vertexStabilizedFloerComplex}
Let ${ G : \Gamma \hookrightarrow S^3 }$ be a spatial graph and suppose ${ \widetilde{G} }$ is obtained from ${ G }$ by an $(m,n)$-vertex stabilization operation at a vertex $v$ where edges incident to $v$ are partitioned as 
$\{e_{i}^{+}\}_{i=1}^m$ and $\{e_{i}^-\}_{i=1}^{n}$. Let $e'$ and $e''$ denote the two new edges inserted during the vertex stabilization. Then we have a homotopy equivalence of curved complexes
$$
\begin{aligned}
\mathrm{CFT}(\widetilde{G}) \simeq \mathrm{CFT}(G)\otimes \left[
\begin{tikzcd}[ampersand replacement=\&, row sep=0pt, column sep=2.2em, baseline=(current bounding box.center)]
\Theta^+ \& \& \Theta^-
\arrow[bend left=15, from=1-1, to=1-3, "{ P_1 }"]
\arrow[bend left=15, from=1-3, to=1-1, "{ P_2 }"]
\end{tikzcd}
\right] 
& \qquad \qquad \begin{matrix}P_1 = Z_{e''} + \prod\limits_{i=1}^{m} Z_{e_{i}^+} \\ P_2 = Z_{e'} + \prod\limits_{i=1}^n Z_{e_{i}^-} 
\end{matrix}
\end{aligned} .
$$
\begin{proof}
We hold off on the full proof of Theorem \ref{thm:vertexStabilizedFloerComplex} until \S
\ref{sec:gridHomologyForSpatialGraphs} where we prove a technical lemma (Lemma \ref{lem:generalizedStabilization}), 
but we provide a sketch here. 

Theorem \ref{thm:edgeStabilization} specifies a manner in which $\mathrm{CFT}(G)$ may be recovered from $\mathrm{CFT}$ of its edge stabilizations. 
Given any spatial graph $G$, one may apply some number of edge stabilizations to obtain a spatial graph $G'$ which admits a Heegaard diagram resembling the grid diagrams of Ozsv\'{a}th, Stipsicz, and Szab\'{o} \cite{OSS15} --- identical Heegaard diagrams for spatial graphs are studied in \cite{HO17, Kub23, Kub25a, Kub25b}.
In this setting we can describe the effect of a vertex stabilization as a concrete transformation of grid diagrams and prove Lemma \ref{lem:generalizedStabilization} which is the restriction of Theorem \ref{thm:vertexStabilizedFloerComplex} to the grid setting.
Combining this with Theorem \ref{thm:edgeStabilization} completes the proof. 
\end{proof}

\end{theorem}
\begin{corollary}\label{cor:StabilizationHtpy}
    If $\widetilde{G}$ is obtained from $G$ by either a vertex or an edge stabilization, then for an appropriate coloring $\mathfrak{c}$, $\mathrm{CFT}(\widetilde{G},\mathfrak{c})$ will be homotopy equivalent to $\mathrm{CFT}(G)$, linearly in the variables corresponding to edges of $G$ but not necessarily the new edges. If $\mathfrak{m}$ is any perfect matching of $G$, and $\widetilde{\mathfrak{m}}$ is the corresponding matching under Lemma~\ref{lem:matchingBijectionEdgeStab}, then this homotopy equivalence intertwines the gradings $gr_\mathfrak{m}$ and $gr_{\widetilde{\mathfrak{m}}}$.
\end{corollary}
\begin{proof}
    Color the edges of $\widetilde{G}$ already present in $G$ with the corresponding variables $Z_i$, color one of the two new edges with an independent variable $U$, and color the other edge with an appropriate monomial in the $Z_i$ so as to make one of the terms $P_1$ or $P_2$ in the other tensor factor from Theorem~\ref{thm:edgeStabilization} or \ref{thm:vertexStabilizedFloerComplex} equal to $0$. Then the new tensor factor becomes a Koszul complex imposing the other of $P_1$ or $P_2$ as a relation involving $U$; the new variable $U$ may then be eliminated via a homotopy equivalence recovering the original $\mathrm{CFT}(G)$.

    For the statement about gradings, the homotopy equivalences of Theorems~\ref{thm:edgeStabilization} and \ref{thm:vertexStabilizedFloerComplex} are easily seen to preserve all of the relative gradings. The absolute lift to form $gr_{\mathfrak{m}}$ then depends only on $\mathrm{CFT}(G)$ as a relatively graded complex, and is not an extra piece of structure: one sets all variables outside $\mathfrak{m}$ to $1$ and puts the resulting tower generator in degree $0$. In particular, since $\mathrm{CFT}(G)$ and $\mathrm{CFT}(\widetilde{G},\mathfrak{c})$ are homotopy equivalent, this procedure will play out identically for each, and $gr_{\mathfrak{m}}$ and $gr_{\widetilde{\mathfrak{m}}}$ will agree.
\end{proof}

\subsection{Square gradings for the tangle Floer complex}\label{sec:squareGradings}
Throughout, we assume that $\Gamma$ is an oriented, bipartite, and balanced graph.
We say that an orientation of ${ \Gamma }$ is bipartite if all edges are oriented from one of the bipartition sets to the other.
Recall that a \textit{circuit} of a graph ${ \Gamma }$ is a non-repeating sequence of edges creating a closed loop--- we remark that although edges do not repeat in a circuit, vertices can. Notice that for a graph equipped with a bipartite orientation, all circuits have even length and, as one travels around the circuit, the orientations of the edges agree and disagree with the direction of travel in an alternating pattern.

\begin{definition}\label{defn:bounding_a_square}
Let ${ G: \Gamma \hookrightarrow S^{3} }$ be a spatial graph. We say that a length-4 circuit of ${ G }$ bounds a square if there is an embedding ${\phi :  I \times I \rightarrow S^{3} \setminus \nu \Gamma }$ such that
 ${ \phi|_{\partial}}$ is the image under ${ G }$ of a push-off of a square circuit ${\square \subset  \Gamma }$. After extending $\phi$ into the interior of $\nu\Gamma$, we require each element of $\{0,1\}^2$ to map to a vertex (not necessarily distinct from the images of the other elements), and the boundary arcs to map to edges. 
\end{definition}

\begin{example}
The above can be motivated by the study of nonorientable knot cobordisms.
Given an oriented knot ${ K \subset Y^{3} }$ with distinct, marked points ${ p }$ and ${ q }$, let ${ \beta : I \times I \hookrightarrow Y^{3} }$ be a nonorientable band whose attaching region in $K$ is disjoint from the marked points and joins the two components of $K \setminus \{p,q\}$. Let ${ v_{1},v_{2}, v_{3},v_{4} }$ denote the images of ${\{0,1\}^2 \subset  I \times I }$ under ${ \beta }$.
Consider the spatial graph ${ G(K,\beta) }$ whose vertex set is ${V= \left\{p,q,v_{1},v_{2},v_{3},v_{4}\right\} }$ and whose edges are the connected components of 
$$ (K \cup \mathrm{Im}(\beta|_{\partial})) \setminus V   $$ 
${ G(K,\beta) }$ is an embedding of a subdivision of ${ K_{4} }$. For an appropriate bipartite orientation of its edges, the map ${ \beta }$ reveals that
${ G(K,\beta) }$ bounds a square whose corresponding circuit is the image of ${ \beta|_{\partial} }$. Taking one of the examples ${ G(K,\beta) }$ described above, one may apply two vertex destabilizations to it to obtain a spatial ${ \Theta_{4} }$-graph in ${ Y^{3} }$ which also bounds a square.
\end{example}

\begin{definition}[Square grading]\label{defn:squareGrading}
Let ${ G : \Gamma \hookrightarrow S^{3} }$ be a spatial graph and suppose ${ \square = (e_{1}, e_{2},e_{3},e_{4}) }$ is a square circuit of ${ \Gamma }$.
The square grading $gr_{\square}$ associated to the circuit $\square$ is the relative grading of $\mathrm{CFT}(G)$ defined by the grading datum
$$\rho_{\square} \in \mathrm{Hom}(\mathcal{G}_{\Gamma},\mathbb{Z}) \qquad \rho_{\square}(1,0) = 0 \qquad \rho_{\square}(e) = \begin{cases} 2 & e \in \{e_1,e_3\} \\ -2 & e \in \{e_2,e_4 \} \\ 0 & \text{else} \end{cases}.$$

\end{definition}

\begin{lemma}\label{lem:squareGradingConstrained}
Let $\mathcal{H}$ be a multipointed Heegaard diagram for some embedded graph $\Gamma \hookrightarrow S^3$. Suppose $\mathcal{H}$ contains a region where basepoints occupy four adjacent regions meeting at a single $\bm{\alpha\beta}$-intersection point $p$. Then there is a Heegaard state $\bm{x}$ such that for any other Heegaard state $\bm{y}$, we have
$$gr_{\square}(\bm{x})- gr_{\square}(\bm{y})\in \{0, 2 \}.$$

\begin{proof}
Consider the local picture in $\mathcal{H}$ about $p$, Figure \ref{fig:localSquareRegion}. 

\begin{figure}[!h]
\includegraphics[width = 0.25\linewidth]{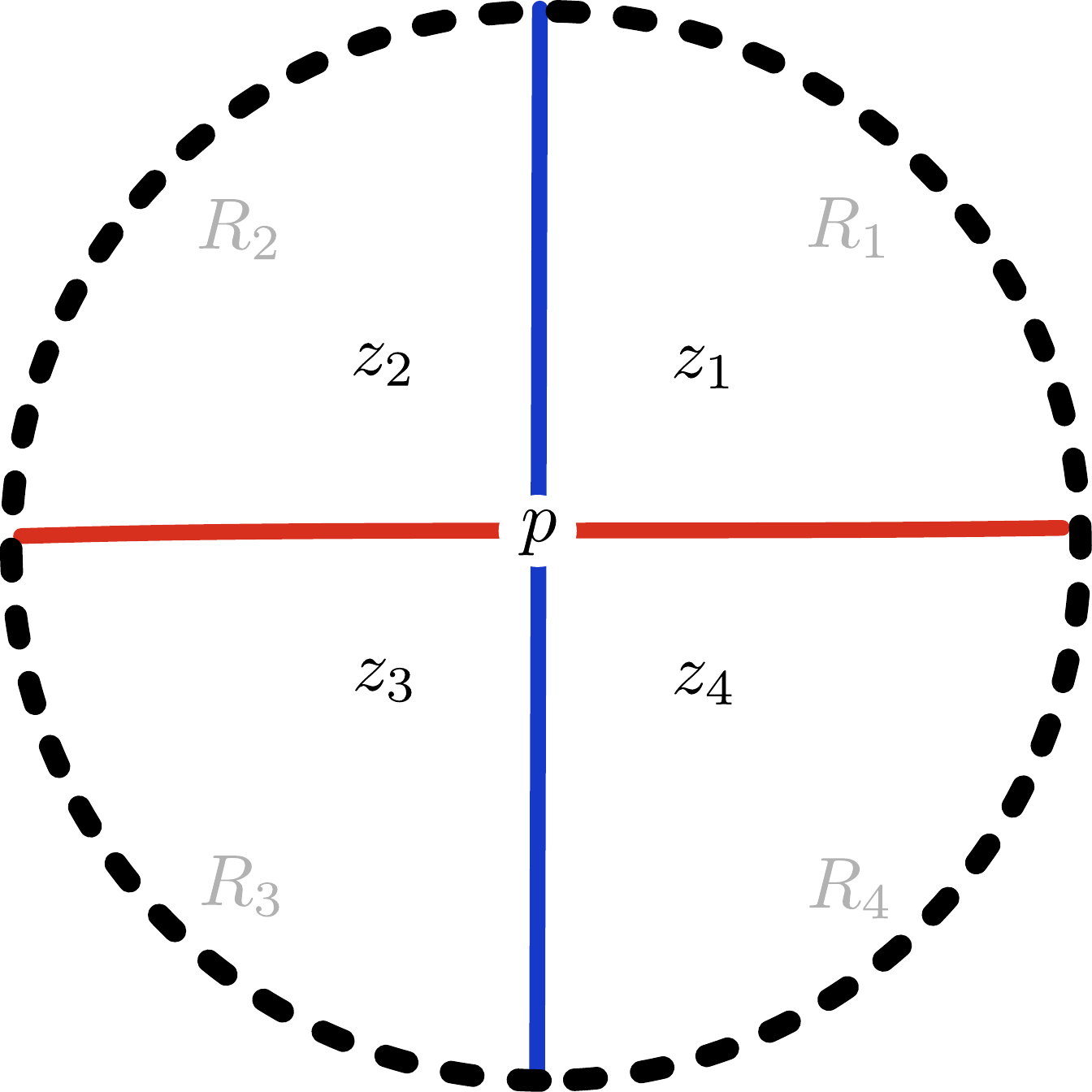}
\caption{
}\label{fig:localSquareRegion}
\end{figure}
Given $\bm{x},\bm{y} \in \mathrm{CFT}(\mathcal{H})$, let $\phi \in \pi_2(\bm{x},\bm{y})$. 
Letting $z_i$ denote the basepoint associated to the edge $e_i$, it is immediate from Definition \ref{defn:squareGrading} and general properties of these relative gradings that
\begin{equation}\label{eqn:alt_grading_difference}
gr_{\square}(\bm{x}) - gr_{\square}(\bm{y}) = -2n_{z_1}(\phi) + 2n_{z_2}(\phi) - 2n_{z_3}(\phi) + 2n_{z_4}(\phi).   
\end{equation}
The Heegaard state $\bm{x} \in \mathbb{T}_{\alpha} \cap \mathbb{T}_{\beta}$ either does or does not contain $p$, partitioning the Heegaard states into two sets.
Recall a well-known constraint on the local multiplicities of $\phi$ about $p$:
$$
n_{z_1}(\phi) - n_{z_2}(\phi) + n_{z_3}(\phi) - n_{z_4}(\phi) = \begin{cases}
1 & p \in \bm{x} \setminus \bm{y} \\ -1 & p \in \bm{y} \setminus \bm{x} \\ 0 & \text{else}
\end{cases}.
$$
 If only one of these sets is non-empty, then the above tells us the relative square grading difference is zero. Otherwise, we may choose $\bm{x}$ such that $p \notin \bm{x}$ and conclude for any  $\bm{y}$ we have
$$gr_{\square}(\bm{x})-gr_{\square}(\bm{y})=\begin{cases} 2 & p \in \bm{y} \\ 0 & \text{else} \end{cases}.$$ 
\end{proof}
\end{lemma}

\begin{remark}\label{rem:squareGrading} 
Observe that any graph admitting a Heegaard diagram with a local picture resembling Figure \ref{fig:localSquareRegion} will bound a square: the four edges bounding the square correspond to the four basepoints distributed about $p$ and the square itself projects to their convex hull.

In the case where $G$ happens to be an embedding of a ${ \Theta_{n} }$-graph (or a suitable stabilization thereof) which bounds a square $\square$, then the relative grading $gr_{\square}$ lifts to an absolute grading. In this case, we can explicitly determine the support of the absolute grading $gr_{\square}$; see \S~\ref{sec:squareGradingComputation}. A necessary step in our argument is that if a suitable stabilization of the $\Theta_n$-graph bounds a square, then there is a Heegaard diagram of exactly the form seen in Figure \ref{fig:localSquareRegion}.
\end{remark}

\section{Surfaces in $B^4 \setminus (S^1 \times B^3)^{\circ}$ and foams} \label{sec:FoamPrelims}

Theorem~\ref{thm:FoamBound1} is a special case (looking only at the endpoint cases of $\mathfrak{c}(t)$ where one edge becomes $0$-colored and may be removed) of the following result, which in turn follows directly from bounds we obtain on the square grading (Corollary~\ref{cor: squareBoundSumUp}):
\begin{proposition}\label{prop:squareUpsilonSlope}
    Suppose $\Gamma$ is a knotted $\Theta_n$-graph, and that some length-4 circuit of its edges $(e_1,e_2,e_3,e_4)$ bounds a square. Fix some coloring $\mathfrak{c}$ of $\Gamma$, and let $\mathfrak{c}(t)$ be the one-parameter family of colorings in which the colors on $e_1$ and $e_3$ are increased by $t$, and colors on $e_2$ and $e_4$ are decreased by $t$. Observe that $\mathfrak{c}(t)$ is defined for $t$ in a closed interval which includes $0$ and, at the endpoints of this interval, some edge is colored zero. Let $b$ be the single-band unoriented cobordism from $e_2 \cup e_4$ to $e_1 \cup e_3$; then $\Upsilon_\Gamma(\mathfrak{c}(t))$ is a piecewise linear function of $t$, with all slopes within $\frac{1}{2}$ of $\frac{1}{4}e(b)$.
\end{proposition}
This should be seen as directly analogous to the fact that the slopes of the $\Upsilon$ invariant of a knot are bounded by the Seifert genus --- indeed, both results are proven by showing that a certain grading on the Floer homology of the graph is supported in a bounded interval, only here the interval is not centered at $0$. 

Although this result is entirely 3-dimensional, its consequences have a 4-dimensional flavor, especially combined with the elementary facts that the following two operations (or their reverses) do not change $\Upsilon$:
\begin{itemize}
    \item Removing a $0$-colored edge
    \item Combining two parallel edges and adding the colors.
\end{itemize}
We will generally call consequences of Proposition~\ref{prop:squareUpsilonSlope} together with these two facts (in particular, bounds on $\Upsilon_\Gamma$ coming from a sequence of applications of Theorem~\ref{thm:FoamBound1}) ``foam inequalities'', since they generally may be interpreted as coming from particular singular surfaces. For example, splitting an edge into two parallel copies corresponds to the birth or collapse of a bubble (a local maximum or minimum along a codimension-1 singularity in the foam), and the special case of Proposition~\ref{prop:squareUpsilonSlope} in which two opposite edges of the square are initially colored $0$ may be seen as a 1-handle attachment between two different edges of the graph. In fact, the operation of varying from $\mathfrak{c}$ to $\mathfrak{c}(t)$ may always be factored into two edge splits, a saddle, and two edge merges. We shall also need the following much more direct analogue of the Seifert genus bound from knots, proven in exactly the same way:
\begin{lemma}\label{lem:threegenusbound}
    Suppose $\Gamma$ is a knotted $\Theta_n$-graph, and that some two-cycle of its edges $e_1,e_2$ bounds  a genus-$g$ surface $\Sigma \subset S^3$ with interior disjoint from $\Gamma$. Fix some coloring $\mathfrak{c}$ of $\Gamma$, and let $\mathfrak{c}(t)$ be the one-parameter family of colorings in which the color on $e_1$ is increased by $t$ and the color on $e_2$ is decreased by $t$. Then all slopes of $\Upsilon_\Gamma(\mathfrak{c}(t))$ are in between $-g$ and $g$.
\end{lemma}

One other special case deserves attention here: in Theorem~\ref{thm:FoamBound1}, if we ignore the orientations and colors on the edges, the graphs have changed as if they were level sets just before and after passing through a tetrahedral codimension-2 singularity in a foam. A foam with such tetrahedral singularities, provided that all of the regular level sets for some Morse function are $\Theta_3$-graphs, gives bounds on $\Upsilon$ at its boundaries, provided that one tracks colors through in this nonlocal way. In \S~\ref{sec:torusThetaWebs}, we will work through many examples of this type, and see that this nonlocality is what ultimately makes our bounds piecewise linear in the coloring and allows them to faithfully reconstruct $\Upsilon$.

Although the bounds from Proposition~\ref{prop:squareUpsilonSlope} may be interpreted as coming from foams, we do not have a fully satisfying way to go in the other direction: in order to read off a bound on $\Upsilon$ from a foam, we generally need some nonlocal information on e.g. the level sets of a Morse function. This comes from two restrictions on the graphs to which our results apply: first, the tangle Floer theories we use always require a bipartite input. This is mostly a bookkeeping issue, but does make the tetrahedral vertex slightly awkward, since the tetrahedron is not bipartite. More seriously, Proposition~\ref{prop:squareUpsilonSlope} requires that $\Gamma$ be a $\Theta_n$-graph. Ultimately, this restriction comes from the need to precisely pin down an absolute grading: for general $\Gamma$, the proof of Proposition~\ref{prop:squareUpsilonSlope} pins all of the slopes to some interval of length $1$, but we cannot specify the exact interval from topological data (e.g. the Euler number in the statement). Below, we will give a situation in which our method does give a reasonably general statement; here, one term in the bound is something like an Euler number of the foam, and this is precisely what we lack a definition for in general.

\begin{theorem}\label{thm:OneBindingCircle}
    Suppose that $\Sigma$ is a properly embedded framed surface inside of $S^2 \times D^2$, possibly with multiple balls removed, such that the portion of $\partial\Sigma$ on $S^2 \times S^1 = \partial (S^2\times D^2)$ is the coherently oriented closure of a framed spherical braid $\beta$ with framing induced from $\Sigma$, each component of $\Sigma$ meets $S^2  \times S^1$, and the portion of $\partial\Sigma$ on the boundary of each removed ball is a knot (with the zero framing induced from $\Sigma$). Then, for any assignment of a color $t_i \in [0,2]$ to components $\Sigma_i$ of $\Sigma$, such that $\sum_i [\partial \Sigma_i] t_i = 2$ (identifying ${[\partial \Sigma_i] \in H_1(S^2 \times S^1) \cong \mathbb{Z}}$ with the integer it represents), we have $$\left|\sum_j \Upsilon_{K_j}(t_{i_j}) + \frac{1}{4}w_{t}(\beta)\right| \leq \frac{1}{2}\sum_i t_i \left( [\partial\Sigma_i] - \chi(\hat{\Sigma_i}) \right)$$ where the $K_j$ are the boundary knots, $i_j$ indexes the component of $\Sigma$ meeting $K_j$, and $\hat{\Sigma_i}$ is the result of filling in all boundary components coming from the $K_j$.
\end{theorem}

Here, $w_t(\beta)$ is a measure of the writhe of the framed, spherical braid $\beta$, weighted by the coloring $t$, and is our analogue of the Euler number for $\Sigma$. To define it, notice that $\beta$ is represented by a conjugacy class in the framed, spherical braid group on $n = [\partial \Sigma]$ many strands, which we may present with generators $\sigma_1,\dots,\sigma_{n-1}$ (crossings) and $\tau_1,\dots,\tau_n$ (framing twists) and relations:
\begin{align*}
    \sigma_i \sigma_{i+1} \sigma_i = \sigma_{i+1}\sigma_i \sigma_{i+1} && [\sigma_i, \sigma_j] = 1 \text{ (for } |i-j| > 1 \text{)} && \text{(braid relations)} \\
    \tau_i\sigma_i = \sigma_i \tau_{i+1}, \tau_{i+1} \sigma_i = \sigma_i \tau_i && [\sigma_i,\tau_j] = 1 \text{ (for } j \neq i,i+1 \text{)} && \text{(twist slides)} \\
    && [\tau_i,\tau_j] = 1 \text{ (for } i \neq j \text{)} && \text{(twists commute)} \\
     \tau_1^2 \cdot \sigma_1\sigma_2\cdots\sigma_{n-1}\sigma_{n-1} \cdots\sigma_2\sigma_1 = 1 && && \text{(wrapping about } S^2 \text{)}
\end{align*}
\begin{definition}
    $$w_t(\beta) = \left(\sum_{\text{twists }\tau\text{ on strand }i} \operatorname{sign}(\tau)t_i(2-t_i)\right) - \left( \sum_{\text{crossings }\sigma\text{ on strands }i,j} \operatorname{sign}(\sigma) t_it_j \right).$$
\end{definition}

Properly speaking, this should be treated as a homomorphism from a groupoid of colored, framed, spherical braids to $\mathbb{R}$, to account for how the differently colored strands trace through the braid. This is straightforwardly unchanged under all but the last relation of the framed spherical braid group, and invariance under the last follows from $\sum_i [\partial \Sigma_i] t_i = 2$.

To treat Theorem~\ref{thm:OneBindingCircle} as a statement about foams, collapse the $S^2 \times S^1$ boundary to $\{\text{pt}\} \times S^1$ to create a foam with a single circle of codimension-1 singularities.
\begin{proof}[Proof of Theorem~\ref{thm:OneBindingCircle}]
    We'll start by collapsing the boundary components, leading to a surface in $S^4$ with two kinds of singularities: first, a single circle (the ``binding circle'') where several sheets of $\Sigma$ come together, coming from collapsing $S^2 \times S^1$ to $S^1$. Since $\Sigma$ meets the $S^2 \times S^1$ boundary component in a braid, the boundary of the smooth part of $\Sigma$ maps to the binding circle as an $n$-to-$1$ covering. Second, from collapsing the $S^3$ boundary components to points, we get finitely many points at which $\Sigma$ fails to be locally flat, with the link of each such singularity being a knot.

    Next, delete small balls around two points on the binding circle. The result is a cobordism in $S^3 \times I$ between two unknotted $\Theta_n$-graphs, with remaining singularities consisting of a pair of arcs between two boundary vertices, and the non-locally-flat points.

    Arrange all of this to be Morse with respect to the projection $S^3 \times I \to I$, push local minima before and local maxima after all other events, and then do a series of handleslides to make the Morse events happen in the following order: first, several births of circles, followed by fusion bands joining the new circles to the $\Theta_n$-graph. Then, a series of events of one of the following forms:
    \begin{enumerate}
        \item A pair of bands, all of whose endpoints are on a single edge (so this one edge undergoes a genus-1 cobordism)
        \item A band connecting two different edges (which will therefore belong to the same component of $\Sigma$ and so have the same color)
        \item A non-locally-flat point.
    \end{enumerate}
    Finally, a series of fission bands and deaths of circles. The first and last steps are concordances and do not change $\Upsilon$ by Theorem~1.1 of \cite{Ali22}, so we need only concentrate on moves of type (1), (2), and (3). In the end, we start and end with an unknotted graph, so the total change in $\Upsilon$ over all of these moves is $0$, which will give our bound.

    Moves of type (1) are handled by Lemma~\ref{lem:threegenusbound}, which shows that if the edge undergoing the cobordism is colored $t$, then $\Upsilon$ changes by at most $t$ in absolute value. Moves of type (3) correspond to tying a copy of $K_j$ into a strand colored $t_{i_j}$; this changes $\Upsilon$ by precisely $\Upsilon_{K_j}(t_{i_j})$ by Lemma~1.2 of \cite{Ali22}, since it corresponds exactly to taking a connected sum with the $\Theta_n$-graph formed by starting with a $\Theta_2$-graph representing $K_j$ and replacing one of the two edges with $n-1$ parallel copies.

    Finally, moves of type (2) are handled using Proposition~\ref{prop:squareUpsilonSlope}. Let $\Gamma$ be the graph appearing immediately before the surgery, and let $\Gamma'$ be the result of adding two new $0$-colored edges to $\Gamma$, corresponding to the new edges that will appear after the surgery. If $x,y$ are the original two edges containing the attaching region for the band and $z,w$ are the two new edges created by the band attachment, then the four-circuit $x,z,y,w$ bounds a square: the union of the band with thin triangles going back to the vertices. Overall, when we replace $x,y$ with $z,w$, we move from one end of the one-parameter family of colorings on this square to the other, and Proposition~\ref{prop:squareUpsilonSlope} shows that $\Upsilon$ changes by an amount within $\frac12 t$ of $\frac{1}{4} t e(b)$, if the common color on these edges is $t$.
    
    Overall, we've shown that each of the Morse moves $m$ making up $\Sigma$ changes $\Upsilon$ by an amount within some constant $d_m$ of a central value $c_m$; since the initial and final graphs are unknotted, the overall change in $\Upsilon$ is $0$ so $\left| \sum_m c_m \right| \le \sum_m d_m$, and we need only to identify these sums with those of the theorem statement. The only tricky part here is identifying the sum over all moves of type (2) of the $t e(b)$ terms with the writhe $w_t(\beta)$.

    To check this, we'll watch how framings of the graph move through the steps in the cobordism. Define a framing of an embedded graph $\Gamma$ as a thickening to an embedded ribbon graph, and similarly a framing for a foam to be a thickening of the original foam, locally an embedding of $\Gamma \times I$ for $\Gamma$ a graph and $I$ an interval, to an embedding of the product of a ribbon-thickening of $\Gamma$ and $I$. Now, fix some framing of our foam that induces the trivial framing on the initial unknotted $\Theta$-graph (the one coming from the graph drawn on a sphere) and the $0$-framing on each boundary knot; this will induce some framing $f_{\text{end}}$ on the final unknotted $\Theta$-graph.

    Now, for any framing $f$ on an embedded colored $\Theta$-graph $\Gamma$, we may define a real number $|f|$ as follows: for each pair $i<j$ the union of the distinct edges $e_i,e_j \in E(\Gamma)$ is a knot $K_{i,j}$, and $f$ induces a framing on $K_{i,j}$; let this be the $n_{i,j}$-framing in the standard numbering. Then $$|f| = \sum_{i<j} t_i t_j n_{i,j}.$$ We will compute $|f_{\text{end}}|$ in two ways. First, if $\beta \subset S^2 \times I$ is any framed spherical braid, we obtain a framed, unknotted $\Theta$-graph in $S^3$ by capping off $S^2 \times I$ with a pair of $B^3$'s, each containing a vertex incident to each of the $n$ strands of $\beta$ (Figure~\ref{fig:FoamyBraidClosures}, top). Call this framing $f_{\beta}$; by counting crossings in a diagram of $\beta$ we have $|f_{\beta}| = w_t(\beta)$. In the current setting, the framing $f_{\text{end}}$ is in fact $f_{\beta}$ for the braid $\beta$ of the theorem statement, because we may obtain the disjoint union of the initial, trivially framed unknotted $\Theta_n$-graph and the final one framed by $f_{\text{end}}$ by surgery on $\beta$ as in the lower part of Figure~\ref{fig:FoamyBraidClosures}. Therefore, $|f_{\text{end}}| = w_t(\beta)$.

    \begin{figure}
        \centering
        \includegraphics[width=0.4\linewidth]{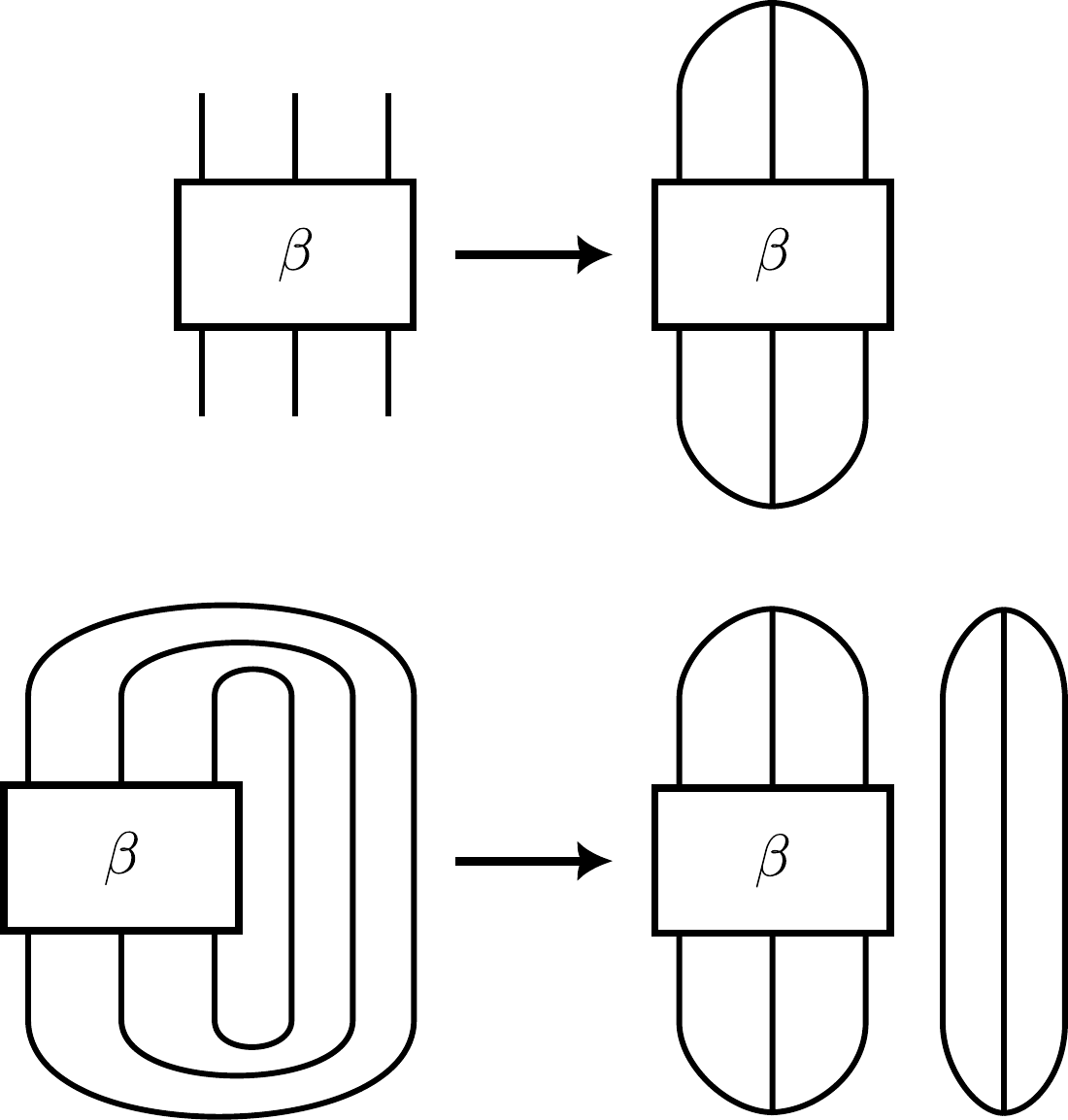}
        \caption{Above: capping off a framed spherical braid $\beta$ to get an $f_\beta$-framed, unknotted $\Theta_n$-graph. Below: surgery from the ordinary closure of $\beta$ to the disjoint union of trivially framed and $f_\beta$-framed unknotted graphs.}
        \label{fig:FoamyBraidClosures}
    \end{figure}

    To compute $|f_{\text{end}}|$ another way, trace the framing through the Morse events giving $\Sigma$. Initially, $|f| = 0$, and most Morse events do not change $|f|$: the initial or final concordance, and moves of type (1) or (3) (the latter because the boundary knots are $0$-framed). For a move of type (2), draw a diagram of $\Gamma$, the graph immediately prior to the move, in which $\Gamma$ is blackboard-framed and the band relating $\Gamma$ to $\Gamma'$, the graph immediately after the move, is contained in a small disk in the plane over which no other arc crosses. Then the framing on $\Gamma'$ that extends over the band is $f'$, the blackboard framing in this drawing. Let the two edges of $\Gamma$ meeting the band be $e$ and $e'$, and divide each into an arc $e_0, e'_0$ before the band and an arc $e_1, e'_1$ after the band. Then the edges of $\Gamma'$ are those of $\Gamma$ not meeting the band, together with $e_0 \cup e_1'$ and $e_0' \cup e_1$. Now, when computing $|f|$ for a blackboard-framed graph, a crossing of sign $\epsilon$ involving a single edge colored $t$ contributes $\epsilon$ to the framing of all cycles involving this edge, so overall contributes $\epsilon t (2-t)$ to $|f|$. By contrast, a crossing of sign $\epsilon$ involving two different edges, colored $t,t'$, contributes to the framing of just one cycle, so overall contributes $-\epsilon t\cdot t'$ to $|f|$ (since, when orienting the cycle to count its writhe, one of the edges is reversed so the crossing appears as the opposite sign). Now, when we pass from $\Gamma$ to $\Gamma'$, the crossings in the diagram stay the same, but those between $e_0$ and $e_1'$ or between $e_1$ and $e_0'$ go from being mixed-edge crossings to same-edge crossings, and those between $e_0$ and $e_0'$ or between $e_1$ and $e_1'$ do the reverse. If the common color on $e$ and $e'$ is $t$, then each such change contributes $\pm 2 t$ to $|f'| - |f|$, with the sign depending on the sign of the crossing and on which portion of the edges it lives in; the sum of all of these terms is precisely $t\cdot e(b)$ since those crossings that change type are precisely those crossings that change sign under the unoriented band.
\end{proof}

\section{Torus ${ \Theta }$-graphs and Batson foams}\label{sec:torusThetaWebs}

\subsection{Relationship with the Farey graph}\label{sec:Farey}

\begin{figure}[!h]
    \centering
    \begin{minipage}{0.7\linewidth}
    \begin{subfigure}[t]{0.3\linewidth}
        \centering
        \includegraphics[width=\linewidth]{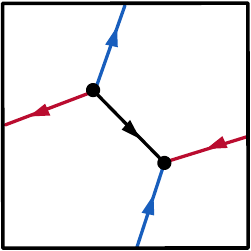}
        \caption{$\Theta_3\begin{psmallmatrix} 1 & 0 \\ 0 & 1 \end{psmallmatrix}$}
        \label{fig:theta_3}
    \end{subfigure}
    \hfill
    \begin{subfigure}[t]{0.3\linewidth}
        \centering
        \includegraphics[width=\linewidth]{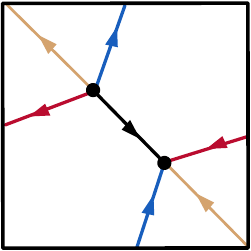}
        \caption{$\Theta_4 \begin{psmallmatrix} 1 & 0 \\ 0 & 1 \end{psmallmatrix}$}
        \label{fig:theta_4}
    \end{subfigure}
    \hfill
    \begin{subfigure}[t]{0.3\linewidth}
        \centering
        \includegraphics[width=\linewidth]{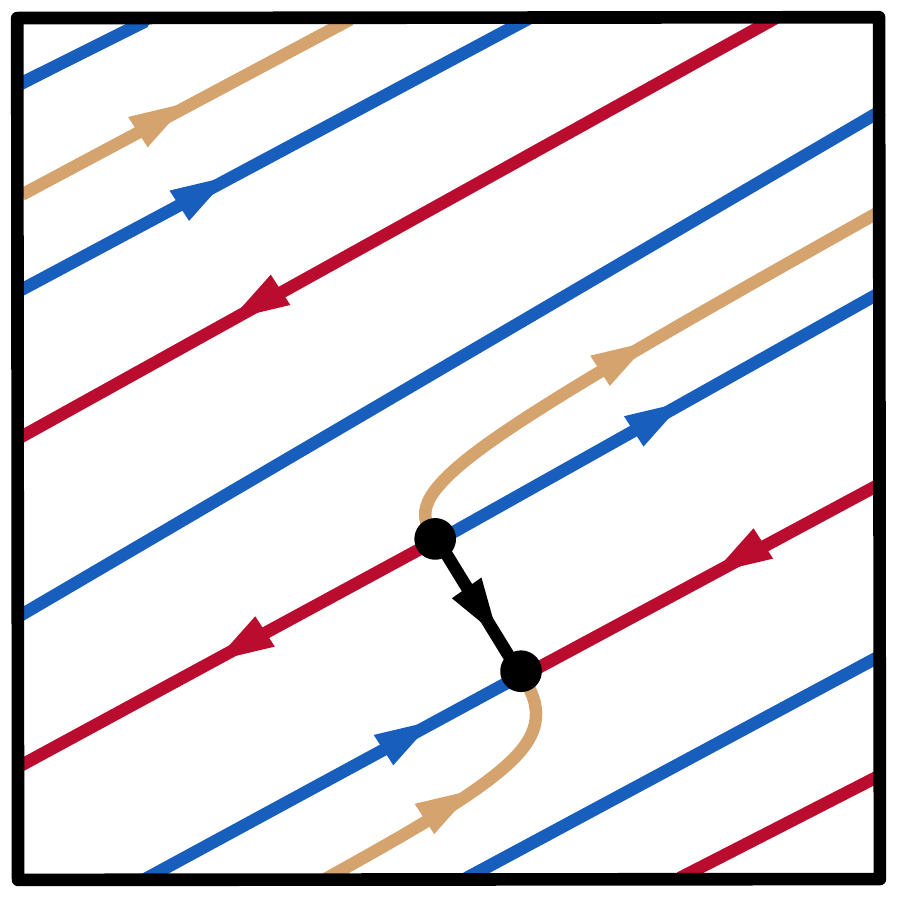}
        \caption{$\Theta_{4}\begin{psmallmatrix} 2 & 1 \\ 3 & 2\end{psmallmatrix}$}
        \label{fig:example_theta_4}
    \end{subfigure}
    \end{minipage}
    \caption{Figures \ref{fig:theta_3} and \ref{fig:theta_4} depict our reference torus $\Theta_3$- \& $\Theta_4$-graphs. All others are obtained by the action of $\mathrm{SL}_2(\mathbb{Z})$ on the torus.}
    \label{fig:thetagraphs}
\end{figure}

\begin{definition}
    The $\Theta_n$-graph is the graph with $2$ vertices and $n$ parallel edges between them. 
\end{definition}
\begin{remark}\label{rem:orientations}
Throughout this note, we orient the edges of $\Theta_n$ so that one of its vertices is a source and the other is a sink.
\end{remark}
A \textit{torus} $\Theta_n$-graph is an embedding of the $\Theta_n$-graph into the Heegaard torus in $S^3$. Some examples are depicted on the square torus in Figure \ref{fig:thetagraphs}. An embedded $\Theta_n$-graph in $S^3$ is \textit{unknotted} if there is a smooth isotopy to the planar $\Theta_n$-graph. 
The following proposition pins down constraints on possible embeddings $\Theta_n \hookrightarrow \mathbb{T}^2$. Intuitively, it says that the complexity of torus $\Theta_n$-graphs plateaus for $n \geq 4$.

\begin{proposition}\label{prop:thetaGraphBigons}
    Suppose that $\Theta_n$ is embedded in the torus such that no face of the complement is a bigon. Then $n \leq 4$, and any two such embeddings of the same $\Theta_n$ are related by a homeomorphism of the torus.  
\end{proposition}
\begin{proof}
    The Euler characteristic of $\Theta_n$ is $2-n$, so the complement must have at least $n-2$ disk regions. Since none of these regions may be a bigon, each has at least four edges on its boundary, and each edge is on the boundary of two regions, so $2n \ge 4(n-2)$ and $n \leq4$. When $n = 1$, any two embeddings of $\Theta_n$ are isotopic, and for $n = 2$ an embedded $\Theta_n$ is an essential, simple closed curve in the torus, any two of which are related by a homeomorphism. For $n = 3$, the complement of the embedding must be a single hexagonal region, which may only happen for one of the two ribbon structures on $\Theta_3$. Similarly, for only one of the six ribbon structures on $\Theta_4$ may the complement be the required two squares.
\end{proof}

\noindent  

The isotopy type of a torus knot is determined by the homology class it represents in the torus, modulo symmetries. Analogously, torus $\Theta_n$-graphs are determined by the classes in $H_1(\mathbb{T}^2;\mathbb{Z})$ represented by their subgraphs.   
For a matrix $\begin{psmallmatrix} p & r \\ q & s \end{psmallmatrix} \in \operatorname{SL}_2(\Z),$ let $\Theta_3\begin{psmallmatrix} p & r \\ q & s \end{psmallmatrix}$ denote the image of the graph in Figure \ref{fig:theta_3} under the automorphism of the torus determined by $\begin{psmallmatrix} p & r \\ q & s \end{psmallmatrix}$. Similarly, $\Theta_4\begin{psmallmatrix} p & r \\ q & s \end{psmallmatrix}$ denotes the image of Figure \ref{fig:theta_4} under this same matrix.

\begin{remark}
Torus $\Theta_3$-graphs are also studied by Calcut--Nieman
\cite{CN25}. If
$\begin{psmallmatrix}p&r\\q&s\end{psmallmatrix}
\in \mathrm{SL}_2(\mathbb Z)$
has strictly positive entries, then, after forgetting edge labels
and orientations, what we call $\Theta_3\begin{psmallmatrix}p&r\\q&s\end{psmallmatrix}$ is isotopic on the
Heegaard torus to their $\theta(p+r,q+s)$.
\end{remark}

\noindent
As objects in the torus, the symmetries of these graphs are such that
$$
\begin{aligned}
\Theta_3\begin{psmallmatrix} p & r \\ q & s \end{psmallmatrix} &\sim_{\,\mathbb{T}^2} \Theta_3\begin{psmallmatrix} -p & -r \\ -q & -s \end{psmallmatrix}\sim_{\,\mathbb{T}^2} \Theta_3\begin{psmallmatrix} r & -p-r \\ s & -q-s \end{psmallmatrix} \\
\Theta_4\begin{psmallmatrix} p & r \\ q & s \end{psmallmatrix} &\sim_{\,\mathbb{T}^2} \Theta_4\begin{psmallmatrix} -p & -r \\ -q & -s \end{psmallmatrix} \sim_{\,\mathbb{T}^2} \Theta_4\begin{psmallmatrix} r & -p \\ s & -q \end{psmallmatrix}.
\end{aligned}
$$
In particular, the $\operatorname{SL}_2(\mathbb{Z})$-action on $\Theta_n$-graphs in the torus descends to a transitive $\operatorname{PSL}_2(\mathbb{Z})$-action. Furthermore, removing an edge from an embedding of $\Theta_4$ gives an embedding of $\Theta_3$ and removing an edge from a $\Theta_3$ gives a simple closed curve; from $\Theta_4\begin{psmallmatrix} p & r \\ q & s \end{psmallmatrix}$ one obtains in this way either of $\Theta_3\begin{psmallmatrix} p & r \\ q & s \end{psmallmatrix}$ or $\Theta_3\begin{psmallmatrix} r & -p \\ s & -q \end{psmallmatrix}$ and from $\Theta_3\begin{psmallmatrix} p & r \\ q & s \end{psmallmatrix}$ one obtains any of the simple closed curves with slope $p/q$, $r/s$, or $(p+r)/(q+s)$. 
These relationships are all summarized by incidences in the infinite Farey tree, $\mathfrak{F}$, defined as the dual of the standard geodesic triangulation of the hyperbolic disk.\footnote{Farey noticed one may label points along the boundary of the hyperbolic disk with elements of $\mathbb{Q} \cup \{\infty\}$ in such a way that two points, corresponding to reduced fractions $\frac{a}{b}$ and $\frac{c}{d}$, are connected by a geodesic arc in the triangulation if and only if $ad - bc = \pm 1$. 
Geodesic triangles can naturally be associated with triples $(\frac{a}{b}, \frac{c}{d}, \frac{a+c}{b+d})$ and each such triangle yields a vertex for $\mathfrak{F}$. A geodesic between $\frac{a}{b}$ and $\frac{c}{d}$ lies in geodesic triangles 
$\triangle(\frac{a}{b},\frac{c}{d}, \frac{a+c}{b+d})$ and $\triangle(\frac{a}{b},\frac{c-a}{d-b},\frac{c}{d}) \sim \triangle(\frac{a}{b},\frac{-c+a}{-d+b},\frac{c}{d})$. This geodesic corresponds to an edge in the Farey graph. The modular group $\mathrm{PSL}_2(\mathbb{Z}) \cong \mathrm{SL}_2(\mathbb{Z}) / \{\pm \mathrm{I}\}$ acts on the boundary points via linear fractional transformations.}

Everything we have said so far pertains to $\Theta_n$-graphs embedded in the torus. There are two additional symmetries when these graphs are regarded as lying on the Heegaard torus in $S^3$. The first is analogous to the fact that $T(1,k)$ represents the unknot for any integer $k$. The second is analogous to the fact that the knots $T(p,q)$ and $T(q,p)$ are isotopic. The proofs of both are obvious.
\begin{proposition}\label{prop:unknottedThetaGraphs}
A torus $\Theta_3$-graph without parallel edges is unknotted in $S^3$ if and only if it admits a representative of the form
$\Theta_3\begin{psmallmatrix} 1 & n \\ 0 & 1 \end{psmallmatrix}$ or $\Theta_3\begin{psmallmatrix} 0 & -1 \\ 1 & n \end{psmallmatrix}$ for some $n \in \mathbb{Z}$.
In particular $\Theta_3$-graphs that are embedded in the Heegaard torus and are unknotted in $S^3$ lie along the subgraph of $\mathfrak{F}$ whose associated geodesic triangles in the hyperbolic disk have a vertex labeled by at least one of $\frac{1}{0}, \frac{0}{1}$.
An identical result holds for torus $\Theta_4$-graphs.
\end{proposition}

\begin{proposition}\label{prop:flippingMeridianAndLongitude}
For $n \in \{3,4\}$, we have $\Theta_n\begin{psmallmatrix} p & r \\ q & s \end{psmallmatrix} \sim_{S^3} \Theta_n\begin{psmallmatrix} s & q \\ r & p \end{psmallmatrix}$. 
\end{proposition}

It is a pleasant consequence of Propositions \ref{prop:unknottedThetaGraphs} and \ref{prop:flippingMeridianAndLongitude} that torus $\Theta_3$/$\Theta_4$-graphs are distinguished among one another by the unoriented isotopy type of their cycles. 
It is also useful to remark that, like torus knots, torus $\Theta$-graphs can be partitioned into ``positive'' and ``negative'' families.
\begin{lemma}\label{lem:positiveRepresentative}
Every torus $\Theta_4$-graph either corresponds to a matrix in $\mathrm{SL}_2(\mathbb{Z})$ with nonnegative coefficients or its mirror does. 
\end{lemma}

\begin{figure}[!h]
    \centering
    \makebox[\linewidth][c]{%
        \includegraphics[width=\linewidth]{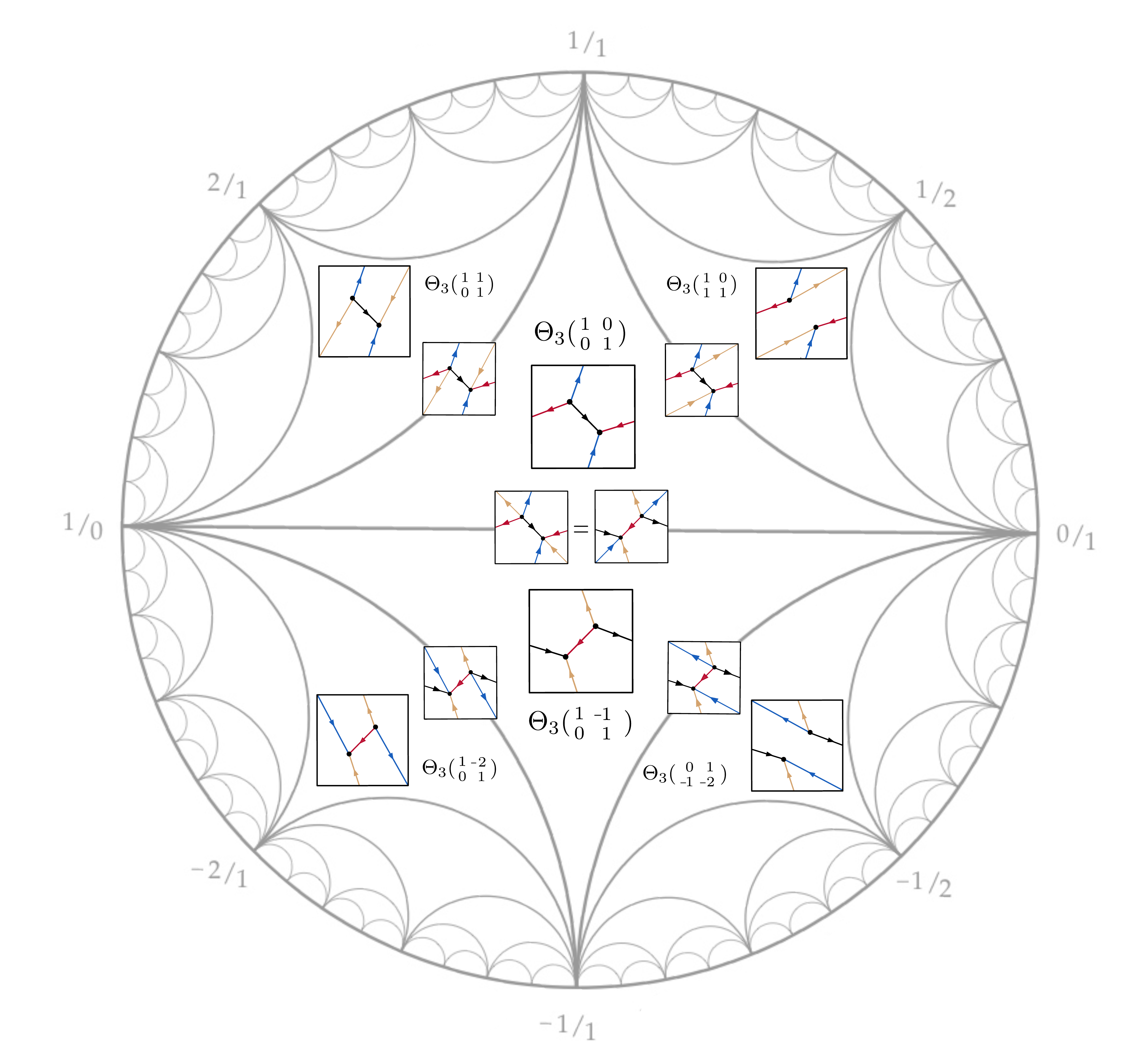}%
    }
    \caption{Each vertex of the Farey tree $\mathfrak{F}$ corresponds to a 2-cell in the depicted triangulation of the hyperbolic disk. There is a unique vertex of $\mathfrak{F}$ determined by the triple $(\frac{p}{q},\frac{r}{s},\frac{p+r}{q+s})$ which can be associated to $\Theta_3\begin{psmallmatrix} p & r \\ q & s \end{psmallmatrix}$. Under this correspondence, $\Theta_4$-graphs correspond to edges of $\mathfrak{F}$ (equivalently, to edges in the hyperbolic triangulation). The above shows this correspondence for the six apparently ``largest'' triangles in the hyperbolic triangulation.}
    \label{fig:FareyTreeAndThetaGraphs}
\end{figure}

\clearpage

\subsection{Foams which generalize Batson's surfaces}

Batson observed \cite{Bat14} that torus knots bound distinguished surfaces, obtained by iteratively attaching nonorientable bands which are framed by the Heegaard torus. Such band attachments are sometimes called \textit{pinch moves} and Batson showed that such moves always describe a nonorientable cobordism from a torus knot to a simpler one. In fact, 
letting $\gamma_4(K)$ denote the smooth, nonorientable 4-genus of $K$, 
Batson showed his surfaces realize ${ \gamma_4(T_{2n,2n-1}) }$ and conjectured they realized ${ \gamma_{4} }$ for all torus knots. Although infinitely many counterexamples to Batson's conjecture have since been discovered \cite{Lobb19, Longo20, BKRST25, Sinha25}, Sabloff showed \cite{Sab23} that for a given torus knot, Batson's surface minimizes $b_1$ across all nonorientable slice surfaces of the same normal Euler number.

We introduce foams which are analogous to Batson's surfaces and obtained via repeated applications of a foam analogue of a pinch move. 

\begin{definition}[Tetrahedron move]\label{defn:tetrahedronMove}
A tetrahedron move on an edge $e$ of a trivalent spatial graph $\Gamma \subset X$ is the local move eliminating $e$, introducing a new edge $e'$, and modifying adjacent edges as displayed in Figure \ref{fig:tetrahedronMove}. Such a move can be interpreted Morse theoretically as a properly embedded 2-complex in $X \times I$ obtained by replacing a neighborhood of $e \times I$ in $\Gamma \times I$ with the local model in Figure \ref{fig:tetrahedronTrace}--- that is, with the cone over the 1-skeleton of a tetrahedron.
\begin{figure}[!h]
    \centering

    \begin{subfigure}[t]{0.48\linewidth}
        \centering
        \includegraphics[width=0.6\linewidth]{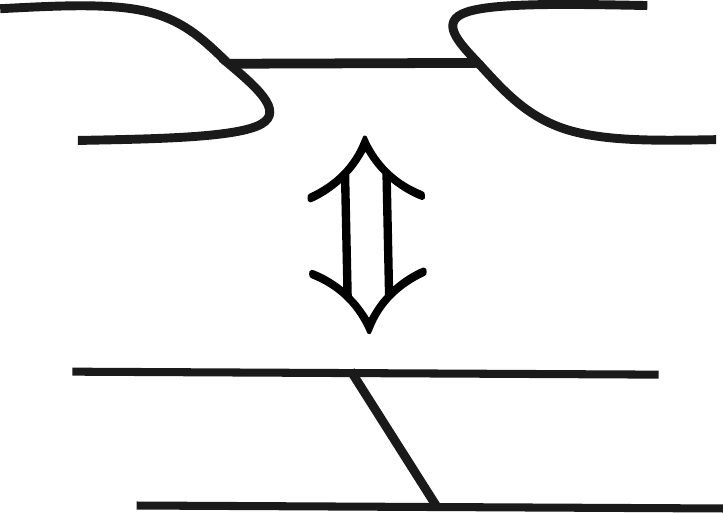}
        \caption{A local picture depicting a tetrahedron move.}
        \label{fig:tetrahedronMove}
    \end{subfigure}
    \hfill
    \begin{subfigure}[t]{0.48\linewidth}
        \centering
        \includegraphics[width=0.7\linewidth]{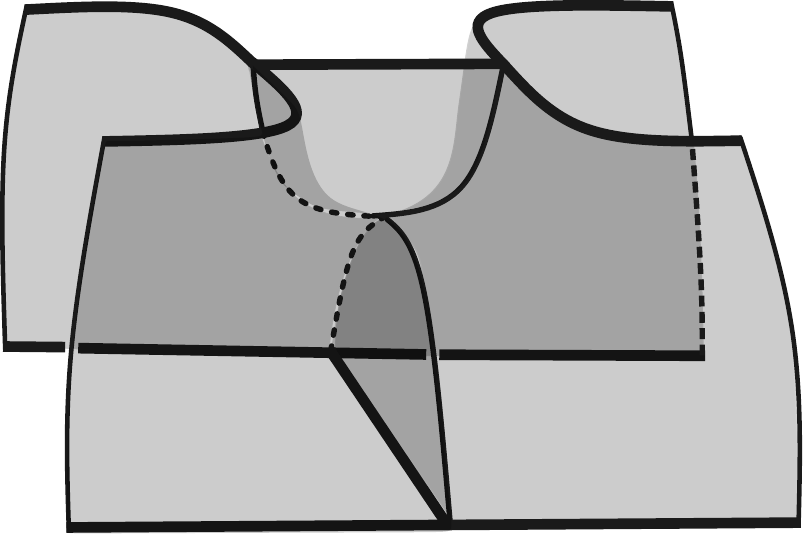}
        \caption{A local picture of a 2-complex obtained as the \\ foam trace of a tetrahedron move.}
        \label{fig:tetrahedronTrace}
    \end{subfigure}

    \medskip 

    \begin{subfigure}[t]{\linewidth}
        \centering
        \includegraphics[width=\linewidth]{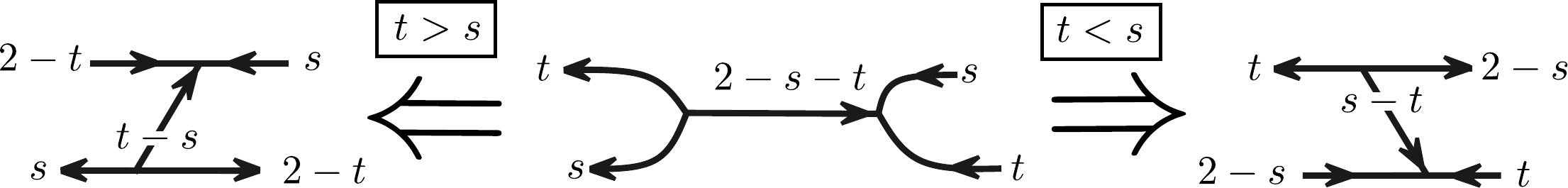}
        \caption{The $\mathbb{R}$-coloring update rule for two level sets of a 2-complex which differ by a tetrahedron move.}
        \label{fig:colorProliferation}
    \end{subfigure}

    \caption{}
    \label{fig:tetrahedronFigures}
\end{figure}

\end{definition}

Our discussion in \S~\ref{sec:Farey} allows us to identify vertices in $\mathfrak{F}$ with isotopy types of $\Theta_3$-graphs embedded in the torus. We say two torus $\Theta_3$-graphs are \textit{Farey-adjacent} if their associated vertices in $\mathfrak{F}$ are connected by an edge. Examining Figure \ref{fig:FareyTreeAndThetaGraphs}, note $\Theta_3 \begin{psmallmatrix} 1 & 0 \\ 0 & 1 \end{psmallmatrix}$ and $\Theta_3 \begin{psmallmatrix} 1 & -1 \\ 0 & 1 \end{psmallmatrix}$ differ by a tetrahedron move. 
Since the edge set of $\mathfrak{F}$ can be obtained as the $\mathrm{PSL}_2(\mathbb{Z})$ orbit of this particular edge, we conclude that all other Farey adjacent torus $\Theta_3$-graphs likewise differ by a tetrahedron move.

  Let $G = \Theta_3\begin{psmallmatrix} p & r \\ q & s \end{psmallmatrix} \in V(\mathfrak{F})$. 
As $\mathfrak{F}$ is a tree, there is a unique shortest path $\gamma$ connecting $G$ to a representative of the unknotted $\Theta_3$-graph and, by the above discussion, each edge traversal corresponds to the application of a tetrahedron move. 
Each such tetrahedron move can be associated to an elementary 2-complex in $T \times I \subset S^3 \times I$ where $T$ denotes the Heegaard torus for $S^3$. 

\begin{definition}[Batson Foam] \label{defn:BatsonFoam}
Let $\gamma$ be the geodesic in $\mathfrak{F}$ from a torus $\Theta_3$-graph $G$ to a representative of the unknotted $\Theta_3$-graph in $S^3$. Given $G' \neq G$ a vertex contained in  $\gamma$, the Batson foam from $G$ to $G'$ is the 2-complex embedded in $T \times I \subset S^3 \times I$ obtained by concatenating the foam traces of the tetrahedron moves associated to each edge in $\gamma' \subset \gamma$ from $G$ to $G'$.
\end{definition}

\begin{figure}
\centering
\includegraphics[width = .9\linewidth]{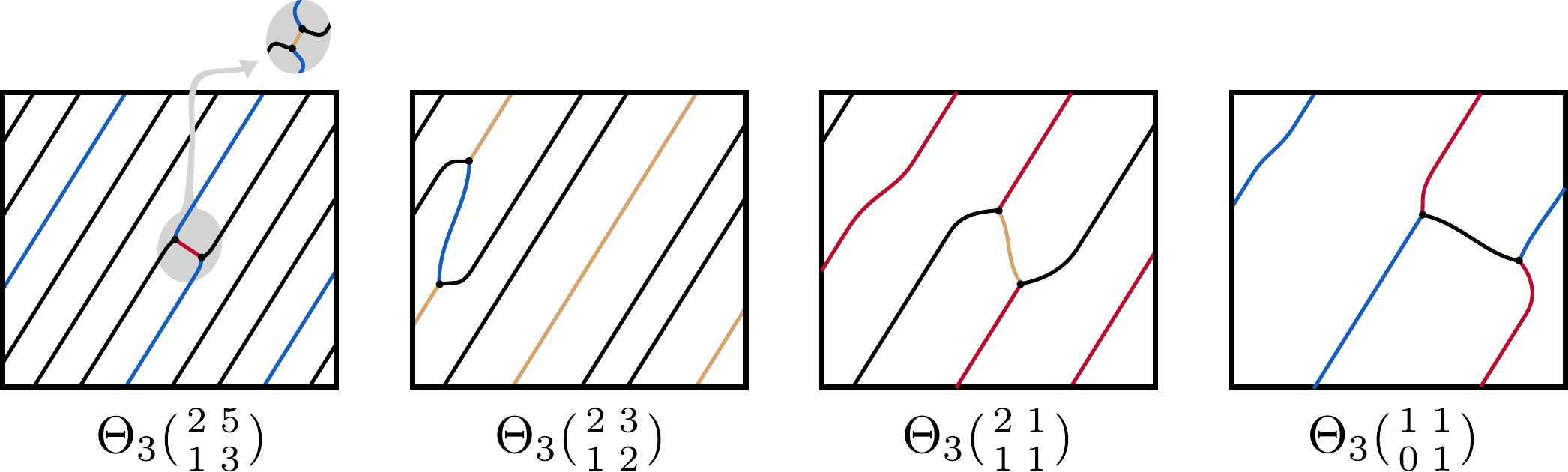}
\caption{A sequence of torus $\Theta_3$-graphs, related by tetrahedron moves (and isotopies). Together, these trace out the Batson foam for $\Theta_3 \begin{psmallmatrix}2 & 5 \\ 1 & 3\end{psmallmatrix}$.}\label{fig:BatsonFoamUncolored}
\end{figure}

\begin{figure}[!h]
\includegraphics[width=\linewidth]{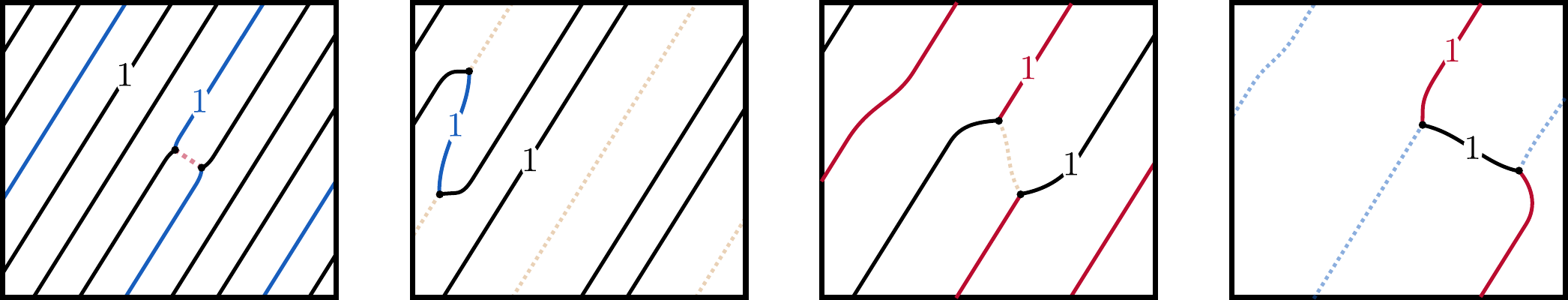}
\caption{A coloring of level sets in the Batson foam for $\Theta_3 \begin{psmallmatrix}2 & 5 \\ 1 & 3\end{psmallmatrix}$. If one ignores the 0-colored (dashed) edges, the level sets are 
$T(4,7) \rightarrow T(2,3) \rightarrow T(0,1)$, exactly the same as in Batson's nonorientable surface for $T(4,7)$.}\label{fig:1colored}
\end{figure}

In particular, each $\Theta_3\begin{psmallmatrix} p & r \\ q & s \end{psmallmatrix}$ cobounds a Batson foam with some representative of the unknotted ${ \Theta_{3} }$-graph; capping off with the cone over the unknotted $\Theta$-graph one obtains a properly embedded 2-complex $\Phi$ in the 4-ball. Introducing colors and orientations to the edges of the regular level sets of $\Phi$ according to the rules set forth in Theorem \ref{thm:FoamBound1}, one may recover and generalize the relationship between $\Upsilon_{T(p,q)}(1)$ and Batson's surfaces. Indeed, Figure \ref{fig:1colored} shows that by coloring two of the three edges of a torus $\Theta_3$-graph with 1, propagating the coloration through tetrahedron moves, and ignoring zero-colored edges, one obtains a foam whose level sets match those of Batson's surfaces.

\begin{figure}[!h]
\includegraphics[width=\linewidth]{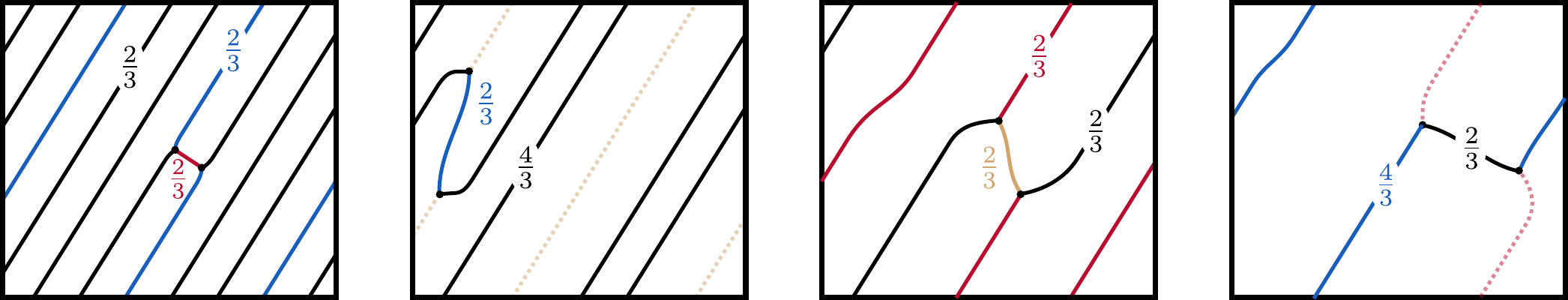}
\caption{The Batson foam bounded by the $(\tfrac{2}{3},\tfrac{2}{3},\tfrac{2}{3})$-colored $\Theta_3\begin{psmallmatrix} 2 & 5 \\ 1 & 3 \end{psmallmatrix}$. Ignoring 0-colored edges, one can think of the level sets of this colored Batson foam as alternating between knots and $\Theta_3$-graphs.}\label{fig:2/3colored}
\end{figure}

\clearpage

\begin{figure}
\includegraphics[width=\linewidth]{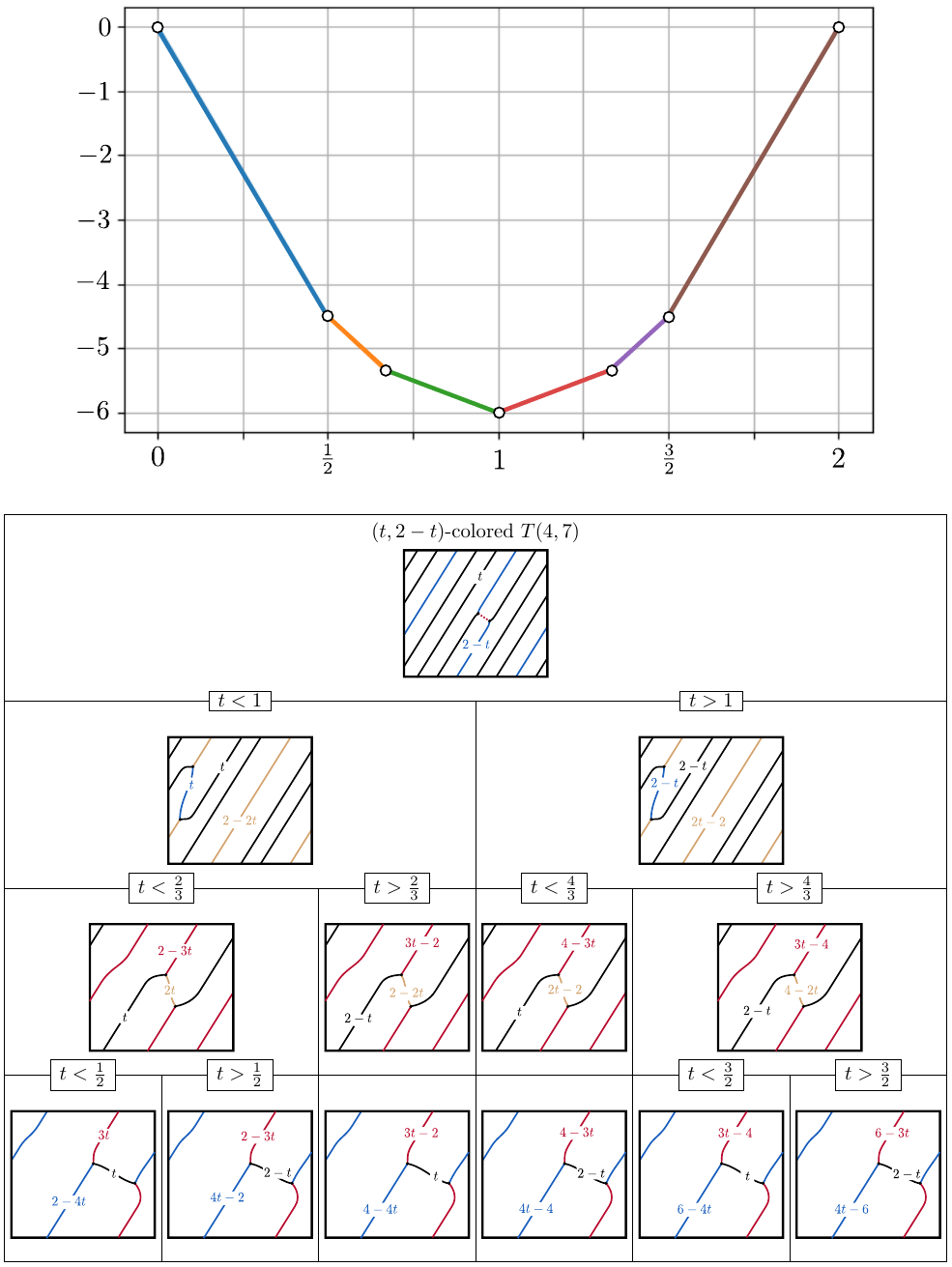}
\caption{
Up top is a plot of $\Upsilon_{T(4,7)}$. The bottom illustrates how an evaluation of $\Upsilon_{T(4,7)}$ at a particular $t\in [0,2]$ can be obtained by propagating the foam inequality through each tetrahedral Morse event in an appropriate $\bm{t}$-coloring of the level sets of the $T(4,7)$-bounded Batson foam. Note that each valid branching in the colored level sets of the foam corresponds to a slope change for $\Upsilon_{T(4,7)}$.
}
\end{figure}

\clearpage

\section{Floer homology of torus $\Theta_n$-graphs}\label{sec:torusThetaGraphHomology} 
By our work in the previous section, it suffices to understand the Floer complexes of torus $\Theta_4$-graphs without parallel edges.
By Lemma \ref{lem:positiveRepresentative} we know that every torus $\Theta_4$-graph is either represented by an element in $\mathrm{SL}_2(\mathbb{Z})$ with nonnegative coefficients or its mirror is. We focus on the former and from the fact that $\mathrm{CFT}(-G) \cong (\mathrm{CFT}(G))^{\ast}$ we are able to obtain models for the tangle Floer complex of any torus $\Theta_4$-graph.

In this section, we construct a model for $\mathrm{CFT}(\Theta_4\begin{psmallmatrix} p & r \\ q & s \end{psmallmatrix})$ when $p,q,r,s \geq 0$. Our construction parallels a well-known association, which we review in \S~\ref{sec:background}, between the standard ``staircase'' model of the knot Floer complex of a positive torus knot $T(p,q)$ and the numerical semigroup 
${\langle p ,q \rangle := \{ap + bq \; | \; a,b \in \mathbb{Z}_{\geq 0}\}}$; another pleasant account may be found in \S~7 of \cite{Al19}. In \S~\ref{sec:submonoids_and_staircase_curves} we describe the combinatorics of certain submonoids of $\mathbb{Z}^2$ and in \S~\ref{sec:HeegaardDiagramsTorusThetas} we show how these submonoids yield Heegaard diagrams for torus $\Theta_4$-graphs and yield a combinatorial description of a staircase model of $\mathrm{CFT}(\Theta_4\begin{psmallmatrix} p & r \\ q & s \end{psmallmatrix})$.

\subsection{Background: staircase knot Floer complexes and numerical semigroups}\label{sec:background}

One may define a graph ${ \Lambda_{p / q}^{+} \subset \mathbb{R}}$ which has vertex set $\langle p, q\rangle$ and an edge for each pair of vertices of unit distance from one another. 
This decomposes ${ \mathbb{Z} }$ into a union of components of ${ \Lambda_{p / q}^{+} }$ and gaps between these components. 
It is not hard to see that ${ \Lambda^{+}_{p/q} }$ has finitely many compact components, each of which has a minimal element with respect to $\mathbb{Z}$'s total order.

\begin{figure}[ht]
\centering

\begin{tikzpicture}[
    x=0.9cm,
    y=1.15cm,
    >=stealth,
]

% -------------------------------------------------
% Semigroup graph for <3,5>
% -------------------------------------------------

% Black edges between consecutive elements of <3,5>
\draw[black, thick] (5,2.5) -- (6,2.5);
\draw[black, thick] (8,2.5) -- (10.5,2.5);

% Draw all vertices white first, with labels underneath
\foreach \x in {-3,...,10} {
    \draw[fill=white, thick] (\x,2.5) circle (3pt);
    \node[below=5pt] at (\x,2.5) {\scriptsize $\x$};
}
\node at (11,2.5) {$\cdots$};
\node at (-4,2.5) {$\cdots$};

% Fill vertices belonging to <3,5> in black
\foreach \x in {0,3,5,6,8,9,10} {
    \draw[fill=black, thick] (\x,2.5) circle (3pt);
}

% -------------------------------------------------
% Half-integer markers corresponding to generators
% -------------------------------------------------

\foreach \x in {-0.5,0.5,2.5,3.5,4.5,6.5,7.5} {
    \draw[thin] (\x,2.35) -- (\x,2.12);
}

% -------------------------------------------------
% Staircase complex
% -------------------------------------------------

% bullet generators: lower row
\node (b0) at (-0.5,0.75) {$\bullet_{0,0}$};
\node (b3) at (2.5,0.75)  {$\bullet_{1,0}$};
\node (b5) at (4.5,0.75)  {$\bullet_{0,1}$};
\node (b8) at (7.5,0.75)  {$\bullet_{1,1}$};

% circle generators: upper row
\node (c0) at (0.5,1.5) {$\circ_{0,0}$};
\node (c3) at (3.5,1.5) {$\circ_{1,0}$};
\node (c5) at (6.5,1.5) {$\circ_{0,1}$};

% differential arrows
\draw[->] (c0) -- node[above left, inner sep = 1pt, font = \scriptsize] {$W$} (b0);
\draw[->] (c0) -- node[above right, inner sep = 1pt, font = \scriptsize] {$Z^{2}$} (b3);

\draw[->] (c3) -- node[above left, inner sep = 1pt, font = \scriptsize] {$W$} (b3);
\draw[->] (c3) -- node[above right, inner sep = 1pt, font = \scriptsize] {$Z$} (b5);

\draw[->] (c5) -- node[above left, inner sep = 1pt, font = \scriptsize] {$W^{2}$} (b5);
\draw[->] (c5) -- node[above right, inner sep = 1pt, font = \scriptsize] {$Z$} (b8);

\end{tikzpicture}

\caption{
The graph $\Lambda_{3 / 5}^+ \subset \mathbb{R}$ for the numerical semigroup $\langle 3,5\rangle$ is seen in black. Underneath, the chain complex $\mathcal{S}_{3,5} \simeq \mathcal{CFK}(T_{3,5})$. The generator $\bullet_{i,j}$ corresponds to the connected component of $\Lambda^+_{3 / 5}$ which contains $3 i + 5j$ and the generator $\circ_{i,j}$ corresponds to the interstice between this component and a gap in the semigroup $\langle 3,5 \rangle$.
}
\label{fig:T35-semigroup-staircase}

\end{figure}
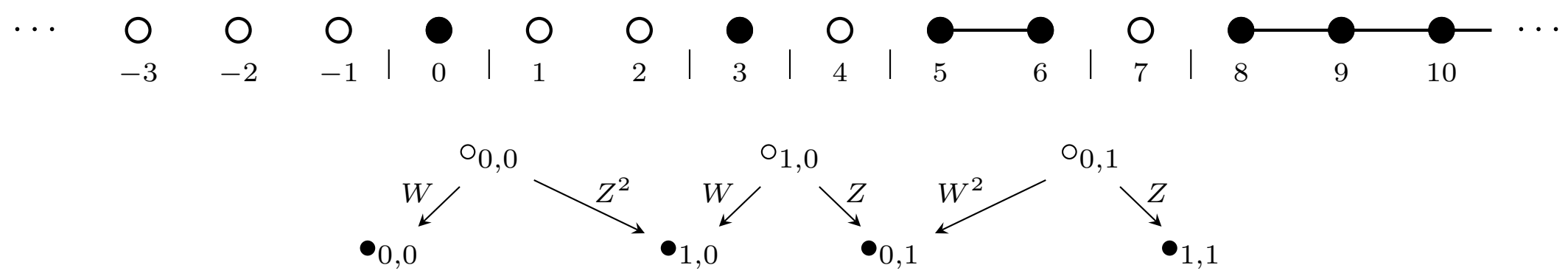

We construct a chain complex ${ \mathcal{S}_{p,q} }$ which is freely generated over ${ \mathbb{F}[W,Z] }$ ---
enumerate components $C_{i,j} \subset \Lambda_{p/q}^{+}$ by the coordinate pairs for their minimal element ${ n = p i + q j }$. We let $\sigma(i,j)$ denote the coordinate pair for the minimal element of the successor component $C_{\sigma(i,j)}$ (the component of $\Lambda_{p / q}^+$ to the immediate right of $C_{i,j}$) and 
let $G_{i,j}$ be the gap in between them.
Associate to ${ C_{i,j} }$ a generator $\bullet_{i,j}$ and, if ${ C_{i,j} }$ is compact, another generator ${ \circ_{i,j} }$. 
There is a bijection between ${ \cup_{i,j} \{\bullet_{i,j}, \circ_{i,j}\} }$ and the interstices between components and gaps of ${ \Lambda^{+}_{p / q} }$; an example is shown in Figure \ref{fig:T35-semigroup-staircase}.
The ``sizes'' of these components and gaps control the coefficients for terms in the differential. 
Construct an endomorphism $\partial$ for $\mathcal{S}_{p,q}$:
$$\partial : \begin{cases} \bullet_{i,j} \mapsto 0 \\ \circ_{i,j}  \mapsto W^{\#(C_{i,j} \cap \mathbb{Z}) }\bullet_{i,j} + Z^{\#(G_{i,j} \cap \mathbb{Z}) } \bullet_{\sigma(i,j)} \end{cases}.$$
It is not hard to see $\partial^2 = 0$.
The \textit{${ H }$-function} of the semigroup ${ \langle p ,q \rangle }$ is
${ H_{p,q}(x) = \# (\langle p ,q \rangle \cap [0,x)) }$.
Letting $c_{p,q} = (p-1)(q-1)$, equip the generators ${ \bullet_{i,j} }$ with a ${ \mathbb{Z}\oplus \mathbb{Z} }$-grading defined by 
$$ (gr_{w}, gr_{z}  )(\bullet_{i,j} ) = \Big( -2H_{p,q}(pi+qj), \;-2H_{p,q}(c_{p,q}-(pi+qj)) \Big)  $$ 
and populate gradings for the $\circ_{i,j}$ according to the rules
$$(gr_w,gr_z)(\partial \bm{x}) = (gr_w,gr_z)(\bm{x}) - (1,1) \qquad (gr_w,gr_z)(W^aZ^b \bm{x}) = (gr_{w},gr_z)(\bm{x})-2(a,b). $$

\begin{proposition}[\cite{BL14, Al19}]\label{prop:standard_torus_knot_gradings}
${ \mathcal{S}_{p,q} }$ is a model for ${ \mathcal{CFK}(T(p,q)) }$ and the gradings ${ (gr_{w},gr_{z}) }$ agree with its absolute Maslov gradings.
\end{proposition}

\subsection{Submonoids of ${ \mathbb{Z}^{2} }$ and staircase curves}\label{sec:submonoids_and_staircase_curves}

Given $\begin{psmallmatrix} p & r \\ q & s \end{psmallmatrix} \in \mathrm{SL}_2(\mathbb{Z})$ with positive coefficients, set 
$$\Sigma_{p/q, r/s} 
:= \{ i (p,r) + j(q,s) \; | \; i,j \in \mathbb{Z}_{\geq 0} \}
.$$ 
In other words, $\Sigma_{p/q,r/s}$ is the submonoid of $\Z^2_{\geq 0}$ generated by the rows $(p,r)$ and $(q,s)$ of the given matrix.

\begin{definition}[Monoid graph]
The \textit{monoid graph} $\Lambda^{+}_{p/q, r/s}$ for $\Sigma_{p /q , r/ s}$ is a planar graph whose vertices are elements of $\Sigma_{p/q, r/s}$ and edges are all length-1 segments connecting elements of $\Sigma_{p/q , r /s}$. Note that all edges are either horizontal or vertical.
\end{definition}

Given points $(x_0,y_0),(x_1,y_1) \in \mathbb{R}^2$, write 
\begin{equation}\label{eqn:relation_on_points_in_R^2}
(x_0,y_0) \leq (x_1,y_1)
\; \Leftrightarrow \;
x_0\leq x_1 \text{ and } y_0\leq y_1.
\end{equation}

\begin{lemma}\label{lem:components_of_monoid_graph}
Let $\begin{psmallmatrix} p & r \\ q & s \end{psmallmatrix} \in \mathrm{SL}_2(\mathbb{Z})$ have all positive entries. 
\begin{enumerate}
\item[(i)]\label{itm:component_count} The connected components of $\Lambda^{+}_{p/q, r/s}$ are in bijection with elements of $([0,s-1
]\times [0,p-1]) \cap \mathbb{Z}^2$. Of these components, $rq= ps-1$ are compact and one is non-compact.
\item[(ii)]\label{itm:uniqueMinimalAndMaximal}
With respect to the partial order $\leq$ on $\mathbb{R}^2$, 
each component of $\Lambda^{+}_{p/q, r/s}$ has a unique minimal element. Moreover, each compact component of $\Lambda^{+}_{p/q, r/s}$ has a unique maximal element. 
\item[(iii)]\label{itm:injectivity_to_collapsed_monoid} The map
    $(x,y) \mapsto x + y$ is injective on the union $K(\Lambda^+_{p/q , r/s})$ of the compact components of $\Lambda^+_{p/q , r/s}$.
\end{enumerate}
\end{lemma}
\begin{proof}
Since $ps-rq = 1$, every pair $(x,y)$ may be uniquely expressed as an integer linear combination of $(p,r)$ and $(q,s)$, so $\Sigma_{p/q, r/s}$ is just the set of $(x,y)$ with nonnegative coefficients in such an expression. Given any $(x,y) \in \Sigma_{p/q, r/s}$, there are $i,j \geq 0$ such that
\begin{equation}\label{eqn:decomposition}
(x,y) = i(p,r) + j (q,s).
\end{equation}
We consider what it should mean for ${ (x\pm 1,y) }$ or ${ (x,y\pm 1) }$ to be in ${ \Sigma_{p/q, r / s} }$:
since ${ ps - rq =1 }$, we have

\begin{equation}\label{eqn:descent}
\begin{aligned}
(x\pm 1,y) &= (i\pm s)(p,r) + (j \mp r)(q,s) \\
(x,y\pm 1) &= (i \mp q)(p,r) + (j \pm p)(q,s)
\end{aligned} \quad 
\text{ and this implies } \quad \begin{aligned} 
(x-1,y) \in \Sigma_{p/q,r/s} &\Leftrightarrow i \geq s.  \\
(x,y-1) \in \Sigma_{p/q,r/s} &\Leftrightarrow j \geq p. \\ 
(x+1,y) \in \Sigma_{p/q,r/s} &\Leftrightarrow j \geq r.  \\
(x,y+1) \in \Sigma_{p/q,r/s} &\Leftrightarrow i \geq q. 
\end{aligned}
\end{equation} 

Observe that if $(x-1,y)$ and $(x,y-1)$ are in $\Sigma_{p/q,r/s}$, then so is $(x-1,y-1)$. 
It follows that each connected component of $\Lambda^{+}_{p/q, r/s}$ has a unique entrywise minimal element.
A similar argument shows that compact components have an entrywise maximal element, proving item (ii).

The above also implies that the vectors occurring as the minimum of a component are precisely those of the form seen in Equation \ref{eqn:decomposition} with $0 \leq i < s$ and $0 \leq j < p$. 
There are $ps$ of these minima and therefore $ps$ components of $\Lambda^+_{p / q , r/s}$. Similarly, the component maxima are enumerated by $(i,j) \in [0,q)\times [0,r)$, of which there are $qr = ps-1$, proving item (i). 

Finally, since the maximum element of a component is one of $i(p,r)+j(q,s)$ with $i<q$ and $j<r$, we should observe that if $(a,b) \in K(\Lambda^{+}_{p/q, r/s})$, then 
$$a+b \leq(q-1)(p+r) + (r-1)(q+s) < (p+r)(q+s).$$ 
For $k \in \{0,1\}$, if the vertices
$(a_k,b_k) = i_k(p,r)+j_k(q,s) \in V(\Lambda^+_{p/q,r/s})$ share a coordinate sum, then the pairs $(i_0,j_0)$ and $(i_1,j_1)$ differ by a nonzero multiple of $(q+s,-p-r)$ and, in particular, their common sum is at least $(p+r)(q+s)$, proving item (iii).
\end{proof}

\begin{definition}\label{defn:staircaseCurves}
Given a coprime pair of positive integers $(a,b)$, its associated \textit{staircase curve} is a piecewise linear curve in $\mathbb{R}^2$ defined by
$$
\Gamma_{a,b}
:=
\bigcup_{n\in\mathbb{Z}}
\left([n,n+1]\times
\left\{\left\lfloor\tfrac{bn}{a}\right\rfloor\right\}\right)
\cup
\bigcup_{n\in\mathbb{Z}}
\left(
\{n\}\times
\left[\,
\left\lfloor\tfrac{bn-b}{a}\right\rfloor,
\left\lfloor\tfrac{bn}{a}\right\rfloor \,
\right]
\right).
$$
\end{definition}

Given $\begin{psmallmatrix} p & r \\ q & s\end{psmallmatrix} \in \mathrm{SL}_2(\mathbb{Z})$ and letting $K(\Lambda^+_{p/q,r/s})$ denote the union of the compact components of $\Lambda^+_{p/q,r/s}$, it is not hard to see that, up to an overall affine shift in the plane, we have $K(\Lambda^+_{p/q,r/s}) \subset \Gamma_{p,r} \cap \Gamma_{q,s}$. 
The assumption $\begin{psmallmatrix} p & r \\ q & s \end{psmallmatrix} \in \mathrm{SL}_2(\mathbb{Z})$ with positive entries implies $\frac{r}{p} < \frac{s}{q}$ and we may translate $\Gamma_{p,r}$ southeast and $\Gamma_{q,s}$ northwest by  $\epsilon \in (0, \tfrac{1}{\sqrt{2}})$. 
Letting $\pi$ denote the $x$-coordinate projection, the vertex set of $\Lambda_{p/q,r/s}^+$ is the intersection of $\mathbb{Z}^2$ with regions of the plane where $\Gamma_{q,s}$ lies above $\Gamma_{p,r}$ with respect to the $y$-coordinate and the edge set is the set of length one segments connecting these points.
We define the \textit{gap graph} $\Lambda^{-}_{p/q,r/s}$ similarly to the description of $\Lambda^+_{p/q,r/s}$ above: its vertex set corresponds to the intersection of $(\mathbb{Z}+\tfrac{1}{2})^2$ with regions of the plane where $\Gamma_{q,s}$ lies below $\Gamma_{p,r}$ and its edges are length-1 segments joining these points.

It will be useful later to enumerate the horizontal and vertical edges of the monoid graph and the gap graph.
We showed in Lemma \ref{lem:components_of_monoid_graph} that each component of ${ \Lambda_{p / q , r / s}^{+} }$ is indexed by a pair ${ (i,j) \in [0,s) \times [0,p) }$; let ${ C_{i,j}^{+} \subset \Lambda_{p / q , r / s}^{+} }$ denote the component whose minimal vertex is ${ (x,y) = i(p,r) + j(q,s) }$. Assuming ${ C_{i,j}^{+}}$ is compact, let ${ C_{i,j}^{-} \subset \Lambda^{-}_{p / q , r/ s} }$ be the component to the immediate right of ${ C_{i,j}^{+} }$.

\begin{lemma}\label{lem:edge_counts}
Let ${ (i,j) }$ be the coordinate pair for  the minimal vertex ${ (x,y)\in C_{i,j}^{+} }$ and let  ${ (i',j') }$ denote the coordinate pair for the minimal vertex ${ (x',y') \in C_{i',j'}^+ }$ of the successor component.  
Letting ${ c_{a,b}:= (a-1)(b-1) }$, we have
$$
\begin{alignedat}{2}
h_{i,j}^{+}
&= H_{p,q}(x'-rq)-H_{p,q}(x-rq)
&\qquad
v_{i,j}^{+}
&= H_{r,s}(y'-rq)-H_{r,s}(y-rq)
\\
h_{i,j}^{-}
&= H_{p,q}(c_{p,q}-x)-H_{p,q}(c_{p,q}-x')
&\qquad
v_{i,j}^{-}
&= H_{r,s}(c_{r,s}-y)-H_{r,s}(c_{r,s}-y')
\end{alignedat}
$$
$$
\begin{aligned}
h_{i,j}^{+} + v_{i,j}^{+} + 1 &= H_{p+r,q+s}(x'+y') - H_{p+r,q+s}(x+y)   \\
h_{i,j}^{-} + v_{i,j}^{-} + 1 &= H_{p+r,q+s}(c_{p+r,q+s} -(x+y)) - H_{p+r,q+s}(c_{p+r,q+s} -(x'+y')).   
\end{aligned}
$$
\begin{proof}
The determinant identity gives $s(p,r) - r(q,s) = (1,0)$
so for any ${ k,\ell \geq 0 }$ we have
$$ (k+s)(p,r) + \ell(q,s) = k(p,r) + (\ell+r)(q,s)+(1,0). $$ 
Hence ${ (k+s)(p,r) + \ell(q,s) }$ and ${ k(p,r) + (\ell+r)(q,s) }$ are joined by a horizontal edge with left endpoint 
$$ X = kp + \ell q + rq \Rightarrow X - rq \in \langle p, q \rangle.  $$
Now, since ${ (x,y) }$ and ${ (x',y') }$ are consecutive component minima, the number of horizontal edges of $\Lambda^+_{p / q , r / s}$ between them coincides with 
$${ h_{i,j}^+ = \#\big( \;\langle p,q \rangle \cap [x-rq,x'-rq)\;\big)= H_{p,q}(x'-rq) - H_{p,q}(x-rq)    . } $$ 
A very similar argument tells us that the bottom vertex of a vertical edge has ${ y }$-coordinate ${ Y }$ satisfying ${ Y - rq \in \langle r ,s \rangle  }$ and the number of vertical edges is therefore equal to the number of such ${ Y }$, which can be enumerated by the formula 
$${ v_{i,j}^{+} =  \#\big(\;\langle r,s \rangle \cap [y-rq, y'-rq)\; \big) = H_{r,s}(y'-rq) - H_{r,s}(y-rq). }$$

We now use the fact that every gap of ${ \langle a,b \rangle }$ is mapped to an element of ${ \langle a,b \rangle  }$ via the reflection 
$$ n \longleftrightarrow c_{a,b}-1 -n  .$$
Hence, to count horizontal edges in ${ C_{i,j}^{-} }$, it suffices to count horizontal and vertical endpoints of the reflected semigroup graph
$$ c_{p,q} - x' \leq X < c_{p,q} - x \qquad c_{r,s}-y' \leq Y < c_{r,s} -y. $$ This gives 
$$
\begin{alignedat}{1}
h_{i,j}^{-}
&= \#\big(\langle p,q\rangle \cap [c_{p,q}-x',\,c_{p,q}-x)\big)
= H_{p,q}(c_{p,q}-x)-H_{p,q}(c_{p,q}-x') \\
v_{i,j}^{-}
&= \#\big(\langle r,s\rangle \cap [c_{r,s}-y',\,c_{r,s}-y)\big)
= H_{r,s}(c_{r,s}-y)-H_{r,s}(c_{r,s}-y') .
\end{alignedat} 
$$
\end{proof}
\end{lemma}

\subsection{Heegaard diagrams and Floer complexes for torus $\Theta_n$-graphs}\label{sec:HeegaardDiagramsTorusThetas}

\begin{figure}[!h]
\includegraphics[width=0.7\linewidth]{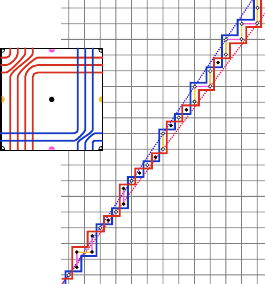}
\caption{In the above, one has translated $\Gamma_{p,r}=\Gamma_{3,4}$ (red) slightly southeast and $\Gamma_{q,s}=\Gamma_{2,3}$ (blue) slightly northwest. Vertices for $\Lambda^+_{p/q,r/s}=\Lambda^+_{3/2,4/3}$ and $\Lambda^-_{3/2,4/3}$ are displayed in white and black respectively and edges for both graphs are seen in orange and pink. Applying the universal covering projection for the torus and identifying four basepoints to be the four components of the image of $\tfrac{1}{2}\mathbb{Z}^2$ under this projection, we obtain the Heegaard diagram $\mathcal{H}_{3/2,4/3}$ (top left) for $\Theta_4\begin{psmallmatrix} 3 & 4 \\ 2 & 3\end{psmallmatrix}$.
}\label{fig:exampleStaircase}
\end{figure}

\begin{definition}
Let ${ p : \mathbb{R}^2 \rightarrow \mathbb{T}^{2} }$ be the universal covering projection and let ${ \Gamma_{p,r}, \Gamma_{q,s} \subset \mathbb{R}^2 }$ be staircase curves associated with ${ \Lambda^{+}_{p/q , r / s} }$. 
Let ${ N_{\alpha} }$ and ${ N_{\beta} }$ be neighborhoods of ${ \Gamma_{p,r} }$ and ${ \Gamma_{q,s} }$ respectively such that
${\mathbb{T}^{2} \setminus  p(N_{\alpha}) }$ and ${ \mathbb{T}^{2} \setminus p(N_{\beta}) }$ are 
closed disks and
${\mathbb{T}^{2} \setminus  p(N_{\alpha} \cup N_{\beta}) }$ is a disjoint union of four closed disks (these are not tubular neighborhoods).
Let ${ \widetilde{\alpha}, \widetilde{\beta} : \mathbb{R}\hookrightarrow \mathbb{R}^2 }$ be smooth embeddings such that 
\begin{itemize}
\item $\mathrm{Im}(\widetilde{\alpha}) \subset N_{\alpha}$ is isotopic to $\Gamma_{p,r}$.
\item $\mathrm{Im}(\widetilde{\beta}) \subset N_{\beta}$ is isotopic to $\Gamma_{q,s}.$
\item ${ \alpha := \mathrm{Im}(p\circ \widetilde{\alpha})  }$ and ${ \beta := \mathrm{Im}(p \circ \widetilde{\beta}) }$ are simple closed curves in the torus. 
\end{itemize}
Note that these embeddings are unique up to isotopy and re-parametrization.
Finally, define the basepoint set to be the image of the vertices and edge-midpoints of $\Lambda^{\pm}_{p/q,r/s}$ under the covering projection;
see Figure \ref{fig:exampleStaircase}. 
We obtain a Heegaard diagram
${ \mathcal{H}_{p / q, r / s} := (\mathbb{T}^{2},\alpha, \beta, \bm{z}) }$.
\end{definition}

\begin{proposition}\label{prop:HeegaardDiagramExistence} 
${ \mathcal{H}_{p / q, r / s}}$ is a Heegaard diagram for ${ \Theta_{4}\begin{psmallmatrix} p & r \\ q & s \end{psmallmatrix} }$. Moreover, the generating set $\mathbb{T}_{\alpha} \cap \mathbb{T}_{\beta}$ contains an element for each component-wise minimum vertex of $\Lambda_{p / q , r / s}^+$ and another element for each component-wise maximum vertex.
\begin{proof}
It is not hard to see that $\mathcal{H}_{p / q, r / s}$ encodes a $\Theta_4$-graph in $S^3$.
Removing all but one basepoint, one can isotope the $\alpha$ and $\beta$ curves until there is a single intersection point. 
The diagram encodes a balanced, sutured manifold with two sutured sphere boundaries. The graph is loopless, has four edges, and hence is a spatial $\Theta_4$-graph. 

As seen in Figure \ref{fig:example_diagram_with_graph}, there is a projection of the encoded $\Theta_4$-graph to the Heegaard torus with no self-intersections; hence the diagram represents a torus $\Theta_4$-graph. By earlier discussion, it is determined by the four distinct torus knots appearing as cycles. Consider the ideal ${ \mathcal{I}_{i,j} = (Z_{i}-1,Z_{j}-1) \leq R_{4} }$. Our conventions are such that
\begin{equation}\label{ref:torusCycleComplexIdentifications}
\begin{alignedat}{3}
\mathrm{CFT}(\mathcal{H}_{p/q,r/s})
    \otimes R_4 / \mathcal{I}_{2,3}
&\simeq\;
&\mathrm{CFT}(\mathcal{H}_{p/q,r/s})
    \otimes R_4 / \mathcal{I}_{1,4}
&\simeq\;
&&\mathcal{CFK}(\,T(p,q)\,)
\\
\mathrm{CFT}(\mathcal{H}_{p/q,r/s})
    \otimes R_4 / \mathcal{I}_{3,4}
&\simeq\;
&\mathrm{CFT}(\mathcal{H}_{p/q,r/s})
    \otimes R_4 / \mathcal{I}_{1,2}
&\simeq\;
&&\mathcal{CFK}(\,T(r,s)\,)
\\
&&
\mathrm{CFT}(\mathcal{H}_{p/q,r/s})
    \otimes R_4 / \mathcal{I}_{1,3}
&\simeq\;
&&\mathcal{CFK}(\,T(|p-r|,|q-s|)\,)
\\
&&
\mathrm{CFT}(\mathcal{H}_{p/q,r/s})
    \otimes R_4 / \mathcal{I}_{2,4}
&\simeq\;
&&\mathcal{CFK}(\,T(p+r,q+s)\,).
\end{alignedat}
\end{equation}

\begin{figure}[!h]
\includegraphics[width=0.7\linewidth]{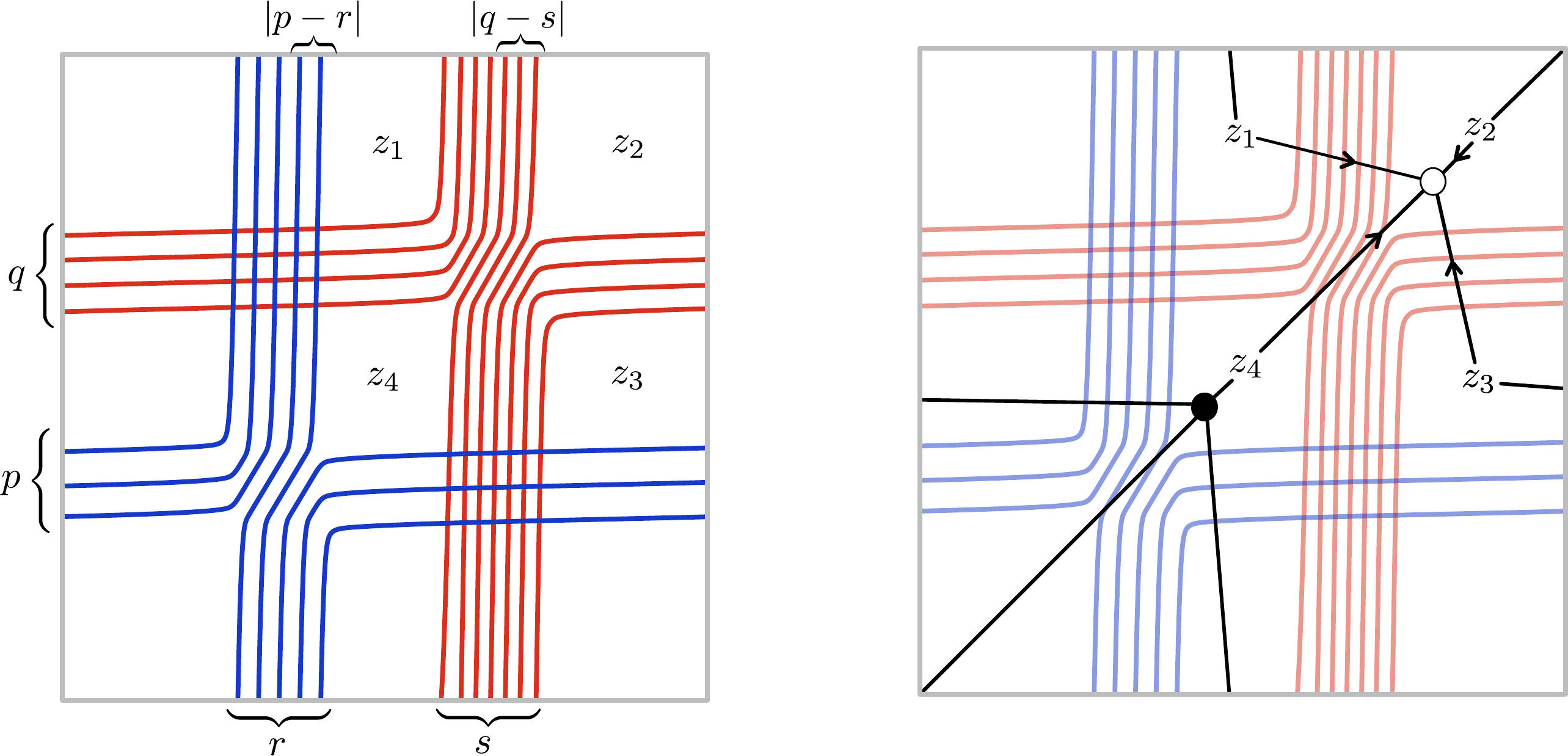}
\caption{The multipointed Heegaard diagram $\mathcal{H}_{p/q, r/s}$ for 
$\Theta_4\begin{psmallmatrix} 
p & r \\ 
q & s
\end{psmallmatrix}
=
\Theta_4\begin{psmallmatrix} 
3 & 5 \\ 
4 & 7
\end{psmallmatrix}
$. 
The torus knot cycles for this $\Theta_4$-graph can be read off from the algebraic intersection of the edges with $\alpha$ and $\beta$.
}\label{fig:example_diagram_with_graph}
\end{figure}

By Lemma \ref{lem:components_of_monoid_graph}(ii), each component $C \subset \Lambda^+_{p /q, r/ s}$ has a minimal element $(x,y) = i(p,r) + j (q,s)$. It is a straightforward consequence of the construction that there is a point in $\Gamma_{p,r} \cap \Gamma_{q,s}$ lying immediately southwest of this minimal element and therefore an associated point in $\widetilde{\alpha} \cap \widetilde{\beta}$. After projecting to $\mathcal{H}_{p/q ,r/s}$, one obtains a generator $\bullet_{i,j} \in \alpha \cap \beta$. 
Lemma \ref{lem:components_of_monoid_graph}(ii) also tells us that a compact component $C$ contains a maximal element $(x',y')=i'(p,r)+j'(q,s)$. It is again straightforward to see there is a point in $\Gamma_{p,r} \cap \Gamma_{q,s}$ lying immediately northeast of this maximal element, yielding a corresponding generator $\circ_{i',j'} \in \alpha \cap \beta$.
\end{proof}
\end{proposition}

Since $\widetilde{\alpha}$ and $\widetilde{\beta}$ approximate the staircase curves $\Gamma_{p,r}$ and $\Gamma_{q,s}$, any bigon contributing to the differential of $\mathrm{CFT}(\mathcal{H}_{p/q, r/s})$ is approximated by a compact region between 
$\Gamma_{p,r}$ and $\Gamma_{q,s}$--- these regions abut components of the graphs $\Lambda^{+}_{p/q,r/s}$ and ${ \Lambda^{-}_{p/q , r / s} }$. 
The following lemma is immediate from the above discussion.

\begin{lemma}\label{lem:graph_component_domains}
For any ${ (i,j) }$ corresponding to a compact component of $\Lambda^+_{p/q,r/s}$, the generator ${ \bullet_{i,j} }$ is a cycle. Moreover
$${ \partial(\circ_{i,j}) =  Z_1^{h_{i,j}^{+}  +v_{i,j}^{+}  +1}Z_2^{h_{i,j}^{+}  } Z_4^{v_{i,j}^{+}  } \bullet_{i,j} + Z_3^{h_{i,j}^{-}  +v_{i,j}^{-}  +1}Z_2^{v_{i,j}^{-}}Z_4^{h_{i,j}^{-}  } \bullet_{\sigma(i,j)}   }$$
where ${ h_{i,j}^{\pm} }$ and ${ v_{i,j}^{\pm} }$ are the counts of horizontal and vertical edges of ${ C_{i,j}^{\pm} }$ which were computed in Lemma \ref{lem:edge_counts}.
\end{lemma}

\begin{proposition}\label{prop:gradings}
Setting ${ (x,y) = i(p,r)+j(q,s) }$ with ${ (i,j) \in [0,s-1] \times [0,p-1] \cap \mathbb{Z}^2 }$, we have
$$
\begin{aligned}
gr_1(\bullet_{i,j})&= - 2H_{p+r,q+s}(x+y) & \qquad &
gr_2(\bullet_{i,j}) &= - 2\Big[ H_{r,s} (c_{r,s} - y )  + H_{p,q} (x -rq) \Big] \\
gr_3(\bullet_{i,j}) &= -2H_{p+r,q+s}(c_{p+r,q+s}-(x+y)) & \qquad &
gr_4(\bullet_{i,j}) &= -2\Big[H_{p,q} (c_{p,q} - x )  + H_{r,s} (y - rq) \Big]
\end{aligned}.
$$

\begin{proof}
Recall ${ (s-1,p-1) }$ is the largest index for a component-minimum vertex of $\Lambda_{p/q,r/s}^+$, yielding the point
$$ (x^{*}, y^{*}  ) = \Big((s-1)p+(p-1)q, \;(s-1)r+(p-1)s\Big) = (c_{p,q} +rq , c_{r,s} + rq  );   $$
this is the minimal vertex of the non-compact component of $\Lambda^{+}_{p / q , r/s}$.

Since each ${ gr_{k} }$ agrees with the absolute Maslov grading for the staircase model of ${ \mathcal{CFK}(T(a,b)) }$ where ${ T(a,b) }$ is a given cycle of ${ \Theta_{4}\begin{psmallmatrix}p & r \\ q & s \end{psmallmatrix} }$, the identifications in
Equation \ref{ref:torusCycleComplexIdentifications} 
together with the grading formulas for staircase models for torus knots (see discussion preceding Proposition \ref{prop:standard_torus_knot_gradings})
imply
\begin{equation}\label{eqn:largest_generator_gradings}
\begin{aligned}
gr_{1}(\bullet_{s-1,p-1} ) &= -c_{p+r,q+s} & \qquad 
gr_{2}(\bullet_{s-1,p-1} ) &= -c_{p,q} \\ 
gr_{3}(\bullet_{s-1,p-1} ) &= 0 &\qquad 
gr_{4}(\bullet_{s-1,p-1} ) &= - c_{r,s}.
\end{aligned}
\end{equation}

Setting $(i',j') = \sigma(i,j)$ and ${ (x',y') = i'(p,r)+j'(q,s) }$, Lemmas \ref{lem:edge_counts} and \ref{lem:graph_component_domains} 
give the relative gradings of homology generators in the staircase complex ${ \mathrm{CFT}(\mathcal{H}_{p / q , r / s}) }$. 
\begin{equation}\label{eqn:relative_gradings_1}
\begin{aligned}
gr_1(\bullet_{i',j'})-gr_1(\bullet_{i,j})&=  -2\Big[H_{p+r,q+s}(x'+y') - H_{p+r,q+s}(x+y) \Big] \\
gr_2(\bullet_{i',j'})-gr_2(\bullet_{i,j}) &= 2\Big[ \left(  H_{r,s} (c_{r,s} - y )-H_{r,s}(c_{r,s} -y' )\right) -\left(H_{p,q} (x' -rq) - H_{p,q} (x -rq) \right) \Big] \\
gr_3(\bullet_{i',j'})- gr_3(\bullet_{i,j}) &= -2\Big[H_{p+r,q+s}(c_{p+r,q+s}-(x'+y')) - H_{p+r,q+s}(c_{p+r,q+s}-(x+y))    \Big]  \\
gr_4(\bullet_{i',j'})-gr_4(\bullet_{i,j}) &= 2\Big[ (H_{p,q} (c_{p,q} - x ) - H_{p,q}(c_{p,q}-x' ))-(H_{r,s} (y' - rq) - H_{r,s} (y -rq) )\Big].
\end{aligned}
\end{equation}

Given a coordinate pair $(i,j)$ describing the minimal vertex of a connected component, its successor $\sigma(i,j)$ is the coordinate pair for the minimal vertex of the connected component immediately to its right. In particular, ${ \sigma^{k}(i,j) }$ is the $k$th successor index pair; observe that the relative grading formulae in Equation \ref{eqn:relative_gradings_1} yield a telescoping sum of differences 
$$
gr(\bullet_{\sigma^k(i,j)}) - gr(\bullet_{i,j}) = \Big(gr(\bullet_{\sigma^k(i,j)}) - gr(\bullet_{\sigma^{k-1}(i,j)})\Big) 
+ 
\cdots  + \Big(gr(\bullet_{\sigma(i,j)}) - gr(\bullet_{i,j})\Big).
$$

In particular, we may use the computation of  ${gr( \bullet_{s-1,p-1}) }$ in Equation \ref{eqn:largest_generator_gradings}, 
the fact that ${ H_{a,b} }$ vanishes on negative values, and the fact that ${ H_{a,b}(c_{a,b}) = \frac{c_{a,b}}{2} }$ to conclude 
as in Proposition \ref{prop:gradings}.
\end{proof}
\end{proposition}

\begin{example}\label{exmp:3_4_2_3_theta}

A Heegaard diagram $\mathcal{H}_{3 / 2 , 4 / 3 }$ for $\Theta_4\begin{psmallmatrix} 3 & 4 \\ 2 & 3 \end{psmallmatrix}$ is seen in Figure \ref{fig:exampleStaircase}. We compute its Floer chain complex to be as below. As in the above discussion, for $(i,j) \in ([0,s-1] \times [0,p-1] )\cap \mathbb{Z}^2$,
$\bullet_{i,j} \in \alpha \cap \beta$ lies just southwest of the $(i,j)$th connected component of $\Lambda_{p / q , r / s}^+$ and $\circ_{i,j}$ lies immediately northeast of the maximal element of this component. 

\[
\begin{tikzcd}
	& {\circ_{0,0}} & {\circ_{0,1}} & {\circ_{1,0}} & {\circ_{0,2}} & {\circ_{1,1}} & {\circ_{2,0}} & {\circ_{1,2}} & {\circ_{2,1}} \\
	\\
	{\bullet_{0,0}} & {\bullet_{0,1}} & {\bullet_{1,0}} & {\bullet_{0,2}} & {\bullet_{1,1}} & {\bullet_{2,0}} & {\bullet_{1,2}} & {\bullet_{2,1}} & {\bullet_{2,2}}
	\arrow["{Z_1}"{description}, from=1-2, to=3-1]
	\arrow["{Z_3^4Z_2^2Z_4}"{description}, from=1-2, to=3-2]
	\arrow["{Z_1}"{description}, from=1-3, to=3-2]
	\arrow["{Z_3}"{description}, from=1-3, to=3-3]
	\arrow["{Z_1}"{description}, from=1-4, to=3-3]
	\arrow["{Z_3^2Z_2}"{description}, from=1-4, to=3-4]
	\arrow["{Z_1}"{description}, from=1-5, to=3-4]
	\arrow["{Z_3}"{description}, from=1-5, to=3-5]
	\arrow["{Z_1}"{description}, from=1-6, to=3-5]
	\arrow["{Z_3}"{description}, from=1-6, to=3-6]
	\arrow["{Z_1^2Z_4}"{description}, from=1-7, to=3-6]
	\arrow["{Z_3}"{description}, from=1-7, to=3-7]
	\arrow["{Z_1}"{description}, from=1-8, to=3-7]
	\arrow["{Z_3}"{description}, from=1-8, to=3-8]
	\arrow["{Z_2Z_1^4Z_4^2}"{description}, from=1-9, to=3-8]
	\arrow["{Z_3}"{description}, from=1-9, to=3-9]
\end{tikzcd}
\]

\end{example}

We conclude this section with the observation that torus $\Theta_4$-graphs without parallel edges always bound squares.

\begin{lemma}\label{lem:alternatingSupport}
There is a square circuit $\square$ of ${ \Theta_{4}\begin{psmallmatrix} p & r \\ q & s\end{psmallmatrix} }$ such that given 
$\bm{x}$ a Heegaard state in $\mathrm{CFT}(\mathcal{H}_{p/q, r / s})$ which generates homology, we have $gr_{\square}(\bm{x}) = -2(ps-1) $. 
\begin{proof}
We choose the circuit $\square$ so that $gr_{\square} = gr_1 - gr_2 + gr_3 - gr_4$.
Recalling that
$$  H_{p,q}(c_{p,q}-s)=H_{p,q}(s)+\frac{c_{p,q}}{2}-s,
$$
the desired result is an immediate consequence of the  grading formulae from 
Proposition \ref{prop:gradings} and the identity
$$H_{p+r,q+s}(x+y) = \frac{1}{2} \left( H_{p,q}(x) + H_{p,q}(x-rq) + H_{r,s}(y) + H_{r,s}(y-rq)\right).$$
\end{proof}
\end{lemma}

This computation makes Theorem~\ref{thm:foamBoundTight} immediate:
\begin{proof}[Proof of Theorem~\ref{thm:foamBoundTight}]
    The bounds from Theorem~\ref{thm:FoamBound1} coming from the Batson foam come from a series of applications of Proposition~\ref{prop:squareUpsilonSlope} to torus $\Theta_4$-graphs. For ${ \Theta_{4}\begin{psmallmatrix} p & r \\ q & s\end{psmallmatrix} }$, the relevant slopes are always in between $ps$ and $ps-1 = qr$, but in fact by Lemma~\ref{lem:alternatingSupport}, only the slope $ps-1$ occurs. Therefore, at each stage of the Batson foam, $\Upsilon$ changes by the smallest amount consistent with Theorem~\ref{thm:FoamBound1}.
\end{proof}

\section{A twisting inequality for the $\Upsilon$ of spatial graphs}\label{sec:twistingInequality}

By a \textit{disk-like foam} for $\Theta_n$, we mean an embedding of the cone of $\Theta_n$.
\begin{remark}
In the following proposition and all subsequent text, we write $\binom{k}{2}$ for $\tfrac{k(k-1)}{2}$ and, in particular, the latter is defined for $k<0$. 
\end{remark}

\begin{proposition}\label{prop:foamTwistingInequality}
Suppose $G: \Theta_n \hookrightarrow S^3$ bounds a $\Theta_n$ disk-like foam in a punctured $\overline{\C P^2}$ with facets $A_1,\dots,A_n$. Then each $[A_i]-[A_j]$ is a cycle rel boundary, and if $a_1,\dots,a_n$ are any integers with $a_i - a_j = [\C P^1] \cdot ([A_i]-[A_j])$, then $\mathrm{CFT}(G)$ contains a cycle representing a non-torsion class in degree $\left(-2\binom{a_1}{2}, \dots, -2\binom{a_n}{2} \right)$.\footnote{Note there are really only $n-1$ independent equations on the $a_i$ so there will always be a one-parameter family of solutions, all given by translating any given solution by constants.}
\end{proposition}
\begin{proof}
    Let $\mathcal{H} = (\Sigma,\bm{\alpha},\bm{\beta},\bm{\gamma},\bm{w}=\{w_1,\dots, w_n\})$ be a multipointed Heegaard triple describing a cobordism $\Sigma_{\mathcal{H}}$ between the unknotted $\Theta_n$-graph $U$ and an arbitrary $\Theta_n$-graph $G$ in $W = \overline{\C P^2}^{\circ,\circ}$. 
    The cobordism map $\mathrm{CFT}(U) \to \mathrm{CFT}(G)$ splits over $\operatorname{Spin}^c$ structures on $W$ and, in particular, each basepoint $w_i$ specifies a correspondence 
    $$\mathfrak{s}_{w_i} : \pi_2(\bm{x},\bm{y},\bm{z}) \rightarrow \mathrm{Spin}^c(W).$$
    We fix a homotopy class $\psi \in \pi_2(\bm{x},\bm{y},\bm{z})$ and set $\mathfrak{s}_i = \mathfrak{s}_{w_i}(\psi)$. 
    Forgetting all basepoints except for $w_i$ and $w_j$, one obtains from $\mathcal{H}$ a doubly pointed Heegaard diagram $\mathcal{H}_{i,j} = (\Sigma, \bm{\alpha},\bm{\beta},\bm{\gamma},w_i,w_j)$. Set $X_{i,j} := X_{\mathcal{H}_{i,j}}$ and $\Sigma_{i,j} := \Sigma_{\mathcal{H}_{i,j}}$ and 
    observe that $X_{i,j}= X_{\mathcal{H}} = W$ and $\Sigma_{i,j}= A_i \cup -A_j$. Lemma 3.9 of \cite{Zem19b} tells us
    $$\mathfrak{s}_i - \mathfrak{s}_j = PD[\Sigma_{i,j}] = PD\Big([A_i] - [A_j]\Big)$$ 
    where $PD([A_i])$ is the Poincar\'{e}-Lefschetz dual of $[A_i]$.

    Forgetting all basepoints except for one, one may obtain from $\mathrm{CFT}(G)$ a copy of $\mathrm{CF}^-(S^3)$. Given 
    ${\mathfrak{s} \in \mathrm{Spin}^c(W)}$, the map
    $F_{W,\mathfrak{s}} : \mathrm{CF}^-(S^3) \to \mathrm{CF}^-(S^3)$ is $U^{-(c_1(\mathfrak{s})^2 + 1)/8}$.     
  Setting $\langle c_1(\mathfrak{s}_i), [\C P^1]\rangle = 2a_i - 1$ for $a_i$ as in the statement, the map $\mathrm{CFT}(U) \to \mathrm{CFT}(G)$ sends a cycle in $\mathrm{CFT}(U)$ to a cycle representing a non-torsion class with $i$th grading $-2 \cdot \frac{(2a_i - 1)^2 - 1}{8} = -2 \cdot \binom{a_i}{2}$.
\end{proof}

\begin{corollary}\label{cor:twistInequality}
    Suppose $G'$ is obtained from $G$ by a positive full twist in a disk that meets the $i$th edge algebraically $a_i$ times. Then 
    $$
    -\min\limits_{k \in \Z} \left( \binom{a_1+k}{2} t_1 + \cdots + \binom{a_n + k}{2}t_n \right) \leq \Upsilon_{G'} - \Upsilon_G, $$ and in fact it suffices to take the minimum over $k \in [-\max\{ a_i\}, -\min \{a_i\}] \cap \mathbb{Z}$.
\end{corollary}
\begin{proof}
Recall that the $\Upsilon$-invariant is additive under connected sum and changes sign under mirroring. We obtain the above by applying
Proposition~\ref{prop:foamTwistingInequality} to the vertex-connect sum $\overline{G} \#_v G'$.
\end{proof}

\begin{remark}
    The $a_i$ are not individually invariant under isotopy of the graph, since when a vertex passes through the twisting disk they all change by $1$. This doesn't change the left-hand side of the above formula.
\end{remark}

\begin{figure}[!h]
\includegraphics[width = 0.7\linewidth]{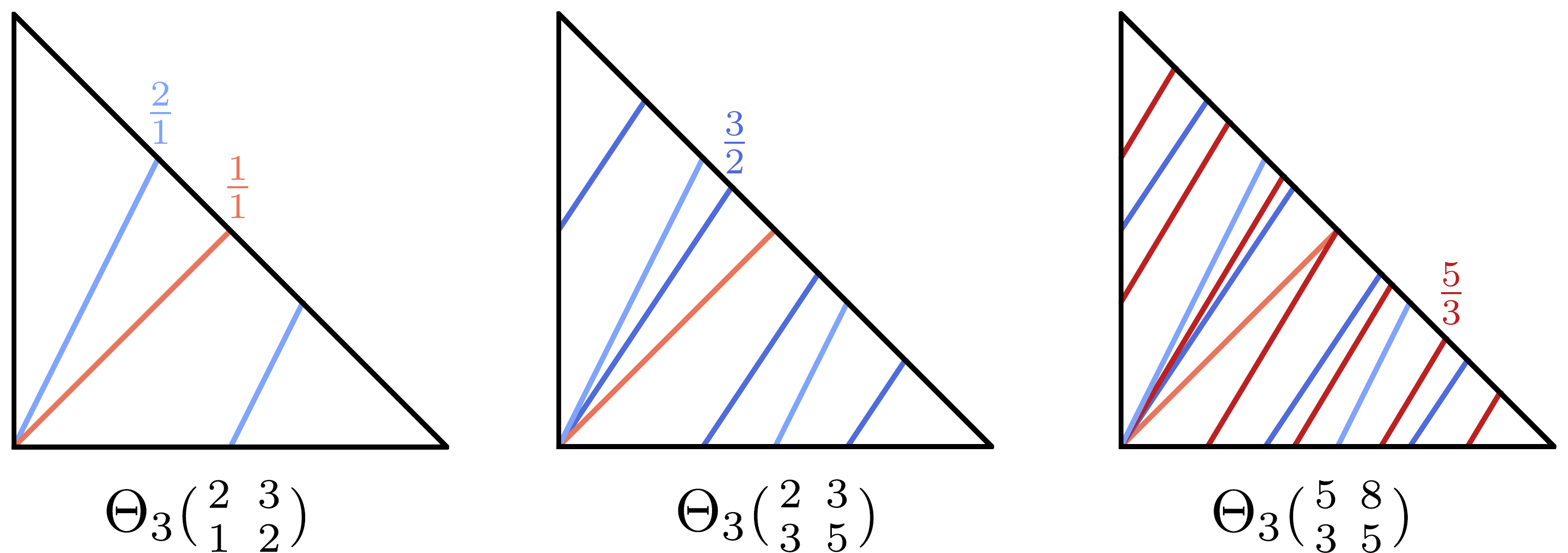}
\caption{
Lines of singularity for the $\Upsilon$-functions of torus $\Theta_3$-graphs differing by full twists. The far-left and center differ by a full twist along the meridian; the center and far-right differ by a twist about the longitude.}\label{fig:upsilon_singular_lines}
\end{figure}

\begin{remark}\label{rem:twistingObtainsAll}
Applying a full twist about a meridian or longitude of the solid torus transforms $\Theta_4\begin{psmallmatrix}p & r \\ q & s\end{psmallmatrix}$ in a predictable manner. In particular, our convention is that a $\pm$-meridional twist corresponds to left multiplication on the representing matrix by $\begin{psmallmatrix}1 & 0 \\ 1 & 1\end{psmallmatrix}^{\pm 1}$ while a $\pm$-longitudinal twist corresponds to left multiplication by  $\begin{psmallmatrix}1 & 1 \\ 0 & 1\end{psmallmatrix}^{\pm 1}$.
Note that it follows as a consequence of the Euclidean algorithm that any $\Theta_4\begin{psmallmatrix}p & r \\ q & s\end{psmallmatrix}$ can be obtained from a torus representative of the unknotted $\Theta_4$-graph by applying a sequence of meridional and longitudinal twists. In fact, there is a symmetry of torus $\Theta_4$-graphs which corresponds to swapping the roles of the meridian and longitude; a longitudinal twist applied to a given torus $\Theta_4$-graph can be interpreted as a meridional twist on an isotopic representative of this same graph.
\end{remark}

\begin{proposition}\label{prop:twistEquality}
Given $\bm{t} = (t_1,t_2,t_3,t_4) \in 2 \Delta^3$ and $p,q,r,s \geq 0$ such that $\Theta_4\begin{psmallmatrix}p & r \\ q+p & s+r\end{psmallmatrix}$ does not represent the unknotted $\Theta_4$-graph. We have
$$\Upsilon_{\Theta_4\begin{psmallmatrix}p & r \\ q+p & s+r\end{psmallmatrix}}(\bm{t}) - \Upsilon_{\Theta_4\begin{psmallmatrix}p & r \\ q & s\end{psmallmatrix}}(\bm{t}) = -\min\limits_{k \in \{-(p+r-1),\dots, 0\}}\left\{\binom{k}{2}t_1 + \binom{r+k}{2}t_2 + \binom{p+r+k}{2}t_3 +  \binom{p+k}{2}t_4  \right\}.$$

\begin{proof} 

Note that $\Gamma_1 := \Theta_4\begin{psmallmatrix}p & r \\ q+p & s+r\end{psmallmatrix}$ may be obtained from $\Gamma_0 := \Theta_4\begin{psmallmatrix}p & r \\ q & s\end{psmallmatrix}$ by applying a twist along a (meridional) disk of the torus which intersects the edges of $\Gamma_0$ algebraically $0,r,p+r,p$ times respectively.

We proceed by induction on the number of tetrahedral Morse events in the Batson foam for our starting torus $\Theta_3$-graph; equivalently, on the number of edges in the Farey tree $\mathfrak{F}$ we must traverse until we arrive at a vertex of $\mathfrak{F}$ corresponding to a representative of the unknotted $\Theta_3$-graph.

For the base case, let $r>0$ and consider the representative 
$\Gamma_0=\Theta_4\begin{psmallmatrix}1 & r \\ 0 & 1\end{psmallmatrix}$ of the unknotted $\Theta_4$-graph. By the remark above, $\Gamma_1 = \Theta_4\begin{psmallmatrix} 1 & r \\ 1 & r+1\end{psmallmatrix}$ differs from $\Gamma_0$ by a full twist and, since $\mathrm{CFT}(\Gamma_1)$ is a staircase complex with homology generators $\bullet_{i,j}$, the fact that $\Upsilon_{\Gamma_0} = 0$ means we have 
$$
\begin{aligned}
\Upsilon_{\Gamma_1}(\bm{t}) &=-\min\limits \Big\{gr_1(\bullet_{i,j})\cdot -\tfrac{t_1}{2} + gr_2(\bullet_{i,j})\cdot -\tfrac{t_2}{2}  + gr_3(\bullet_{i,j})\cdot -\tfrac{t_3}{2}  + gr_4(\bullet_{i,j})\cdot -\tfrac{t_4}{2}   \Big\}
\end{aligned}
$$
where, as in \S~\ref{sec:submonoids_and_staircase_curves}, $\bullet_{i,j}$ corresponds to the semigroup element $(x,y) = i(1,r) + j(1,r+1)$ where $0 \leq i \leq r$ and $j = 0$.
It suffices to show the gradings of the $\bullet_{i,j}$ are such that the above is equal to
\begin{equation}\label{eqn:goalFormula}
-\min\limits_{k \in \{-r,\dots, 0\}}\left\{\binom{k}{2}t_1 +\binom{r+k}{2}t_2 + \binom{1+r+k}{2}t_3 + \binom{1+k}{2}t_4  \right\}.
\end{equation}
Direct substitution into the grading formulae from Proposition \ref{prop:gradings} gives 
$$
\begin{aligned}
gr_1(\bullet_{i,j})&= - 2H_{r+1,r+2}(i(r+1)) &= -2\binom{i+1}{2} \\
gr_2(\bullet_{i,j}) &= - 2\Big[ H_{r,r+1} (c_{r,r+1} - ir )  + H_{1,1} (i -r ) \Big] &= -2\binom{r-i}{2} \\
gr_3(\bullet_{i,j}) &= -2H_{r+1,r+2}(c_{r+1,r+2}-i(r+1)) &= -2\binom{r-i+1}{2} \\
gr_4(\bullet_{i,j}) &= -2\Big[H_{1,1} (- i )  + H_{r,r+1} (ir - r) \Big] &= -2\binom{i}{2}.
\end{aligned}
$$
Setting $k = -i$ and substituting appropriately into the above using $\binom{-i}{2} = \binom{1+i}{2}$ and $\binom{1-i}{2}=\binom{i}{2}$, we obtain exactly the expression from Equation \ref{eqn:goalFormula}.

For the inductive step, observe that if the stated equality holds for a given torus $\Theta_4$-graph, then one obtains the analogous result for its two distinct torus $\Theta_3$-graphs by restricting values of $\bm{t}$ to appropriate facets of $2\Delta^3$. Let
$G_0 = \Theta_3\begin{psmallmatrix}p & r-p \\ q & s-q\end{psmallmatrix}$ be some $\Theta_3$-graph for which the twisting equality holds. There is a $\Theta_4$-graph $\Gamma := \Theta_4\begin{psmallmatrix}p & r \\ q & s\end{psmallmatrix}$ into which $G_0$ embeds and we let $G_1 := \Theta_3\begin{psmallmatrix}p & r \\ q & s\end{psmallmatrix}$ denote the other embedded $\Theta_3$-graph. $\Gamma$ corresponds to an edge in the Farey tree $\mathfrak{F}$ connecting vertices associated to $G_0$ and $G_1$. The action of a meridional full twist on $\Theta_3$- \& $\Theta_4$-graphs corresponds to the automorphism of the Farey tree associated with the matrix $\begin{psmallmatrix} 1 & 0 \\ 1 & 1  \end{psmallmatrix}$. If $\Gamma'$ denotes the result of meridionally twisting $\Gamma$, then $\Gamma'$ corresponds to an edge in $\mathfrak{F}$ with vertices associated to meridional twists $G_i'$ of the $G_i$.

Recall the construction of the square grading from \S~\ref{sec:squareGradings}. Setting $\bm{t}= (t_1,t_2,t_3,t_4) \in 2\Delta^{3} \subset \mathbb{R}_{+}^{4}$, note that as one translates ${ \bm{t} }$ along $\bm{v} = (1,-1,1,-1)$ the change in the $\bm{t}$-modified grading $gr_{\bm{t}}$ is linear with respect to the square grading. 
Explicitly, for $t \in [-\min\{t_1,t_3\},\min\{t_2,t_4\}]$, we have $\bm{t} + t\bm{v} \in 2\Delta^3$ and
\begin{align*}
gr_{\bm{t}+t \bm{v}} &= \tfrac{1}{2}\Big((t_1+t) \cdot gr_1 + (t_2-t)\cdot gr_2 + (t_3+t)\cdot gr_3 + (t_4-t)\cdot gr_4\Big) \\
&= gr_{\bm{t}} + \tfrac{t}{2}\cdot gr_{\square}.
\end{align*}
Set $T := \min\{t_1,t_3\}$ and set $\underline{\bm{t}} := \bm{t}-T \cdot \bm{v}$. There is a facet of the 3-simplex along which the restriction of $\Upsilon_{\Gamma}$ agrees with $\Upsilon_{G_0}$ and our convention is that $\underline{\bm{t}}$ lies in the facet over which $\Upsilon_{\Gamma}$ restricts to
$\Upsilon_{G_0}$. The inductive hypothesis can then be restated as
\begin{equation}\label{eqn:indHyp}
\Upsilon_{\Gamma'}(\underline{\bm{t}}) - \Upsilon_{\Gamma}(\underline{\bm{t}}) = -\min\limits_{k \in \mathbb{Z}}\left\{ \binom{k}{2}(t_1-T) + \binom{r+k}{2}(t_2+T) + \binom{p+r+k}{2}(t_3-T) +  \binom{p+k}{2}(t_4+T)  \right\}.
\end{equation}

The observation $gr_{\bm{t}}=gr_{\underline{\bm{t}}}+\tfrac{T}{2}\cdot gr_{\square}$ together with the computation of $gr_{\square}$ for torus $\Theta_4$-graphs in Lemma \ref{lem:alternatingSupport} tells us
\begin{equation}\label{eqn:upsFormulae_twists}
\begin{split}
\Upsilon_{\Gamma}\left(\bm{t}\right) &= \Upsilon_{\Gamma}(\underline{\bm{t}}) - T(ps-1) \\ 
\Upsilon_{\Gamma'}\left(\bm{t}\right) &= \Upsilon_{\Gamma'}(\underline{\bm{t}}) - T(ps+pr-1).
\end{split}
\end{equation}

Taking their difference and applying the inductive hypothesis, as formulated in Equation \ref{eqn:indHyp}, we have 
$$
\Upsilon_{\Gamma'}(\bm{t}) - \Upsilon_{\Gamma}(\bm{t})  
= \Upsilon_{\Gamma'}(\underline{\bm{t}}) - \Upsilon_{\Gamma}\left(\underline{\bm{t}}\right)  - Tpr 
$$
\begin{equation}\label{eqn:terminal}
=-\min\limits_{k \in \mathbb{Z}}\left\{ \begin{matrix}
 \binom{k}{2}t_1 + \binom{r+k}{2}t_2 +\binom{p+r+k}{2} t_3 + \binom{p+k}{2}t_4 +T\left[-\binom{k}{2} + \binom{r+k}{2} - \binom{p+r+k}{2} + \binom{p+k}{2}   \right] \end{matrix}
\right\} -Tpr.
\end{equation}
One may leverage the identity $\binom{a+b}{2} = \binom{a}{2}+\binom{b}{2}+ab$ to show  
$$-\binom{k}{2} +  \binom{r+k}{2} - \binom{p+r+k}{2} + \binom{p+k}{2}  = -pr,$$
so the $-Tpr$ in Equation \ref{eqn:terminal} cancels; this leaves the desired formula
$$\Upsilon_{\Gamma'}(\bm{t}) - \Upsilon_{\Gamma}(\bm{t}) = -\min\limits_{k \in \mathbb{Z}} \left\{ \binom{k}{2}t_1 + \binom{r+k}{2}t_2 + \binom{p+r+k}{2}t_3 +  \binom{p+k}{2}t_4\right\}.$$
\end{proof}
\end{proposition}

\begin{proof}[Proof of Theorem \ref{thm:twistBoundTight}]
Per Remark \ref{rem:twistingObtainsAll}, every torus $\Theta_4$-graph can be obtained by a sequence of longitudinal and meridional twists starting at a representative of the unknotted $\Theta_4$-graph embedded in the torus. Moreover, the symmetry from Proposition \ref{prop:flippingMeridianAndLongitude} intertwines meridional and longitudinal twists, so a version of Proposition~\ref{prop:twistEquality} also holds for the longitudinal twists $\Theta_4\begin{psmallmatrix}p & r \\ q  & s\end{psmallmatrix} \to \Theta_4\begin{psmallmatrix}p+q & r+s \\ q  & s\end{psmallmatrix}$.

Therefore, induction using Proposition \ref{prop:twistEquality} establishes Theorem \ref{thm:twistBoundTight} for torus $\Theta_3$- and $\Theta_4$-graphs with nonnegative matrix coefficients. One obtains the analogous statement for the nonpositive coefficient case by using Lemma \ref{lem:positiveRepresentative} and the fact that mirroring negates $\Upsilon$.
Given $n\geq 3$, it is a consequence of Proposition \ref{prop:thetaGraphBigons} that any $G: \Theta_n \hookrightarrow S^3$ whose image lies on the Heegaard torus shares an $\Upsilon$ invariant with some torus $\Theta_3$ or $\Theta_4$-graph;  parallel colored edges function, from the perspective of $\Upsilon$, like a single edge whose color is the sum across the parallel edge colors. Moreover each of the parallel edges has the same algebraic intersection number with the twisting disk. Hence, the saturation of the difference inequality for torus $\Theta_n$-graphs is established.

We may obtain a recursive formula for $\Upsilon$ of torus $\Theta_n$-graphs using the above; this generalizes the recursive formula for $\Upsilon$ of torus knots from \cite{FK17}. 
Suppose a sequence of full twists obtains  $G$ from a representative $G_0$ of the unknotted $\Theta_n$-graph $U_{\Theta_n}$:
\begin{equation}\label{eqn:twisting_sequence}
U_{\Theta_n}=G_0 \rightarrow G_1 \rightarrow \cdots \rightarrow G_k = G.
\end{equation}
Since $\Upsilon_{G_0}\equiv 0$, sharpness of the difference formula at each stage of (\ref{eqn:twisting_sequence}) tells us ${\Upsilon_G = \sum_{i=0}^{k-1} (\Upsilon_{G_{i+1}} - \Upsilon_{G_i})}$ and, by the above, this is equal to the sum across the corresponding binomial terms obtained as in the statement of Proposition \ref{prop:twistEquality}.
\end{proof}

\section{Grid homology for spatial graphs}\label{sec:gridHomologyForSpatialGraphs}

The purpose of this section is to prove our formula for the square grading of $\Theta_n$-graphs which bound squares (Corollary \ref{cor: squareBoundSumUp}) using  techniques from grid homology. 
In \S~\ref{sec:gridDiagrams} we recall convenient facts about grid diagrams and in \S~\ref{sec:generalStabilizationLemma} we prove Theorem \ref{thm:vertexStabilizedFloerComplex}; this allows us to compute the square grading for a destabilization using a computation for its stabilization. Finally, in \S~\ref{sec:squareGradingComputation} we obtain our desired formula for the square gradings.

\subsection{Grid complexes}\label{sec:gridDiagrams}

\begin{definition}\label{defn:gridDiagram}
A \textit{grid diagram} is a genus-1 multipointed Heegaard diagram $\mathcal{G} = (\mathbb{T}^2,\bm{\alpha},\bm{\beta},\bm{z})$
such that 
\begin{enumerate}
\item $|\pi_0(\bm{\alpha})| = |\pi_{0}(\bm{\beta})|$. The quantity $n = |\pi_0(\bm{\alpha})|=|\pi_0(\bm{\beta})|$ is called the \textit{grid number} of $\mathcal{G}$.
We adopt the convention in the square torus that the $\bm{\alpha}$ curves are horizontal and the $\bm{\beta}$ curves are vertical.

\item Each ``row'' in $\mathbb{T}^2 \setminus \bm{\alpha}$ contains $\geq 1$ basepoints and likewise for each ``column'' of $\mathbb{T}^2 \setminus \bm{\beta}$. In particular, a given row (or column, or grid cell) may have more than two basepoints and two different rows (or columns) may have different numbers of basepoints.
\end{enumerate}
\end{definition}

\begin{remark}
Each row and column of a grid diagram $\mathcal{G}$ corresponds to a vertex of the encoded bipartite, balanced graph, with the row/column-parity giving the bipartition. 
In the classical setting where a grid diagram encodes a link, each row and column will have exactly two basepoints. In the general setting, having $n$ basepoints in a given row or column corresponds to a valence-$n$ vertex.
\end{remark}

\begin{figure}[htbp]
    \centering
    \begin{subfigure}[t]{0.45\textwidth}
        \centering
        \includegraphics[width=\linewidth]{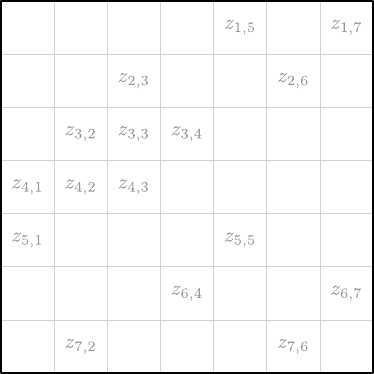}
        \caption{}
        \label{fig:GridDiagram}
    \end{subfigure}
    \hfill
    \begin{subfigure}[t]{0.45\textwidth}
        \centering
        \includegraphics[width=\linewidth]{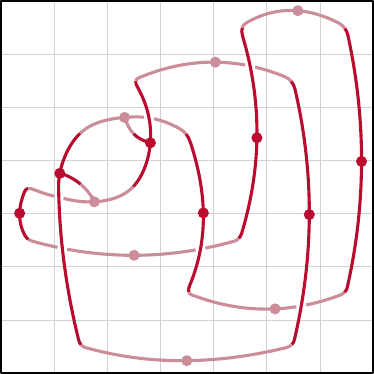}
        \caption{}
        \label{fig:GridDiagramWithGraph}
    \end{subfigure}
    \hfill

    \caption{A Heegaard diagram of genus one for a subdivision of $K_4$ embedded in $S^3$. The right-hand side illustrates how this graph sits in bridge position with respect to the Heegaard torus. This particular graph can be interpreted as arising from a nonorientable band attachment giving a cobordism between the figure-eight knot and a trefoil and, by earlier discussion, contains a cycle which bounds a square. This is reflected in the basepoint configuration $\{z_{3,2},z_{3,3}, z_{4,2}, z_{4,3}\}$ immediately left-of-center in Figure \ref{fig:GridDiagram}.
    }
    \label{fig:ExampleGridDiagram}
\end{figure}

It is well known that every link has a grid diagram which represents it. It is not hard to likewise see that every embedded balanced bipartite graph $\Gamma \subset S^3$ can be presented by a grid diagram after some number of edge stabilizations, starting from any rectilinear planar drawing. 

We concern ourselves with such diagrams because they allow us to obtain constraints on the square grading discussed in \S~\ref{sec:squareGradings}.

\begin{lemma}\label{lem:plaidRegion}
If $\Gamma$ is a bipartite, balanced, spatial graph  and contains a cycle which bounds an embedded square, then there is a Heegaard diagram for some edge stabilization of $\Gamma$ in which the square projects to a subsurface of the Heegaard surface. Moreover, the image under projection of the square contains a single $\alpha\beta$-intersection point with four basepoints arranged around it.
\begin{proof}
The argument is essentially summarized by Figure \ref{fig:GridDiagramWithGraph}. The idea is that one can stabilize the edges of $\Gamma$ to obtain a graph for which a rectilinear diagram yields a grid diagram. It is not hard to see that this rectilinear diagram can be chosen so that the cycle bounding the square is very small.
\end{proof}
\end{lemma}

\begin{definition}
Let $\mathcal{G}$ be a grid diagram for some graph $\Gamma$ in $S^3$. 
The \textit{grid complex} is the tangle Floer complex of the grid diagram  $\mathcal{G}$. Explicitly, it is generated by \textit{grid states} in $\mathbb{T}_{\alpha}\cap \mathbb{T}_{\beta}$ and the differential counts rectangles between grid states and whose interiors are disjoint from members of said states.
\end{definition}

Notice that the existence of a grid diagram for $\Gamma$ forces it to be balanced and bipartite. The complex $\mathrm{CFT}(\mathcal{G})$ possesses a (possibly empty) family of absolute gradings in bijection with the perfect matchings $\mathfrak{M}_{\Gamma}$ of the graph $\Gamma$.
One advantage is that these gradings can be computed combinatorially.

\begin{definition}
Equip $\mathbb{R}^2$ with the partial ordering defined by $(x_0,y_0) < (x_1,y_1)$ if $x_0 < x_1$ and $y_0 < y_1$. Given finite sets $P,Q \subset \mathbb{R}^2$, define 
$$\mathcal{J}(P,Q) = \frac{\#\{(p,q) \in P\times Q \;|\; p < q \} + \#\{(p,q) \in P\times Q \;|\; p > q\}}{2}.$$
\end{definition}

Fix a perfect matching $\mathfrak{m}\in \mathfrak{M}_{\Gamma}$ and let $S_{\mathfrak{m}} \subset \bm{z}$ denote the set of basepoints associated to the matching. 
We think about the intersection points of $\bm{\alpha} \cap \bm{\beta}$ as lying in $\mathbb{Z}^2$ and the basepoints $\bm{z}$ as lying in $(\mathbb{Z}+\frac{1}{2})^2$.

\begin{definition}\label{defn:MaslovFunction}
Given a perfect matching $\mathfrak{m}\in \mathfrak{M}_{\Gamma}$, its \textit{Maslov function} $M_{\mathfrak{m}} : \mathbb{T}_{\alpha} \cap \mathbb{T}_{\beta} \rightarrow \mathbb{Z}$  is
defined by the equation
$$M_{\mathfrak{m}}(\bm{x}) = \mathcal{J}(\bm{x},\bm{x}) - 2\mathcal{J}(\bm{x},S_{\mathfrak{m}})+\mathcal{J}(S_{\mathfrak{m}},S_{\mathfrak{m}})+1.$$
\end{definition}

\begin{proposition}\label{prop:MaslovGradingProof}
 $M_{\mathfrak{m}}$ is the unique function satisfying the following conditions:
\begin{enumerate}
\item\label{itm:Maslov1} Let $\bm{x}^{NW\mathfrak{m}}$ denote the grid state whose components are the upper-left corners of regions in $\mathbb{T}^2 \setminus (\bm{\alpha} \cup \bm{\beta})$ marked with some basepoint in $S_{\mathfrak{m}}$. Then 
$$M_{\mathfrak{m}}(\bm{x}^{NW\mathfrak{m}})=0.$$
\item\label{itm:Maslov2} If $\bm{x}, \bm{y} \in \mathbb{T}_{\alpha}\cap \mathbb{T}_{\beta}$ and are connected by a rectangle $r \in Rect(\bm{x},\bm{y})$, then 
$$M_{\mathfrak{m}}(\bm{x}) - M_{\mathfrak{m}}(\bm{y}) = 1- 2\#(r \cap S_{\mathfrak{m}}) + 2\# (\bm{x} \cap \mathrm{Int}(r)).$$
\end{enumerate}

\begin{proof}
Fix a matching ${ \mathfrak{m} \in \mathfrak{M}_{\Gamma} }$ and note that its associated basepoint set ${ S_{\mathfrak{m}} \subset \bm{z} }$ will contain a single basepoint in each row and column of the grid. 
Hence, the proof of the above proceeds in the same manner as in the proof of  \cite[Lemma 4.3.5]{OSS15}. 

\end{proof}
\end{proposition}

\subsection{A general stabilization lemma for grid complexes}\label{sec:generalStabilizationLemma}

Let $\mathcal{G}$ be a grid diagram for a spatial graph $G$.
Performing a vertex stabilization of $G$ corresponds to selecting a row or column corresponding to the vertex to be stabilized. Having chosen our row or column, we introduce a \textit{generalized vertex stabilization diagram} which contains more basepoints than necessary for any particular vertex stabilization. By choosing an appropriate coloring of the basepoints in this diagram, one may recover the curved complex associated to any stabilization at the chosen vertex, allowing us to treat all vertex stabilizations simultaneously. 

\begin{figure}[htbp]
\centering

\begin{subfigure}{0.31\linewidth}
    \centering
    \includegraphics[width=0.8\linewidth]{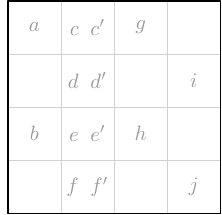}
    \caption{}
    \label{fig:grid1}
\end{subfigure}
\hfill
\begin{subfigure}{0.31\linewidth}
    \centering
    \includegraphics[width=0.8\linewidth]{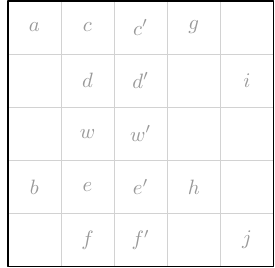}
    \caption{}
    \label{fig:grid2}
\end{subfigure}
\hfill
\begin{subfigure}{0.31\linewidth}
    \centering
    \includegraphics[width=0.8\linewidth]{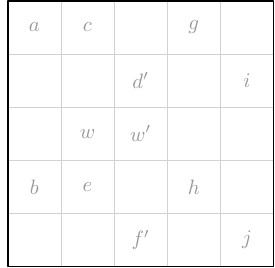}
    \caption{}
    \label{fig:grid3}
\end{subfigure}
\caption{A diagram associated to a graph $G$ with formal basepoint pairs in the second column (\ref{fig:grid1}), the associated generalized vertex stabilization diagram (\ref{fig:grid2}), and a diagram for a particular vertex stabilization of $G$ (\ref{fig:grid3}). 
The diagrams in (\ref{fig:grid2}) and (\ref{fig:grid3})
determine the same curved complex when the former uses a coloring which assigns $\mathfrak{c}: c',d,e',f \mapsto 1$; from the perspective of $\mathrm{CFT}$ this is identical to omitting these four basepoints. The final diagram represents the vertex stabilization in which the new edge associated to $w$ is incident to the edges for $c$ and $e$, while the edge associated to $w'$ is incident to those for $d'$ and $f'$.
}
\label{fig:gridStabilization}
\end{figure}

We describe the process of stabilizing a column vertex; the process for row vertices is identical. See Figure \ref{fig:gridStabilization}. Place a pair of formal basepoints in each grid cell contained in a given column.
When computing the differential of the curved complex associated to a diagram, the coefficient contributed by a holomorphic curve which meets a given grid cell is the product of colorings of the basepoints within said cell and, in particular, a grid cell containing no basepoints may equivalently be thought of as having a basepoint which is colored 1. 
Split the chosen column into two adjacent columns, distributing each pair of basepoints across the two split columns. Insert a new row with $w$ and $w'$ lying in the cells corresponding to the left and right columns and with all other cells labeled 1. 

The result is the generalized vertex stabilization diagram. Using a coloring which assigns some of the formal basepoints to $1$, we obtain a curved complex associated to a particular vertex stabilization of the original graph. In this manner, the curved complexes for any vertex stabilization at the given column vertex are obtained by coloring the basepoints in the generalized vertex stabilization diagram.

We conclude this section with a technical lemma which leads to the proof of Theorem \ref{thm:vertexStabilizedFloerComplex}.

\begin{figure}[!h]
    \centering
    \includegraphics[width = .65\textwidth]{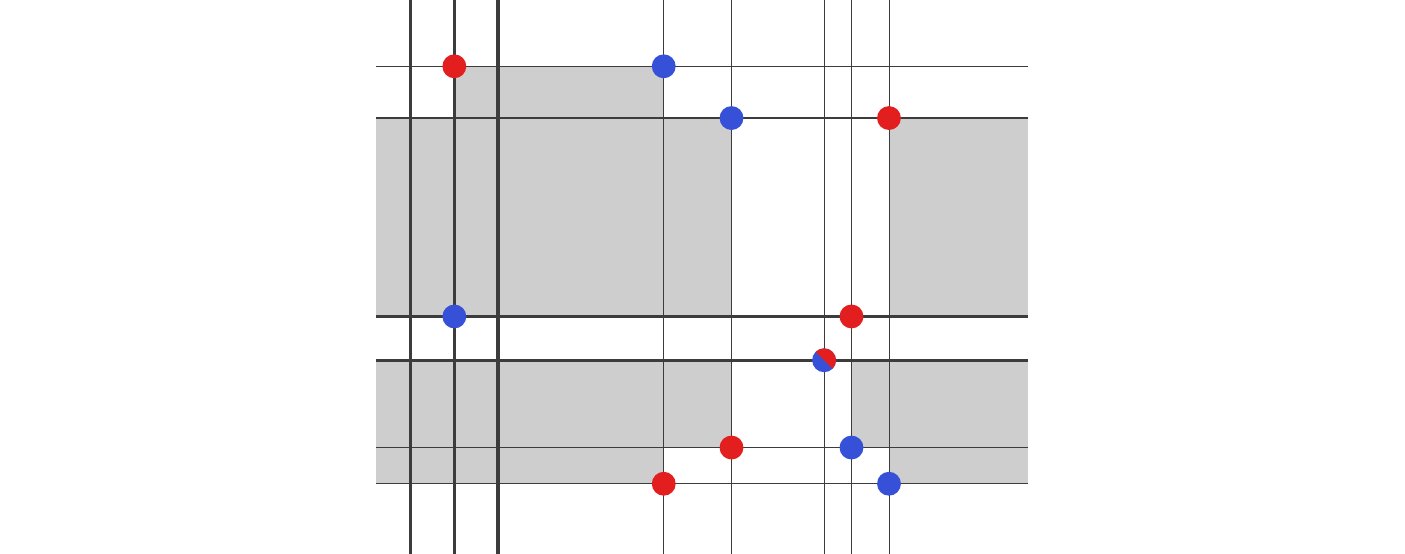}
    \caption{A generic term in the image (drawn in red) of an initial good state (blue) under $(hd_*)^n$. Here, $n=4$. The shaded areas must contain no red or blue marked intersections other than the ones shown, and all other intersections are unchanged between the two states. The feature of this diagram needed in the proof is that all of the red corners apart from the one in the stabilization column lie along a southwest-northeast curve.}\label{fig: GridZigzagGeneric}
\end{figure}

\begin{figure}[!h]
    \centering
    \includegraphics[width = .65\textwidth]{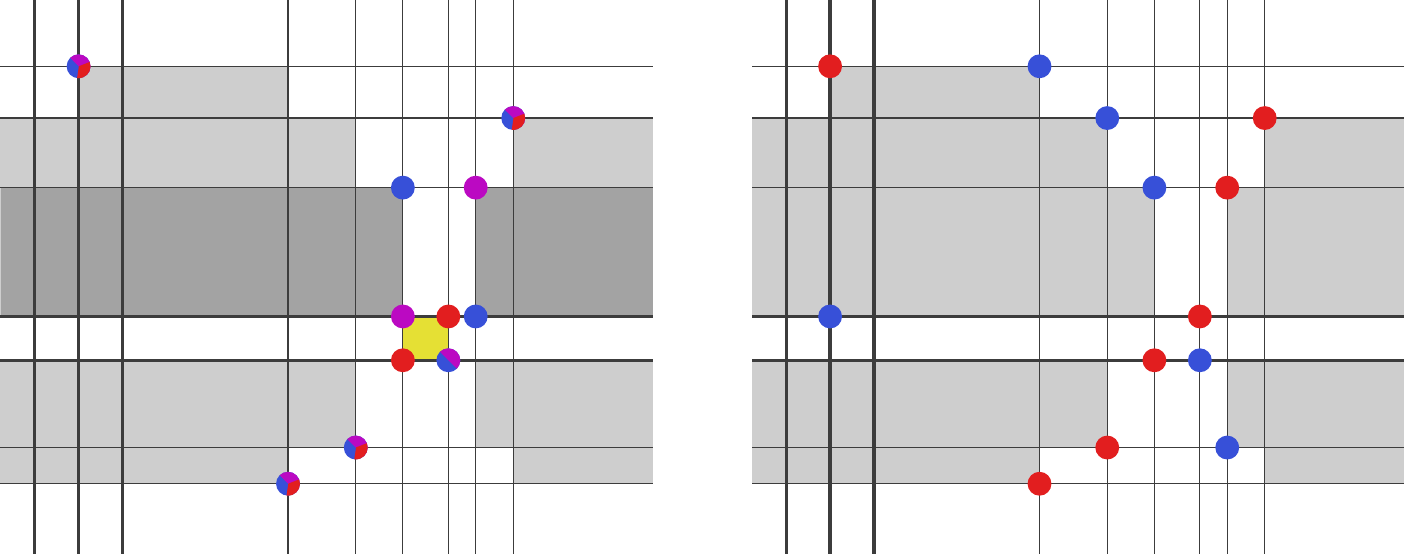}
    \caption{The inductive step in showing that Figure~\ref{fig: GridZigzagGeneric} describes $(hd_*)^n$: On the left, we have the image of a state under $(hd_*)^n$ in blue. In order to acquire a nonzero image under $h$, the relative position of the two marked corners in the stabilization row must swap, so the state drawn in purple is a generic term in the image of the blue state under $d_*$ that will have a nonzero image under $h$, and in red is its image under $h$. Putting this together with Figure~\ref{fig: GridZigzagGeneric}, on the right we have a good state (blue) and a term (red) in its image under $(hd_*)^{n+1}$, showing that the southwest-northeast arrangement of the red points is preserved.}\label{fig: GridZigzagInduction}
\end{figure}

\begin{figure}[!h]
    \centering
    \includegraphics[width = .65\textwidth]{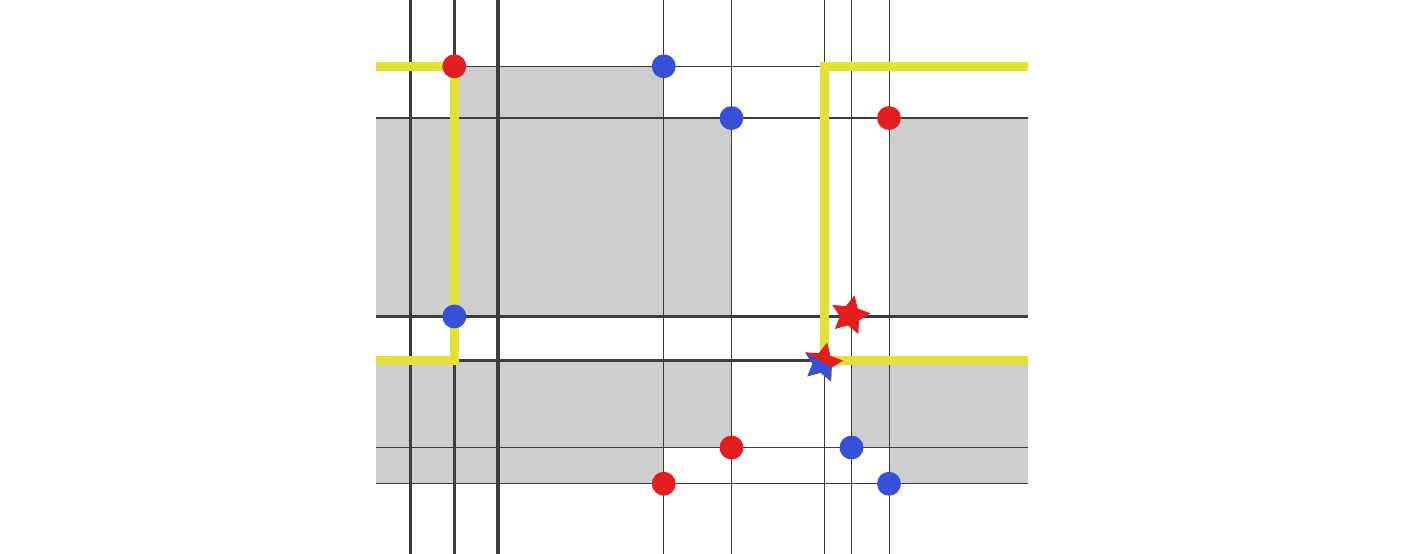}
    \caption{A nonzero term in the image of a state under $\pi d_* (hd_*)^n \iota$ must come from applying $d_*$ to a state as in Figure~\ref{fig: GridZigzagGeneric} and arriving at a good state; this will require moving along a rectangle connecting the red point on the vertical stabilization line to one of the two starred red points. One such rectangle is outlined here, showing that as soon as $n \ge 3$ such a rectangle will not be empty and so not contribute to $d$.}\label{fig: GridZigzagFinal}
\end{figure}

\begin{figure}[p]
    \centering
    \includegraphics[width = .6\textwidth]{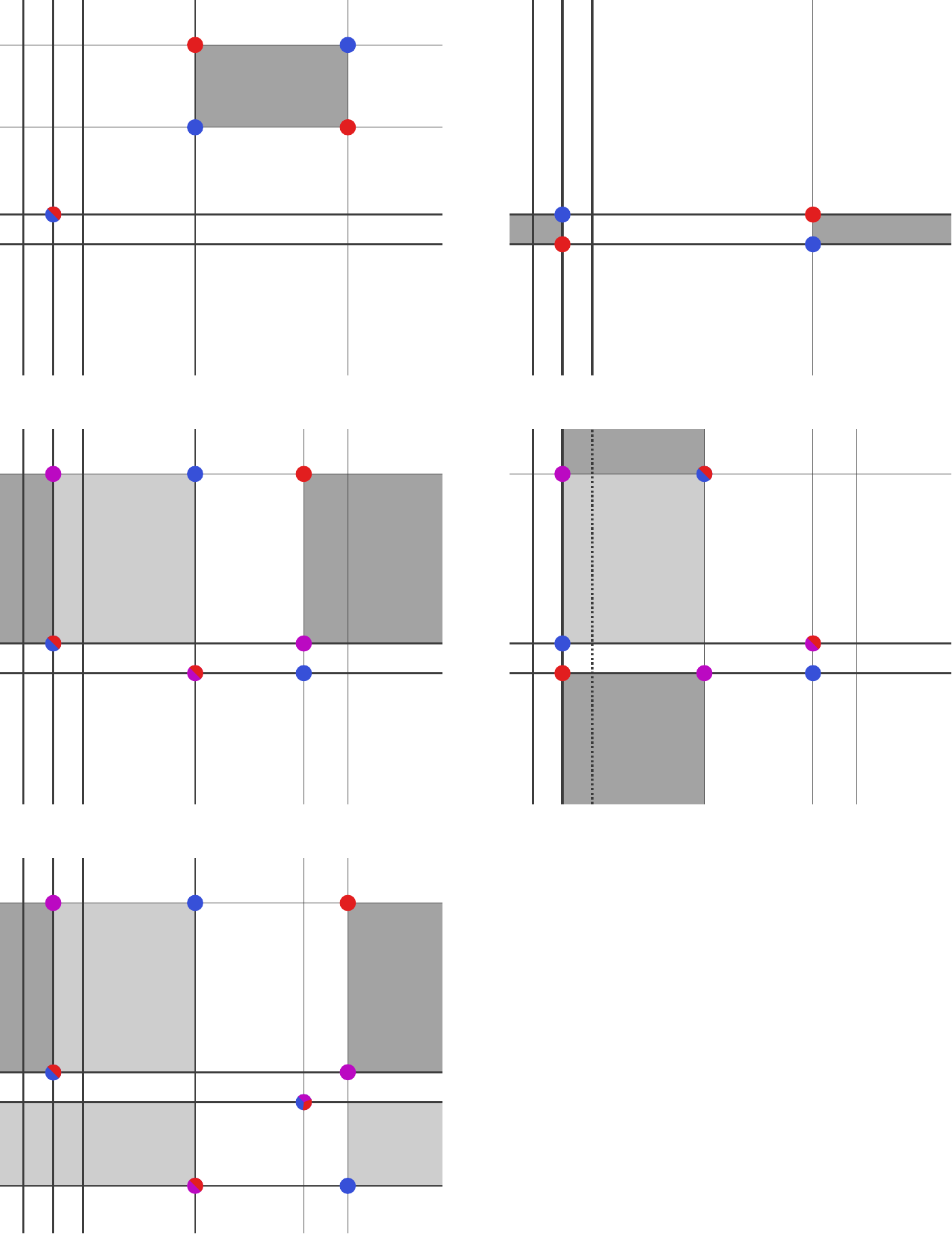}
    \caption{The five types of terms in $\pi d_* (hd_*)^n \iota$. In each image, the initial state is drawn in blue, its image under $(hd_*)^n$ in purple, and the final state in red. The first row has $n = 0$, the second $n=1$, and the last $n=2$. For the second type in the $n=1$ row, the right-hand edge of the rectangles must in fact be immediately to the right of the vertical stabilization edge, since that line contains some marked corner that would otherwise block one of the two rectangles from being empty.}\label{fig:GridZigzagTerms}
\end{figure}

\begin{lemma}\label{lem:generalizedStabilization}
    Suppose that $\tilde{\mathcal{G}}$ is a grid diagram in which all but two adjacent cells in a single row are labeled $1$ and the two excluded cells are labeled $w,w'$. All other cells are labeled by $n^2-n$ distinct variables $z_{i,j}$. 

    Let $\mathcal{G}$ be the destabilized diagram, formed by deleting the row containing the $1$ entries and combining the two columns containing $w$ and $w'$ to form a single column, multiplying their labels; see Figure \ref{fig:gridStabilization}. Then $$\operatorname{CFT}(\tilde{\mathcal{G}}) \simeq \operatorname{CFT}(\mathcal{G}) \otimes K,$$ where $K$ is the two-dimensional curved complex with off-diagonal entries in the differential $$w - \prod_{(i,j) \text{ in the column of }w'} z_{i,j} \qquad \text{ and } \qquad w' - \prod_{(i,j) \text{ in the column of }w}z_{i,j}.$$
\end{lemma}
The general vertex stabilization result will then follow by specializing some of the $z_{i,j}$ to $1$; we work in this ``universal'' context for the present argument since the extra gradings offer rigidity we find useful.
\begin{proof}
    A Heegaard state for $\tilde{\mathcal{G}}$ marks some intersection along the vertical grid line between $w$ and $w'$; call the state \emph{good} if this marked intersection is one of the two corners shared by the cells of $w$ and $w'$, and \emph{bad} if it is anywhere else along the line. In any bad state, the two horizontal grid lines immediately above and below the row containing $w$ and $w'$ will both carry marked intersections away from the vertical grid line between $w,w'$. The proof proceeds in two stages: first, use homological perturbation to find a complex homotopy equivalent to $\operatorname{CFT}(\tilde{\mathcal{G}})$ generated only by the good states, and then identify this complex with $\operatorname{CFT}(\mathcal{G}) \otimes K$.

    Decompose the differential of $\operatorname{CFT}(\tilde{\mathcal{G}})$ as $d = d_0 + d_*$, where $d_0$ counts only rectangles contained entirely within the row of $w$ and $w'$, but not meeting $w$ or $w'$, and $d_*$ counts all other rectangles. Then all nonzero matrix entries of $d_0$ are $1$ and all nonzero entries of $d_*$ carry variables.

    Note that $d_0$ exactly pairs up the bad states according to the relative position of the marked intersections in the row containing $w$ and $w'$: If $x$ is a bad state for which $d_0 x = 0$, then there is a unique bad state $y$ with $d_0y = x$, so we may define an operator $h$ on a state $x$ by $hx = 0$ if $x$ is good or if $d_0x \neq 0$, and otherwise $hx = y$ for the unique $y$ with $d_0y = x$. Then $d_0 h + hd_0$ is exactly the projection operator onto the bad states. 

    Now, provided that the composite $hd_*$ is nilpotent, homological perturbation tells us that $\operatorname{CFT}(\tilde{\mathcal{G}})$ is homotopy equivalent to a complex spanned by just the good states, with differential $$D = \sum_{n \ge 0} \pi d_* (hd_*)^n \iota,$$ where $\iota$ and $\pi$ are the inclusion of good states into all states, and the projection of all states onto good states, respectively. For degree reasons, $hd_*$ is indeed nilpotent: $hd_{*}$ is homogeneous of degree 0 and each of its nonzero matrix entries carries a variable, so each nonzero matrix entry in $(hd_*)^n$ consists of products of at least $n$ variables. By considering a relative grading in which every variable has strictly negative degree (e.g. the sum of all gradings coming from matchings), such an operator on a finite complex must be zero for large enough $n$.

    It remains only to compute the perturbed differential. First, we will show inductively that Figure~\ref{fig: GridZigzagGeneric} depicts a generic term in the image of a good state $x$ under $(hd_*)^n \iota$. To see this, we need only analyze the image of the end state of Figure~\ref{fig: GridZigzagGeneric} under $d_*$, and see when we may arrive at a state not killed by $h$. To do this, a rectangle counted in $d_*$ must exchange the relative position of the two marked points in the row containing $w$ and $w'$, but if this rectangle is tall enough to violate the drawn monotonicity then it will not be empty (this inductive step is shown in Figure~\ref{fig: GridZigzagInduction}).

    Finally, terms in the perturbed differential $\pi d_* (hd_*)^n \iota$ correspond exactly to good states that may be reached from the final state of Figure~\ref{fig: GridZigzagGeneric} by a single application of $d_*$. As shown in Figure~\ref{fig: GridZigzagFinal}, this is possible only when $n \leq 2$, and the resulting terms are enumerated in Figure~\ref{fig:GridZigzagTerms}. Under the correspondence between good states and states of $\operatorname{CFT}(\mathcal{G}) \otimes K$, these become respectively the following types of terms in the differential of $\operatorname{CFT}(\mathcal{G}) \otimes K$:
    \begin{enumerate}
        \item Rectangles not meeting the stabilization point.
        \item The terms $w,w'$ in the differential of $K$.
        \item Rectangles meeting the stabilization point in their boundary.
        \item The remaining terms in the differential of $K$.
        \item Rectangles meeting the stabilization point in their interior.
    \end{enumerate}
\end{proof}

\clearpage

\subsection{A formula for ${ gr_{\square} }$ of ${ \Theta_{n} }$-graphs in ${ S^{3} }$ which bound squares}\label{sec:squareGradingComputation}

The multipointed Heegaard diagrams obtained from Lemma \ref{lem:plaidRegion} contain a local picture satisfying the hypotheses of Lemma \ref{lem:squareGradingConstrained}. These diagrams give models for the tangle Floer complex in which the Heegaard states $\bm{x} \in \mathbb{T}_{\alpha} \cap \mathbb{T}_{\beta}$ have square gradings whose relative support lies in the set $\{0,2\}$. 
For $\Theta_n$-graphs which bound squares, these relative gradings have unique lifts to absolute gradings and, in this section, we determine their support on Heegaard states.

Given a $\Theta_n$-graph in $S^3$ which bounds a square, one may apply two vertex stabilizations to obtain a new graph we call $G'$ which also bounds a square. 
Let ${ v_{1}, v_{2}, v_{3}, v_{4} }$ be the vertices in $G'$ contained in the square. Using the notational convention ${ (\text{source}, \text{target}) }$ for oriented edges, the boundary of the square consists of the oriented edges
$$\square  = \Big(\; e_{1} = (v_{1} , v_{2} ),\; e_{2} = (v_{3},v_{2} ), \; e_{3} = (v_{3},v_{4}  ),\; e_{4} = (v_{1}, v_{4} )\; \Big).$$
We may assume that the two vertex stabilizations used to obtain $G'$ are such that $v_1$ and $v_2$ each have valence $n-1$ while $v_3$ and $v_4$ have valence $3$. In other words, all of the edges of the original $\Theta_n$-graph which do not contribute to the square correspond, in $G'$, to edges connecting $v_1$ to $v_2$.
Once this is done, we apply some number of edge stabilizations to $G'$ to obtain a new graph $G$ which admits a grid diagram. 
By a \textit{long edge} of $G$, we mean a union of edges in $E(G) \setminus \square$ which form a loopless path connecting either ${ v_1 }$ to ${ v_3 }$ or ${ v_2 }$ to ${ v_4 }$; our construction of $G$ guarantees there is exactly one of each type.

We have shown that there's a concrete relationship between the tangle Floer complex of a graph and that of its edge and vertex stabilizations. This relationship is such that if we can determine the square grading on $\mathrm{CFT}(G) \simeq \mathrm{CFT}(\mathcal{G})$, then this grading is determined for the tangle Floer complex of the original $\Theta_n$-graph.

\begin{figure}[htbp]
    \centering

\begin{minipage}[t]{0.23\linewidth}
    \centering
    $\boxed{\mathfrak{m}_1}$\par\vspace{0.4em}
    \includegraphics[width=\linewidth]{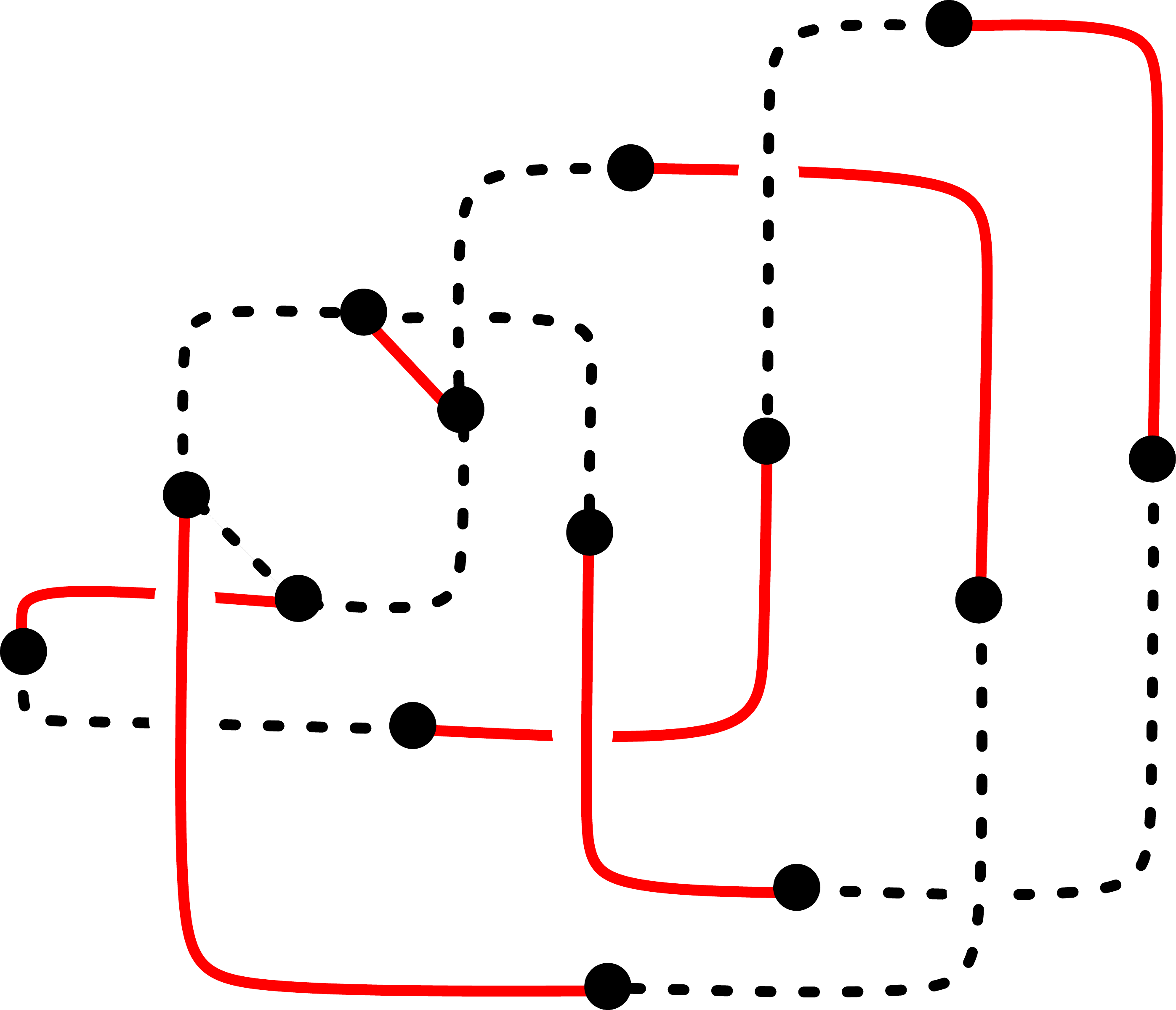}
\end{minipage}
\hfill
\begin{minipage}[t]{0.23\linewidth}
    \centering
    $\boxed{\mathfrak{m}_2}$\par\vspace{0.4em}
    \includegraphics[width=\linewidth]{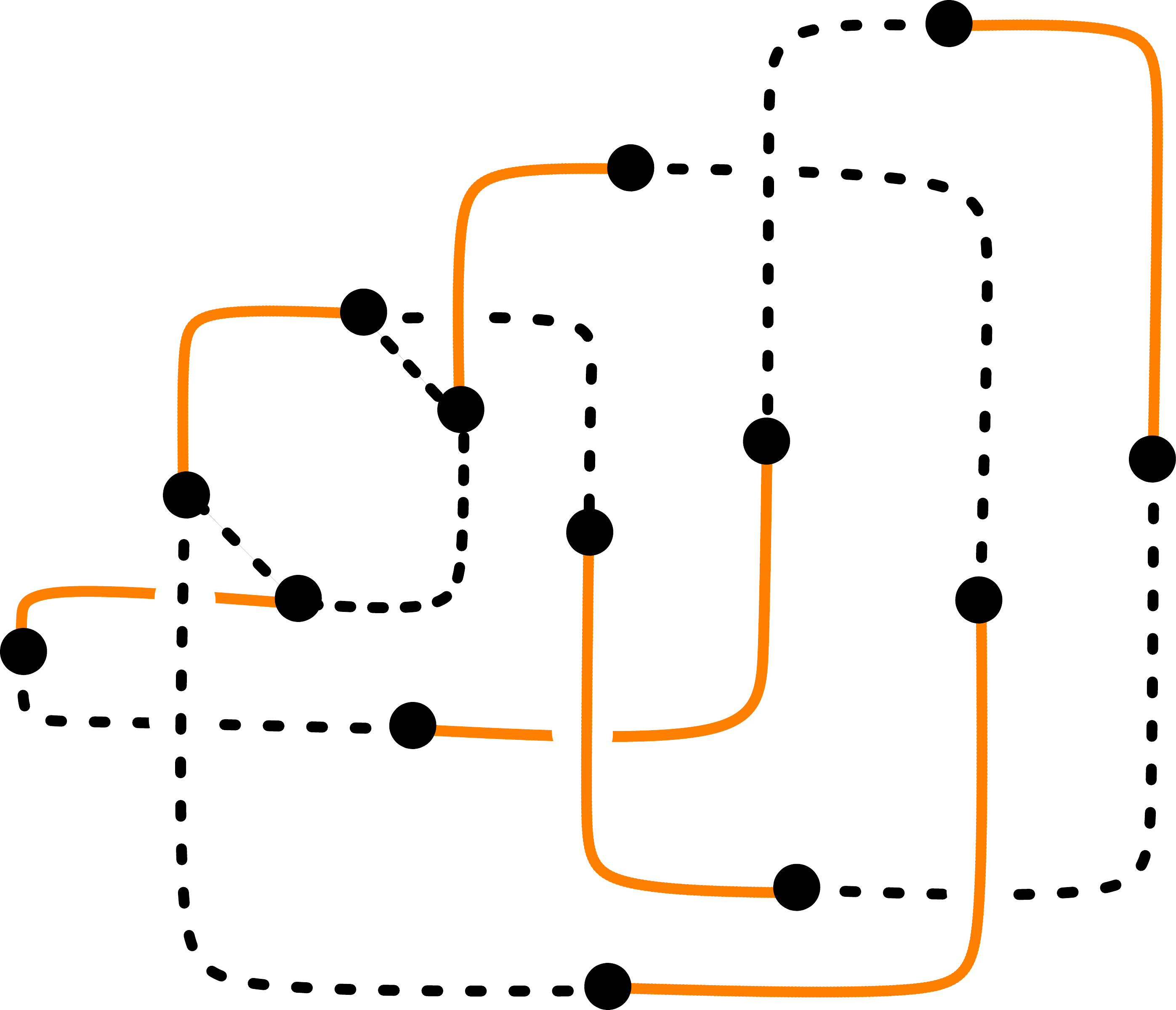}
\end{minipage}
\hfill
\begin{minipage}[t]{0.23\linewidth}
    \centering
    $\boxed{\mathfrak{m}_3}$\par\vspace{0.4em}
    \includegraphics[width=\linewidth]{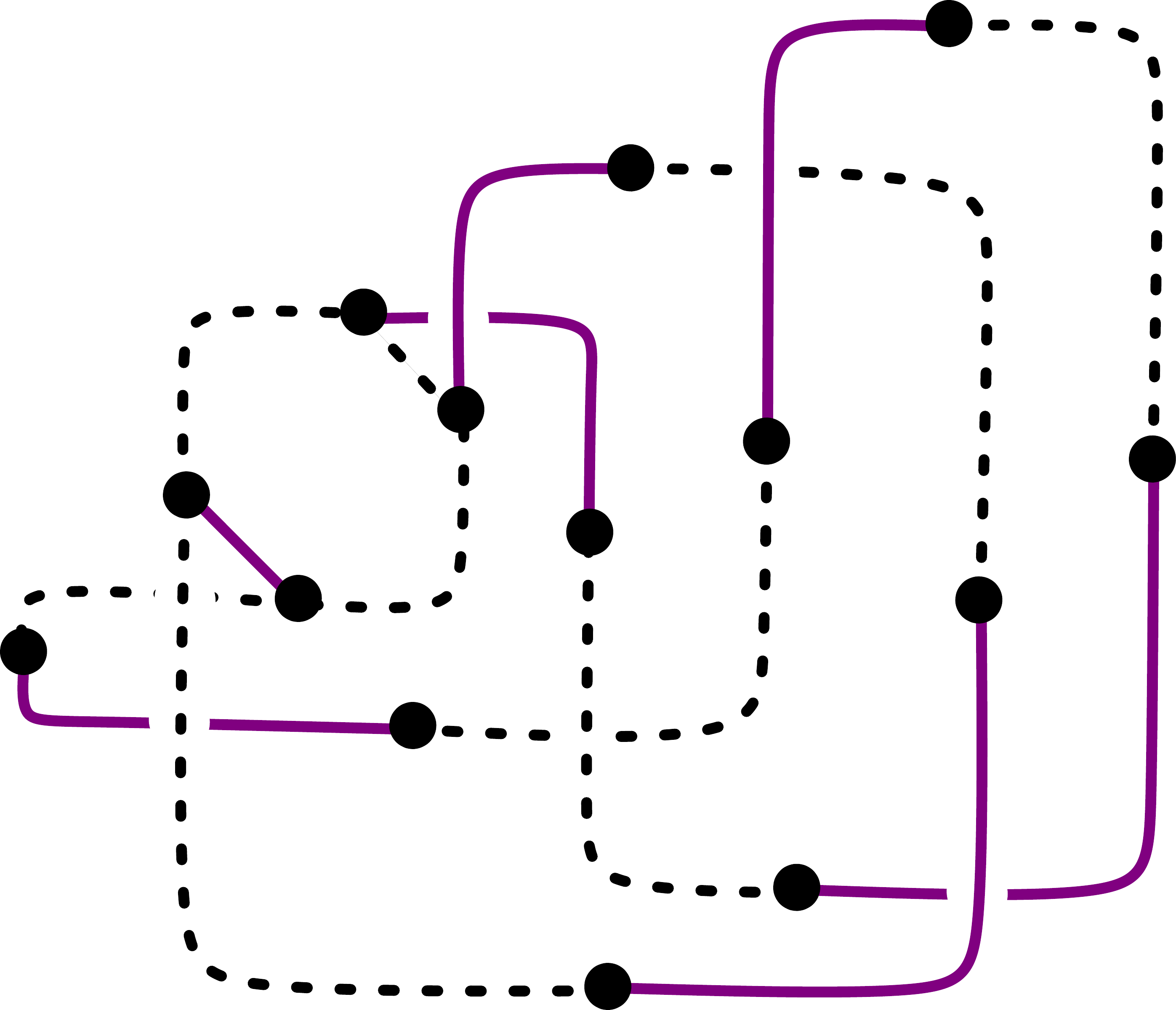}
\end{minipage}
\hfill
\begin{minipage}[t]{0.23\linewidth}
    \centering
    $\boxed{\mathfrak{m}_4}$\par\vspace{0.4em}
    \includegraphics[width=\linewidth]{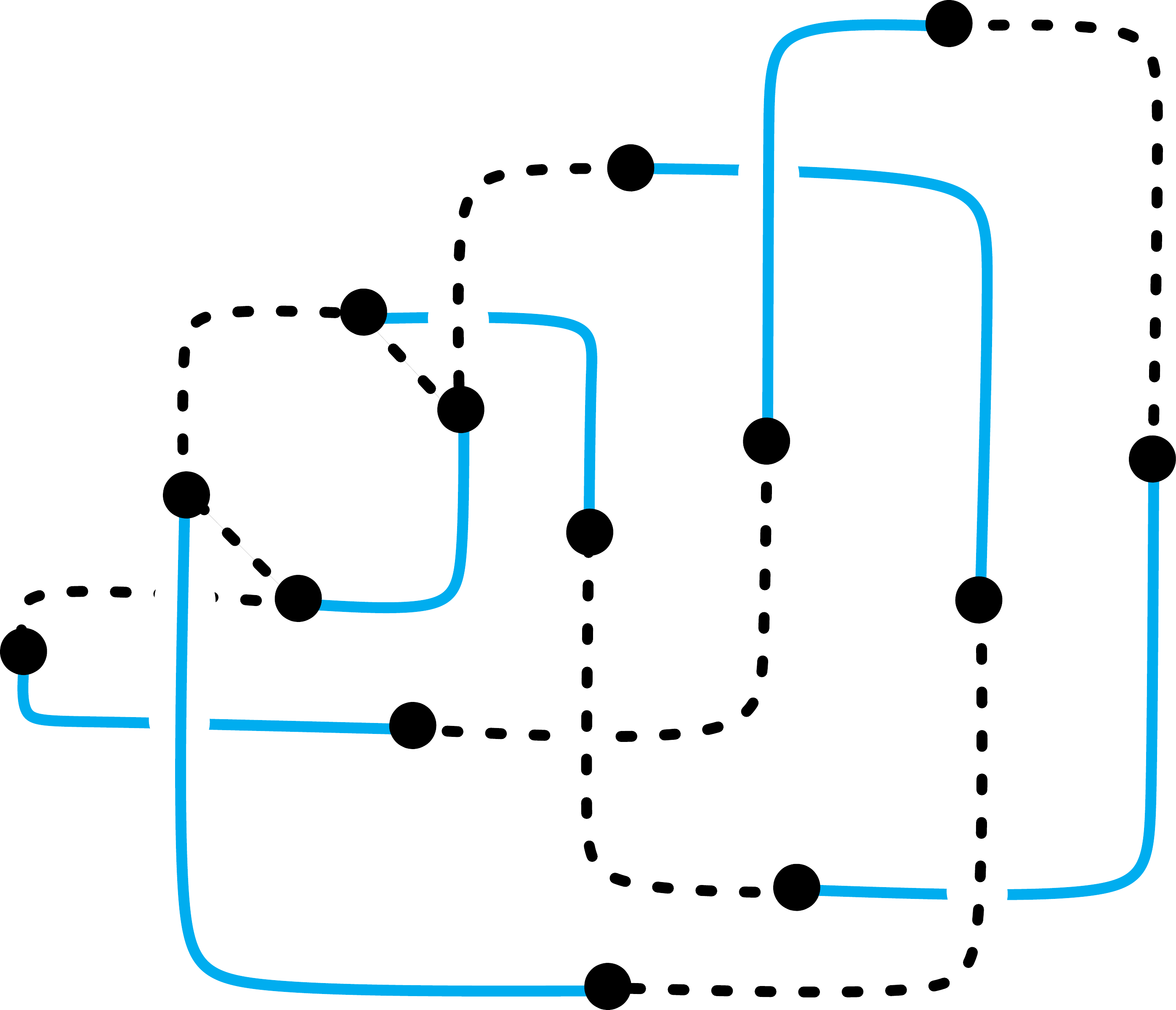}
\end{minipage}

    \vspace{1em}

    \includegraphics[width=0.23\linewidth]{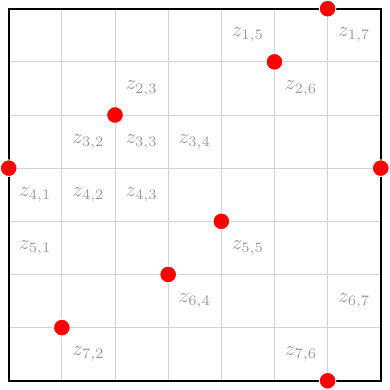}
    \hfill
    \includegraphics[width=0.23\linewidth]{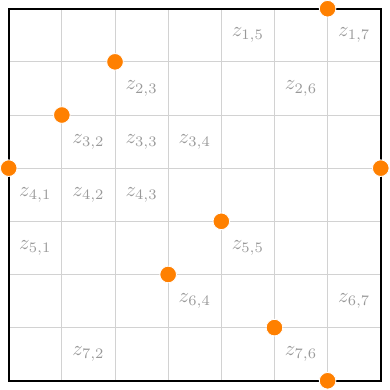}
    \hfill
    \includegraphics[width=0.23\linewidth]{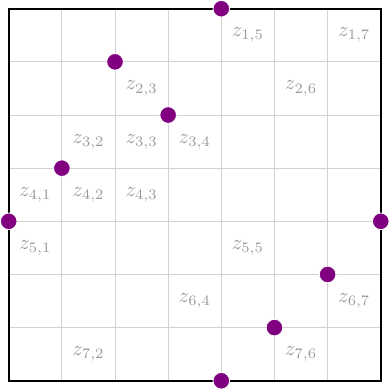}
    \includegraphics[width=0.23\linewidth]{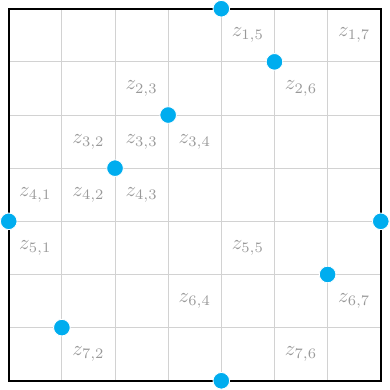}
    \hfill

    \par\vspace{1em}

\makebox[\linewidth]{%
    \hfill
    \begin{subfigure}[t]{0.32\linewidth}
        \centering
        \includegraphics[width=0.9\linewidth]{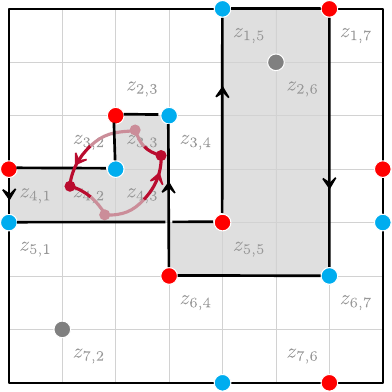}
        \caption{$D_{4,1} : \bm{x}_4 \to \bm{x}_1$}
        \label{fig:D41}
    \end{subfigure}
    \hfill
    \begin{subfigure}[t]{0.32\linewidth}
        \centering
        \includegraphics[width=0.9\linewidth]{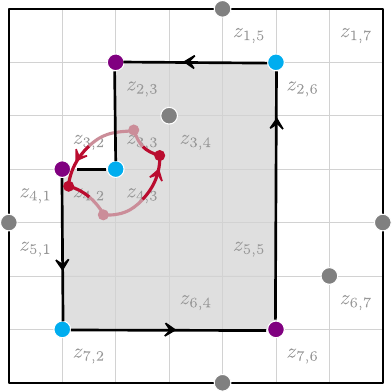}
        \caption{$D_{4,3}: \bm{x}_4 \to \bm{x}_3$}
        \label{fig:D43}
    \end{subfigure}
    \hfill
}

    \caption{
    Recall the stabilized $\Theta_4$-graph from Figure \ref{fig:ExampleGridDiagram}.
    The comparison with the local model seen in Figure \ref{fig:localSquareRegion} is $(z_1,z_2,z_3,z_4)\to (z_{3,3},z_{3,2},z_{4,2},z_{4,3})$.
    The top row displays perfect matchings $\mathfrak{m}_1$ through $\mathfrak{m}_4$ for $G$. The middle row displays the corresponding ``northwest'' Heegaard states in the grid diagram $\mathcal{G}$. 
    The final row shows domains associated to Seifert surfaces for $C_{4,1}^{NW}$ and $C_{4,3}^{NW}$ together with the oriented cycle $Q = C_{1,2} - C_{4,3} = C_{3,2}-C_{4,1}$ which bounds a square. 
    }
    \label{fig:MatchingsAndGenerators}
\end{figure}

By our previous work, the subdivision of the long edges is such that ${ G }$ admits $n$ perfect matchings; each of the edges ${ e_{i} }$ belonging to the boundary of the square belongs to a unique matching ${ \mathfrak{m}_{i} }$. It is not hard to see that for any ${ i,j \in \{1,2,3,4\} }$ such that $i\neq j$, the union ${ \mathfrak{m}_{i} \cup \mathfrak{m}_{j} }$ contains a single cycle 
$${ C(\mathfrak{m}_{i},\mathfrak{m}_{j}) = \mathfrak{m}_{i} \triangle \mathfrak{m}_{j} }$$ where ${ A \triangle B = (A \cup B) \setminus (A \cap B) }$ is the symmetric difference operation. In particular, letting $H \subset G$ denote the subdivided $K_4$ consisting of the union of the $e_i$ together with the long edges, one should observe that 
$$ C(\mathfrak{m}_{1}, \mathfrak{m}_{3}   ) = H \setminus (e_{2} \cup e_{4}  ) \qquad C(\mathfrak{m}_{2}, \mathfrak{m}_{4}  ) = H \setminus (e_{1} \cup e_{3} )  $$
are both cycles.
Each of these cycles contains both of the long edges of the graph. On the other hand, the remaining cycles each only contain one of the two long edges. 
We give our cycles orientations; write ${ C_{i,j} }$ for the cycle obtained by retaining the orientation of the ${ \mathfrak{m}_{i} }$ edges and reversing the orientation of the ${ \mathfrak{m}_{j} }$ edges.

The construction of a grid diagram ${ \mathcal{G} }$ for ${ G }$ requires a projection ${ \pi: (S^{3},G)  \rightarrow \mathbb{T}^{2}}$ where ${ \mathbb{T}^{2}\subset S^{3} }$ is the Heegaard torus with respect to which ${ G }$ is in bridge position. Recall that this projection is such that the arcs comprising ${ G }$ are projected to horizontal and vertical segments in ${ (\mathbb{T}^{2}, \bm{\alpha}, \bm{\beta}) }$ which are centered in each row and column. In particular, the intersections of ${ G }$ with the Heegaard torus coincide with the basepoints in the grid diagram.

We give $G$ a framing $f$ so that, letting $G^{NW}$ denote the framed push-off of $G$, ${ \pi(G) }$ and ${\pi(G^{NW}) }$ differ from one another by a small affine shift in the northwest direction in the grid diagram. 
In particular, the edges of $G$ which formerly projected to segments centered within a given row or column are shifted slightly northwest to lie on a segment of the $\alpha$ and $\beta$ curves. 

Each cycle $C \subset G$ has a corresponding push-off cycle $C^{NW}\subset G^{NW}$. Thinking of ${ \pi(C_{i,j}^{NW}) }$ as a knot diagram in the plane, one may apply Seifert's algorithm to obtain a projection ${ D_{i,j} }$ of a Seifert surface bounded by the oriented cycle ${ C_{i,j}^{NW} }$. We think of ${ D_{i,j} }$ as a domain which starts at ${ \bm{x}_i := \bm{x}^{NW\mathfrak{m}_{i}} }$ and ends at ${ \bm{x}_j := \bm{x}^{NW\mathfrak{m}_{j}} }$. The $D_{i,j}$ do not necessarily contribute to the differential of the grid complex, but they determine the gradings ${ gr_{i} }$ of the northwest generators via properties of the relative Maslov grading.

\begin{proposition}\label{prop:squareGradingComputation}
Let ${ G }$ be a stabilized ${ \Theta_{n} }$-graph that bounds a square and let ${ \mathcal{G} }$ be a grid diagram for ${ G }$ containing a region as in Lemma \ref{lem:plaidRegion}. Then the square grading for a Heegaard state $\bm{x} \in \mathbb{T}_{\alpha} \cap \mathbb{T}_{\beta}$ is valued in the set 
$${ \left\{-2\ell k(C_{1,4},C_{4,3}^{NW}), -2\ell k(C_{1,4},C_{4,3}^{NW})-2\right\}. }$$
\begin{proof}

Let $\mathcal{G}$ be a grid diagram for a graph obtained by some number of stabilizations of a $\Theta_n$-graph which bounds a square. 
By Lemma
\ref{lem:squareGradingConstrained}, it suffices to compute the square grading of the generator ${ \bm{x}_{1} }$. Applying the grading formula to the projected Seifert domains and some basic algebra yields 
$$
\begin{aligned}
gr_{\square}(\bm{x}_{1} )&=gr_{1}(\bm{x}_{1}  ) - gr_{2} (\bm{x}_{1} ) + gr_{3} (\bm{x}_{1} ) - gr_{4} (\bm{x}_{1} ) \\
&= gr_{1} (\bm{x}_{1} ) - gr_{2} (\bm{x}_{1} ) + gr_{3} (\bm{x}_{3} ) - gr_{4} (\bm{x}_{3} ) + \underbrace{\left( gr_{3} (\bm{x}_{1} ) - gr_{3}(\bm{x}_{3} )   \right)}_{= \mu(D_{1,3} ) - 2n_{3} (D_{1,3} )}  - \underbrace{\left( gr_{4} (\bm{x}_{1} ) - gr_{4}(\bm{x}_{3} )   \right)}_{= \mu(D_{1,3} ) - 2n_{4} (D_{1,3} )}    \\
&= \underbrace{gr_{1} (\bm{x}_{1} )}_{=0=gr_{2}(\bm{x}_{2} ) }  - gr_{2} (\bm{x}_{1} ) + \underbrace{gr_{3} (\bm{x}_{3} )}_{= 0 = gr_{4} (\bm{x}_{4} )}  - gr_{4} (\bm{x}_{3} ) + 2\Big(n_{4}(D_{1,3} ) - n_{3}(D_{1,3} )    \Big)     \\
%&=  gr_{2}(\bm{x}_{2} ) - gr_{2} (\bm{x}_{1} ) +gr_{2} (\bm{x}_{4} ) - gr_{2}(\bm{x}_{3} ) + \underbrace{\left( gr_{4}(\bm{x}_{4} ) - gr_{4} (\bm{x}_{3} )  \right)}_{= \mu(D_{4,3} ) - 2n_{4}(D_{4,3} ) }  - \underbrace{\left( gr_{2} (\bm{x}_{4} ) - gr_{2}(\bm{x}_{3} )  \right)}_{= \mu(D_{4,3} ) - 2n_{2}(D_{4,3} ) }  + 2\Big(n_{4}(D_{1,3} ) - n_{3}(D_{1,3} )    \Big)     \\
&= \underbrace{gr_{2} (\bm{x}_{2} )-gr_{2} (\bm{x}_{3} )}_{=\mu(D_{2,3} ) - 2n_{2} (D_{2,3} )}  -  \underbrace{\Big(gr_{2} (\bm{x}_{1} ) - gr_{2} (\bm{x}_{4} )\Big)}_{=\mu(D_{1,4} ) - 2n_{2} (D_{1,4} ) }  + 2\Big(n_{2} (D_{4,3}  ) - n_{4}(D_{4,3} )  \Big) + 2\Big( n_{4}(D_{1,3} ) - n_{3}(D_{1,3} )   \Big) \\
&= \underbrace{\mu(D_{2,3} ) - \mu(D_{1,4} )}_{A}  + 2\Big(\underbrace{n_{2} (D_{1,4} ) - n_{2} (D_{2,3} )}_{B}\Big)  + 2\Big(\underbrace{n_{2} (D_{4,3}  ) - n_{4}(D_{4,3} )  +  n_{4}(D_{1,3} ) - n_{3}(D_{1,3} )   }_{C} \Big). 
\end{aligned}
$$
\noindent
To determine the values $A,B,C$, we require some lemmas.
\begin{lemma}\label{lem:linkingFormula}
Let $a,b,c,d \in \{1,2,3,4\}$, $a \neq b$ and $c\neq d$. We have 
$$n_{a}(D_{c,d}) - n_b(D_{c,d}) = \ell k (C_{b,a}, C_{c,d}^{NW} ) = \langle C_{b,a}, PD[D_{c,d}] \rangle  .$$
\begin{proof}

This is an immediate consequence of the relative grading formulae, the fact that $D_{c,d}$ is the projection of a Seifert surface, and the fact that basepoints for a matching correspond to edges in the cycle $C_{b,a}$ which intersect the Seifert circles of $D_{c,d}$ in the multiplicities $n_a(D_{c,d})$ and $n_{b}(D_{c,d})$. Our orientation convention yields the claimed relationship with the linking number. 
\end{proof}
\end{lemma}

\begin{lemma}\label{lem:cycleIdentifications} 
Let ${ \mathcal{G} }$ be a grid diagram of standard square form. Let ${ g_{NW} }$ denote the grid cell containing the basepoint $z_2$. 
\begin{enumerate}[label=\roman*.]
\item\label{itm:isotopicCycles} ${ C_{4,1} }$ and ${ C_{3,2} }$ are isotopic as oriented knots in ${ S^{3} }$. Similarly, ${ C_{4,3}}$ and ${ C_{1,2} }$ are isotopic.
\item\label{itm:domainRelation} We have ${ D_{4,3}+g_{NW} = D_{1,2} }$ and ${ D_{4,1} + g_{NW} = D_{3,2} }$.
\end{enumerate}

\begin{proof}
For \ref{itm:isotopicCycles} one should observe that 
$${ C(\mathfrak{m}_{1},\mathfrak{m}_{4}) = \mathfrak{m}_{1} \triangle \mathfrak{m}_{4} = (\mathfrak{m}_{2} \triangle \mathfrak{m}_{3}) \triangle (\cup_{i=1}^{4} e_{i}   ) = C(\mathfrak{m}_{2},\mathfrak{m}_{3})\triangle (\cup_{i=1}^{4} e_{i}   ) }$$
and ${ \cup_{i=1}^{4}e_{i} }$ bounds a square in ${ S^{3} }$ which can be used to isotope the arc ${ e_{1} \cup e_{4} }$ to the arc ${ e_{2} \cup e_{3} }$. Our orientation conventions are such that this implies ${ C_{4,1} }$ and ${ C_{3,2} }$ are isotopic as oriented knots in ${ S^{3} }$. A very similar argument shows that ${ C_{4,3} }$ and ${ C_{1,2} }$ are isotopic. For item (\ref{itm:domainRelation}), observe that ${ \pi }$ maps the square bounded by ${ (\cup_{i=1}^{4} e_{i})^{NW} }$ to ${ g_{NW} }$. 
\end{proof}
\end{lemma}

Lemma \ref{lem:cycleIdentifications}.\ref{itm:domainRelation} immediately yields 
$${ B= n_{2}(D_{1,4}) - n_{2}(D_{2,3}) = n_{2}(D_{1,4}) - n_{2}(D_{1,4}-g_{NW}) = 1 .}$$
For ${ A }$, recall that the Maslov index of a domain ${ D_{i,j} }$ is given by the formula (see e.g. \S~4 of \cite{Lip06}) 

$$ \mu(D_{i,j} ) = e(D_{i,j}) + n_{\bm{x}_{i} }(D_{i,j}) + n_{\bm{x}_{j} }(D_{i,j})   . $$
Lemma \ref{lem:cycleIdentifications}.\ref{itm:domainRelation} together with linearity of the Euler measure gives
\begin{equation}\label{eqn:Astep1}
\mu(D_{2,3} ) - \mu(D_{1,4} ) = n_{\bm{x}_{2} }(D_{2,3} ) - n_{\bm{x}_{1} }(D_{1,4} ) + n_{\bm{x}_{3} }(D_{2,3} ) - n_{\bm{x}_{4} }(D_{1,4} ) .
\end{equation}
Actually more is true: there is an explicit bijection 
$$ \text{corners of }D_{1,4} \text{ contributed by } \left\{\begin{matrix} \bm{x}_{1} \\ \bm{x}_{4}   \end{matrix}\right\} \longleftrightarrow \text{corners of }D_{2,3} \text{ contributed by } \left\{\begin{matrix} \bm{x}_{2} \\ \bm{x}_{3}  \end{matrix}\right\}  $$ 
which respects the status of each corner being acute or obtuse as well as the local multiplicity of the domain about said corner. 
The relevant consequences are:
\begin{enumerate}[label=(\alph*)]
\item\label{itm:a} All local multiplicities in ${ n_{\bm{x}_{2}}(D_{2,3}) }$ and  ${n_{\bm{x}_{1}}(D_{1,4}) }$ coming from corners of their respective domains cancel in their difference. The surviving terms correspond to grid intersection points contributing to $\bm{x}_2$ and $\bm{x}_1$ which lie in the interior of $D_{2,3}$ and $D_{1,4}$ respectively. 
\item\label{itm:b} Recall that every grid intersection point contributing to $\bm{x}_1$ and $\bm{x}_2$ is northwest of a basepoint associated to our given  matching. In particular, this means a grid intersection point $x \in \bm{x}_i$ only lies in the interior of a domain $D$ if the basepoint $z \in \mathfrak{m}_i$ which $x$ is northwest of also lies in the interior of $D$.
\end{enumerate}

\noindent
For brevity, we adopt the notation 
$$\Lambda(C_{a,b},C_{c,d}) = \ell k (C_{a,b},C_{c,d}^{NW}) = \langle C_{a,b}, PD[D_{c,d}] \rangle \qquad \text{ and } \qquad Q = C_{1,2}-C_{4,3} = C_{3,2}-C_{4,1}. $$
One should observe that $Q$ is exactly the oriented cycle which bounds a square; see the bottom of Figure \ref{fig:MatchingsAndGenerators}.
We are now ready to compute the values of $A$ and $C$.

$$
\begin{aligned}
A &= \mu(D_{2,3} ) - \mu(D_{1,4} ) & \\
  &= n_{\bm{x}_{2} }(D_{2,3} ) - n_{\bm{x }_{1} }(D_{1,4} ) + n_{\bm{x }_{3} }(D_{2,3} ) - n_{\bm{x}_{4} }(D_{1,4} ) & \text{Equation \ref{eqn:Astep1}} \\ 
  &= n_2(D_{2,3}) - n_4(D_{1,4}) + n_3(D_{2,3}) - n_1(D_{1,4}) & \text{Items \ref{itm:a} \& \ref{itm:b} } \\ 
  &= n_2(D_{2,3}) - n_4(D_{2,3}) + n_3(D_{2,3}) - n_1(D_{2,3}) & \text{Lemma \ref{lem:cycleIdentifications}.\ref{itm:domainRelation} } \\
  %= \langle C_{4,2} + C_{1,3}, PD[D_{2,3}] \rangle & \text{Lemma \ref{lem:linkingFormula}} \\
  &= \Lambda(C_{4,2}, C_{2,3} ) + \Lambda(C_{1,3}, C_{2,3}) & \text{Lemma } \ref{lem:linkingFormula} \\
  &= \Lambda (2C_{4,3}+Q, -C_{4,1}-Q) & \left\{\begin{aligned} C_{4,2}+C_{1,3} &= 2C_{4,3}+Q \\ C_{2,3} &= -C_{4,1}-Q \end{aligned}\right. \\
  &= -2 \Lambda(C_{4,3},C_{4,1}) - 2 \Lambda(C_{4,3},Q) - \Lambda(Q,C_{4,1}) - \Lambda(Q,Q) .
\end{aligned}
$$

For ${ C }$, note that Lemma \ref{lem:linkingFormula} tells us 
$$ 
\begin{aligned}
C &= n_{2}(D_{4,3}) - n_{4}(D_{4,3}) + n_{4}(D_{1,3}) - n_{3}(D_{1,3}) \\ 
&= \ell k (C_{4,2}, C_{4,3}^{NW}) + \ell k (C_{3,4}, C_{1,3}^{NW}) \\
&= \ell k (C_{4,2}, C_{4,3}^{NW}) - \ell k (C_{4,3}, C_{1,3}^{NW}) \\
&=\Lambda(C_{4,1},C_{4,3}) + \Lambda(Q,C_{4,3}) + \Lambda(C_{4,3},C_{4,1}) & \left\{\begin{aligned} C_{4,2} &= C_{4,3} + C_{4,1} + Q \\ C_{1,3} &= C_{4,3}-C_{4,1} \end{aligned} \right. . \\
\end{aligned}
$$ 
Finally we have
$$
gr_{\square}(\bm{x}_1) = A + 2B + 2C = 2\Lambda(C_{4,1},C_{4,3}) +2 + 2\Lambda(Q,C_{4,3}) - 2\Lambda(C_{4,3},Q) - \Lambda(Q,C_{4,1}) - \Lambda(Q,Q).
$$
Observe we always have $\Lambda(Q,Q)=1$ and
$\Lambda(C_{4,3},Q) = \ell k (C_{4,3}, Q^{NW}) = 0$. These hold because $Q^{NW} = \partial g_{NW}$ and $C_{4,3}$ does not intersect the domain $g_{NW}$. On the other hand, $\ell k (Q,C_{4,3}^{NW})=\ell k(Q,C_{4,1}^{NW})=-1$; see Figures \ref{fig:D41} and \ref{fig:D43}. 
Since $-C_{4,1} = C_{1,4}$, we conclude 
$$gr_{\square}(\bm{x}_1) = 2\ell k (C_{4,1},C_{4,3}^{NW}) = -2\ell k (C_{1,4}, C_{4,3}^{NW}). $$
Lemma \ref{lem:squareGradingConstrained} therefore tells us that for any $\bm{x} \in \mathbb{T}_{\alpha} \cap \mathbb{T}_{\beta}$ we have
$$gr_{\square}(\bm{x}) \in \{-2\ell k (C_{1,4},C_{4,3}^{NW}),-2\ell k (C_{1,4},C_{4,3}^{NW})-2\}.$$

\end{proof}
\end{proposition}

\begin{corollary}\label{cor: squareBoundSumUp}
    Suppose $\Gamma$ is a spatial $\Theta_n$-graph in which some four-circuit bounds a square which we identify with a nonorientable band $b$. Then $\mathrm{CFT}(\Gamma)$ admits a model whose basis generators have square gradings in $\{\frac{1}{2}e(b) \pm 1\}.$
\end{corollary}
\begin{proof}
    Perform two vertex stabilizations on $\Gamma$ to make the square embedded (separating the vertices at its corners), and then a series of edge stabilizations to make the resulting graph admit a grid diagram. Call this graph $\Gamma'$; Corollary~\ref{cor:StabilizationHtpy} gives us a coloring $\mathfrak{c}$ for which $\mathrm{CFT}(\Gamma',\mathfrak{c})$ is homotopy equivalent to $\mathrm{CFT}(\Gamma)$. Proposition~\ref{prop:squareGradingComputation} then gives us a chain model for $\mathrm{CFT}(\Gamma')$ in which all generators are in two specific degrees; specializing certain variables gives a model for $\mathrm{CFT}(\Gamma',\mathfrak{c})$ with the same property, and so one for $\mathrm{CFT}(\Gamma)$ using the homotopy equivalence. 

    It remains only to identify the two gradings in the statement of Proposition~\ref{prop:squareGradingComputation} with $\frac{1}{2}e(b) \pm 1$. When computing $-2\ell k(C_{1,4},C_{4,3}^{NW})$, one adds up a sign $\pm 1$ for each crossing between $C_{1,4}$ and $C_{4,3}^{NW}$ in the grid diagram, counting all of the crossings between the long diagonal edges as well as one distinguished positive crossing in the square. To find $e(b)$, one tracks which crossings in the cycle formed by the two long diagonal edges change sign between the two resolutions within the square and counts $\pm 2$ for each; this counts precisely the crossings between the two different long edges. Overall, then, $-2\ell k(C_{1,4},C_{4,3}^{NW}) = \frac12 e(b) + 1$, and so $-2\ell k(C_{1,4},C_{4,3}^{NW}) - 2 = \frac12 e(b) - 1$. 
\end{proof}
Note that one has to interpret this statement somewhat carefully, since the variables assigned to edges may have nonzero square grading --- the correct interpretation, and the one needed to e.g. deduce bounds on $\Upsilon$, is that $\mathrm{CFT}(\Gamma)$ is homotopy equivalent to a complex of free $R_n$-modules with every generator in square grading $\frac{1}{2}e(b) \pm 1$.

\bibliographystyle{alpha}
\bibliography{bibliography}

\end{document}